\documentclass[a4paper,12pt, reqno]{amsart} 
\usepackage{amssymb,amsthm,amsmath}

\usepackage{amsfonts}
\usepackage{amscd}
\usepackage{amssymb}
\usepackage{enumerate}
\usepackage{graphicx}
\allowdisplaybreaks
\usepackage{mathabx}
\usepackage{color}
\usepackage{amsbsy}
\usepackage{graphicx}
\usepackage{amsthm}
\usepackage{amsmath}
\usepackage{amsxtra}
\usepackage{mathrsfs}
\usepackage{bbm}
\usepackage{dsfont}

\usepackage{ifthen}

\usepackage[colorlinks,citecolor=red,pagebackref,hypertexnames=false]{hyperref} 

\newtheorem{thm}{Theorem}[section]
\newtheorem{Cor}[thm]{Corollary}
\newtheorem{lem}[thm]{Lemma}

\newtheorem{prop}[thm]{Proposition}

\theoremstyle{definition}
\newtheorem{defn}[thm]{Definition}
\theoremstyle{remark}
\newtheorem{rmk}[thm]{Remark}
\numberwithin{equation}{section}
\definecolor{red}{rgb}{1.0, 0.0, 0.0}
\newcommand{\Bea}{\begin{eqnarray*}}
	\newcommand{\Eea}{\end{eqnarray*}}
\newcommand{\Be} {\begin{equation*}}
	\newcommand{\Ee} {\end{equation*}}
\newcommand{\be} {\begin{equation}}
	\newcommand{\ee} {\end{equation}}
\newcommand{\bea} {\begin{eqnarray}}
	\newcommand{\eea} {\end{eqnarray}}

\title[]
{Non-harmonic $M$-elliptic pseudo differential operators on manifolds 
}
\author[Aparajita Dasgupta]{Aparajita Dasgupta}
\address{
	Aparajita Dasgupta:
	\endgraf
	Department of Mathematics
	\endgraf
	Indian Institute of Technology, Delhi, Hauz Khas
	\endgraf
	New Delhi-110016 
	\endgraf
	India
	\endgraf
	{\it E-mail address:} {\rm adasgupta@maths.iitd.ac.in}
}
\author{Vishvesh Kumar} 

\address{Vishvesh Kumar, Ph. D.  \endgraf Department of Mathematics: Analysis, Logic and Discrete Mathematics
	\endgraf Ghent University
	\endgraf Krijgslaan 281, Building S8,	B 9000 Ghent,
	Belgium .} 
\email{vishveshmishra@gmail.com}

\author[Lalit Mohan]{Lalit Mohan}
\address{
	Lalit Mohan:
	\endgraf
	Department of Mathematics
	\endgraf
	Indian Institute of Technology, Delhi
	\endgraf
	India
	\endgraf
	{\it E-mail address:} {\rm mohanlalit871@gmail.com}
}
\author[Shyam Swarup Mondal]{Shyam Swarup Mondal}
\address{
	Shyam Swarup Mondal:
	\endgraf
	Department of Mathematics
	\endgraf
	Indian Institute of Technology, Delhi
	\endgraf
	India
	\endgraf
	{\it E-mail address:} {\rm mondalshyam055@gmail.com}
}
\date{\today}

\begin{document}
	
	\maketitle

	\allowdisplaybreaks
	
	\begin{abstract} 
		In this article, we introduce and study $M$-elliptic pseudo-differential operators in the framework of non-harmonic analysis of boundary value problems on a manifold $\Omega$ with boundary $\partial \Omega$,   introduced by Ruzhansky and Tokmagambetov ( Int. Math. Res. Not. IMRN, (12), 3548-3615, 2016) in terms of a model operator $\mathfrak{L}$. More precisely, we consider a weighted $\mathfrak{L}$-symbol class   $M_{\rho, 0, \Lambda}^{m}, m\in \mathbb{R},$ associated to a suitable weight function  $\Lambda$ on a countable set $\mathcal{I} $ and study elements of the symbolic calculus for pseudo-differential operators associated with $\mathfrak{L}$-symbol class   $M_{\rho, 0, \Lambda}^{m},$  by deriving formulae for the composition, adjoint, and transpose. 
		Using the notion of $M$-ellipticity for symbols belonging to $\mathfrak{L}$-symbol class   $M_{\rho, 0, \Lambda}^{m}$, we  construct the parametrix of $M$-elliptic  pseudo-differential operators.  	Further, we investigate the minimal and maximal extensions for $M$-elliptic pseudo-differential operators and show that they coincide when the symbol $\sigma\in M_{\rho, 0, \Lambda}^{m}, $ is $M$-elliptic. 	
		We provide a necessary and sufficient condition to ensure that the pseudo-differential operators $T_{\sigma}$  with symbol in the  $\mathfrak{L}$-symbol class $M_{\rho, 0,\Lambda}^{0} $  is a  compact operator in        $L^{2}(\Omega)$   or a Riesz operator in $L^{p}(\Omega).$ 
		Finally,  we prove G\"arding's inequality for pseudo-differential operators associated with symbol from    $M_{\rho, 0,\Lambda}^{0} $ in the setting of non-harmonic analysis.

	\end{abstract}
	\tableofcontents 	
	
	\section{Introduction} \label{Introduction}
The theory of pseudo-differential operators holds a crucial position in modern mathematics as it draws significant inspiration from partial differential equations, signal processing, and time-frequency analysis, as discussed in \cite{Hor, MR&VT book, fis1}. Pseudo-differential operators, which act on functions defined on smooth manifolds, provide a fundamental generalization of differential operators. The study of pseudo-differential operators originated in the 1960s with the pioneering works of Kohn and Nirenberg \cite{Niren} and H\"ormander \cite{Hor}. Their focus was on singular integral differential operators, particularly in the context of inverting differential operators to solve elliptic differential equations. Since then, this theory has become an indispensable tool, with notable connections to mathematical physics, harmonic analysis, quantum field theory, and index theory.

	Consider an n-dimensional smooth manifold $\Omega$ with boundary $\partial \Omega$ such that  $\bar{\Omega}=\Omega \cap \partial \Omega$ represents a compact manifold (not necessarily smooth at the boundary). The main focus of this work is to explore $M$-elliptic pseudo-differential operators within the framework of non-harmonic analysis, specifically in the context of boundary value problems on $\Omega$ with $\partial \Omega$ as the boundary.

	The study of symbols in the context of boundary value problems on smooth manifolds with boundaries is an active and fruitful research area, posing numerous open questions and concerns. The symbol of a pseudo-differential operator plays a crucial role in understanding the operator's spectrum, particularly in the context of boundary value problems on smooth manifolds with boundaries. It helps analyze the behavior of solutions, estimate their regularity and decay properties near the boundary, and determine the expected decay rate.

Various authors have delved into boundary value problems for pseudo-differential operators on manifolds or within domains such as $\Omega\subseteq \mathbb{R}^n$ with boundaries \cite{Ruz16,MR&JPVR,ruzhansky2017nonharmonic}. For a global perspective, references \cite{DMT17,KT15} provide insightful analyses of pseudo-differential operators on manifolds. The investigation of boundary value problems in non-harmonic analysis often revolves around the use of pseudo-differential operators. In this work, our primary objective is to delve into the pseudo-differential calculus associated with a boundary-value problem $\mathfrak{L}_\Omega$ defined by a differential operator $\mathfrak{L}$ on $\Omega$ with a discrete spectrum and appropriate boundary conditions, all within the framework of non-harmonic analysis.

It is noteworthy that the non-harmonic theory presented in \cite{Ruz16,ruzhansky2017nonharmonic} has fascinating applications in partial differential equations, particularly in the context of the wave equation, as demonstrated in \cite{ruzhansky2017very, ruzhansky2017wave}. The non-harmonic analysis of boundary value problems bears resemblance to global pseudo-differential analysis on closed manifolds, as outlined in \cite{delgado1,delgado2}, which builds upon the non-harmonic analysis of compact manifolds by Seeley \cite{seeley1,seeley2}. This type of analysis proves effective in various problems, such as characterizing classes of functions and distributions of the Komatsu type \cite{apara}. Additionally, the analysis serves as a mathematical framework for studying linear operators in the context of non-self-adjoint quantum mechanics \cite{mostafazadeh2010pseudo,inoue2016non}.
	
	The unit circle $\mathbb{T}$ serves as a simple and illustrative example of global analysis. By considering the system of eigenfunctions $\left\{e^{i x \cdot k}\right\}_{k \in \mathbb{Z}}$ associated with the differentiation operator $-i \frac{d}{d x}$, the Fourier series on the unit circle $\mathbb{T}$ can be interpreted as a unitary transform in the Hilbert space $L^2(-\pi, \pi)$. This transform is generated by the differentiation operator $-i \frac{d}{d x}$ subject to periodic boundary conditions.

This concept can be extended to a broader framework, even when the problem lacks symmetries. The generalization, as presented in \cite{Ruz16,ruzhansky2017nonharmonic}, involves considering a differential operator $\mathfrak{L}$ of order $m$ with smooth coefficients, rather than the specific differential operator $-i \frac{d}{d x}$. Additionally, the authors consider that $\Omega$ is equipped with certain boundary conditions, resulting in a discrete spectrum. The family of eigenfunctions corresponding to this spectrum forms a Riesz basis in $L^2(\Omega)$.
	
	The presence of such a basis allows for the emulation of harmonic analysis constructions and facilitates the realization of a global analysis framework akin to the case of the unit circle. The term ``non-harmonic analysis" originates from the work of Paley and Wiener \cite{paley1934fourier}, who studied exponential systems $\left\{e^{2 \pi i x \cdot \lambda}\right\}_{\lambda \in \Lambda}$ on $L^2(0,1)$ for a discrete set $\Lambda$. Paley and Wiener employed the term ``non-harmonic Fourier series" to emphasize the distinction from the conventional (harmonic) Fourier series that arises when $\Lambda=\mathbb{Z}$. For further elucidation and advancements in non-harmonic analysis, we suggest referring to survey papers by Sedletskii \cite{saki1}. Additionally, in \cite{Ruz16}, numerous examples of the model operator $\mathfrak{L}$ and corresponding boundary conditions (BC) are presented.

	The Mikhlin-H\"ormander theorem for Fourier multipliers states that, in general, pseudo-differential operators associated with the class $S_{1, 0}^{0}$ are $L^{p}$-bounded. Similarly, it is well known that pseudo-differential operators with symbols from the class $S_{\rho, \delta}^{0}$, where $0 \leq \delta<\rho \leq 1$, are $L^{p}$-bounded when $p=2$ \cite{beals, Hor}. However, the situation changes significantly when $p\neq2$. Fefferman \cite{feff} observed that pseudo-differential operators with symbols in $S_{\rho, 0}^{0}$, where $0<\rho<1$, are not generally $L^{p}$-bounded for $p \neq 2$. To address this issue, Taylor introduced a suitable symbol subclass $M_{\rho}^{m}$ of $S_{\rho, 0}^{0}$ and developed a symbolic calculus for the associated pseudo-differential operators in \cite{taylor81}. Subsequently, Garello and Morando introduced the subclass $M_{\rho, \Lambda}^{m}$ of $S_{\rho, \Lambda}^{m}$, which are weighted versions of the symbol class introduced by Taylor, and have found various applications in studying the regularity of multi-quasi-elliptic operators \cite{GM, mor05}. Over the past sixteen years, numerous researchers have investigated various aspects of pseudo-differential operators associated with the symbol class $M_{\rho, \Lambda}^{m}$ on $\mathbb{R}^n$. Detailed studies on $M$-elliptic pseudo-differential operators in various frameworks can be found in \cite{kal, GM, MWMelliptic, AM, MWspectral, AD&LM, SSM&VK, MWMelliptic}.

Motivated by previous research, we consider a weighted symbol class $\mathfrak{L}$-symbol class $M_{\rho, 0, \Lambda}^{m}$ associated with a suitable weight function $\Lambda$ defined on a countable set $\mathcal{I}$ and investigate the symbolic calculus for pseudo-differential operators associated with the $\mathfrak{L}$-symbol class $M_{\rho, 0, \Lambda}^{m}$. We derive formulae for composition, adjoint, and transpose operations and further construct the parametrix of $M$-elliptic pseudo-differential operators using the concept of $M$-ellipticity for symbols belonging to the $\mathfrak{L}$-symbol class $M_{\rho, 0, \Lambda}^{m}$.

	The study of pseudo-differential operators involves investigating the conditions under which the minimal and maximal extensions of these operators on $L^p(\mathbb{R}^n)$ spaces, where $1 < p < \infty$, coincide. Wong demonstrated in \cite{MWMelliptic} that for $M$-elliptic operators, the minimal and maximal extensions on $L^p(\mathbb{R}^n)$ coincide. Wong also determined the domains of the minimal and maximal operators. Several authors have studied the properties of the minimal and maximal extensions, as well as Fredholmness and symbolic calculus of $M$-elliptic pseudo-differential operators in different settings \cite{vish, MR&JPVR, mor16, AM, kal}. In this work, we investigate the minimal and maximal extensions for $M$-elliptic pseudo-differential operators in the framework of non-harmonic analysis and prove that they coincide when the symbol belongs to the class $M_{\rho, 0, \Lambda}^{m}$ and is $M$-elliptic.

Another important result in the analysis is the Gohberg lemma, which has applications in spectral theory and singular integral equations. The Gohberg lemma was first derived by Gohberg in \cite{Goh} while investigating integral operators. In 1970, Gru\v{s}hin obtained the Gohberg lemma for pseudo-differential operators with bounded symbols \cite{VVG}. Later, an analog of Gohberg's lemma was proved for symbols in the H\"ormander class $S^0_{1, 0}(\mathbb{ S}^1\times \mathbb{ Z })$ in \cite{SM&MW} to study spectral invariance (see also \cite{dasgupta1}). This theory was further extended to compact Lie groups by the first author and Ruzhansky in \cite{dasgupta}. Recently, Ruzhansky and Velasquez-Rodriguez investigated a non-harmonic version of Gohberg's lemma for symbols in the H\"ormander class $S^0_{1, 0}(\Omega \times \mathcal{I})$ \cite{MR&JPVR}. In this paper, we establish a version of the Gohberg lemma for the weighted symbol class $M_{\rho, \Lambda}^0(\bar{\Omega} \times \mathcal{I})$. As an application of the Gohberg lemma, we provide a sufficient and necessary condition to ensure that the corresponding pseudo-differential operator is compact on $L^2(\bar{\Omega})$.

	The G\r{a}rding inequality is a significant result in the theory of pseudo-differential operators, particularly in the study of problems related to parabolic-type initial value problems. It provides valuable information about the regularity, existence, and uniqueness of solutions to such problems. The investigation of the G\r{a}rding inequality for strongly elliptic operators was initially conducted by G\r{a}rding \cite{LG}. His work established the existence of solutions to the Dirichlet problem for elliptic operators and examined the distribution of eigenvalues. Since then, many researchers have focused on studying the G\r{a}rding inequality for pseudo-differential operators in various contexts. In particular, this was further improved and generalised by several researchers including   Agmon \cite{Agmon}, Smith \cite{Smith} and Schechter \cite{Sche, Sche2}.

The G\r{a}rding inequality has been explored for H\"{o}rmander symbols on $\mathbb{R}^n$ with $0 \leq \delta \leq \rho \leq 1$, for SG-class operators, for compact Lie groups with matrix-valued symbols, and in the context of non-harmonic analysis on general smooth manifolds with applications to PDE \cite{taylor81,MR&JW,DC&MR}. In this paper, we investigate the G\r{a}rding inequality for $M$-elliptic pseudo-differential operators associated with symbols from the $M_{\rho, 0,\Lambda}^{0}(\bar{\Omega} \times \mathcal{I})$ class in the framework of non-harmonic analysis. Our goal is to establish the inequality and analyze its implications for the properties of the operators and the solutions to associated problems.

The structure of the manuscript is organized into six sections, including the introduction:
 \begin{itemize}
     \item {\bf Section \ref{Introduction}. Introduction}: This section provides an overview of the motivation and objectives of the study, as well as a brief introduction to the concepts and background of the topic.

    \item  {\bf Section \ref{sec2}. Basics of non-harmonic analysis}: In this section, the fundamentals of harmonic analysis in the framework of non-harmonic analysis are recalled. The authors review the concepts and techniques introduced by M. Ruzhansky and N. Tokmagambetov \cite{Ruz16} in their work, focusing on the model operator $\mathfrak{L}$.

   \item  {\bf Section \ref{sec3}. Weighted $\mathfrak{L}$-symbol classes and $M$-elliptic pseudo-differential operators}: This section delves into the investigation of weighted  $\mathfrak{L}$-symbol class $M_{\rho, 0, \Lambda}^{m}(\bar{\Omega} \times \mathcal{I})$. We study various important properties of this class, including symbolic calculus for pseudo-differential operators associated with $M_{\rho, 0, \Lambda}^{m}$. We derive formulae for composition, adjoint, and transpose operations and construct the parametrix of $M$-elliptic pseudo-differential operators.

    \item {\bf Section \ref{sec4}. Minimal and maximal extension operators}: Here, we focus on the construction of the minimal and maximal extensions of $M$-elliptic pseudo-differential operators. we prove the coincidence of these extensions on $L^{2}(\Omega)$ under the assumption of $M$-ellipticity for the symbol and determine the domains of the minimal operators.

   \item  {\bf Section \ref{sec5}. Gohberg's lemma, Compact operators and Riesz operators}: This section introduces a version of Gohberg's lemma for the weighted symbol class $M_{\rho, 0, \Lambda}^0(\bar{\Omega} \times \mathcal{I})$. We establish the lemma and provide a necessary and sufficient condition for the pseudo-differential operators with symbols from this class to be compact operators in $L^{2}(\Omega)$ or Riesz operators in $L^{p}(\Omega)$, where $1 < p < \infty$.

    \item {\bf Section \ref{sec6}. The G\r{a}rding inequality}: In the final section, the G\r{a}rding inequality for $\mathfrak{L}$-elliptic pseudo-differential operators with symbols from $M_{\rho, 0, \Lambda}^0(\bar{\Omega} \times \mathcal{I})$ is investigated. We  discuss the inequality and its implications, providing sufficient conditions for the existence and uniqueness of strong solutions in $L^{2}(\overline{\Omega})$ for the pseudo-differential operator $T_{\sigma}$.
    \end{itemize}

Overall, the manuscript covers a comprehensive analysis of various aspects of pseudo-differential operators associated with the weighted symbol class $M_{\rho, 0, \Lambda}^m$ in the framework of non-harmonic analysis. We explore symbolic calculus, minimal and maximal extensions, Gohberg's lemma, and G\r{a}rding's inequality, providing insights into the properties and applications of these operators in different contexts.

	\section{Basics of Non-Harmonic Analysis}\label{sec2}
	In this section, we recall   basics of harmonic analysis  in the framework of   non-harmonic analysis of boundary value problems on a manifold $\Omega$ with boundary $\partial \Omega$  in terms of a model operator $\mathfrak{L}$ to make the paper self contained. A complete account of analysis  in the framework of  non-harmonic analysis of boundary value problems on a manifold can be found in \cite{DMT17,Ruz16,ruzhansky2017nonharmonic,MR&JPVR}. However, we mainly adopt the notation and terminology given in \cite{MR&JPVR}.
	
	We will be working with discrete sets of eigenvalues and eigenfunctions indexed by a countable set $\mathcal{I} $. In different problems it may be more convenient to make different choices for this set, e.g. $\mathcal{I}=\mathbb{N}$ or $\mathbb{Z}$ or $\mathbb{Z}^{k},$ etc. In order to allow different applications we will be denoting it by $\mathcal{I}$ and without loss of generality one can  assume $ \mathcal{I}$ is a subset of $\mathbb{Z}^{k}$ for some $k \geq 1$. For simplicity, one may  think of $\mathcal{I}=\mathbb{Z}$ or $\mathcal{I}=\mathbb{N} \cup\{0\}$ throughout this article. Let us consider the assumptions throughout the paper.
	
	\textbf{Assumption 1}: Let $\Omega$ be a smooth n-dimensional manifold with boundary $\partial \Omega$
	such that $\bar{\Omega}$ is a compact (not necessarily smooth in the boundary) manifold. Consider the  differential operator $\mathfrak{L}$ of order $m$ with smooth bounded coefficients in $\Omega$, equipped with some	fixed linear boundary conditions. Assume that the boundary conditions  lead to a discrete spectrum of $\mathfrak{L}$ with a family of eigenfunctions yielding a Riesz basis in $L^2(\Omega)$. The discrete sets of eigenvalues and eigenfunctions will be indexed by a countable set $\mathcal{I}$ and	without loss of generality, will be always a subset of $\mathbb{Z}^k$ for some $k\in \mathbb{N}$. We consider the	spectrum $\{\lambda_\xi\in \mathbb{C}: \xi \in \mathcal{I}\}$ of $\mathfrak{L}$ with corresponding eigenfunctions in $L^2(\Omega)$ denoted by $u_\xi$, i.e., $$\mathfrak{L}u_\xi=\lambda_\xi u_\xi~\text{ in $\Omega$, for all $\xi \in \mathcal{I},$}~$$
	and the eigenfunctions $u_\xi$ satisfy the boundary conditions (BC). The conjugate spectral problem is
	$$\mathfrak{L}^*v_\xi=\overline{\lambda_\xi} v_\xi~\text{ in $\Omega$, for all $\xi \in \mathcal{I}$},$$
	which we equip with the conjugate boundary conditions which is denoted by (BC)$^*$.

	Here we auusme that the functions $u_\xi, v_\xi$ are normalised, i.e. $\|u_\xi\|_{L^2(\Omega)}=1$ and $ \|v_\xi\|_{L^2(\Omega)}=1$ for all $\xi \in \mathcal{I}$ and there exist a constant $C>0$ and $\mu_0>0$ such that 
	$$\sup_{x\in \Omega} |u_\xi(x)|\leq C\langle \xi \rangle^{\mu_0}$$ for every $\xi \in \mathcal{I}.$
	Here we  define the weight $$\langle \xi \rangle:=(1+|\lambda_\xi|^2)^{\frac{1}{2m}},$$
	where $m$ is the order of the differential operator $\mathfrak{L}.$ 
	The  systems $ \{u_\xi\}_{\xi \in \mathcal{I}}$ and $ \{v_\xi\}_{\xi \in \mathcal{I}}$ are  biorthogonal i.e., $$\langle u_\xi, v_\eta\rangle_{L^2(\Omega)}=\delta_{\xi\eta},$$ where the inner product $\langle\cdot, \cdot \rangle_{L^2(\Omega)}$ is the usual inner product on the Hilbert space $L^2(\Omega)$ given by
	$$\langle f, g\rangle_{L^2(\Omega)}=\int_{\Omega}f(x)\overline{g(x)}dx.$$ 	\textbf{Assumption 2:} The system $\{u_\xi\}_{\xi \in \mathcal{I}}$ is a basis for $L^2(\Omega)$.
	
	\begin{defn}\label{wz}
		The system $\{u_\xi: \xi \in \mathcal{I}\}$ is called a WZ-system if the functions
		$u_\xi(x), v_\xi(x)$ do not have zeros on the domain 
		for all $\xi \in \mathcal{I}$ and if there exist a constant $C > 0$
		and a non-negative integer $N $ such that
		$$\inf_{x\in \Omega}|u_\xi(x)|\geq C\langle \xi\rangle^{-N},$$
		$$\inf_{x\in \Omega}|v_\xi(x)|\geq C\langle \xi\rangle^{-N},$$
		as $|\xi|\to \infty.$ We refer to \cite{Ruz16} for examples and an extensive list of references in this subject.
	\end{defn}

	\subsection{Test functions for $\mathfrak{L}$ }
	In this subsection we recall the spaces of distributions generated by the differential operator $\mathfrak{L}$ and by its adjoint $\mathfrak{L}^*$  from \cite{DMT17,Ruz16,ruzhansky2017nonharmonic,MR&JPVR}.
	
	\begin{defn}
		The space $C_{\mathfrak{L}}^{\infty}(\overline{\Omega}):=$ $\operatorname{Dom}\left(\mathfrak{L}^{\infty}\right)$ is called the space of test functions for the differential operator $\mathfrak{L}.$ From \cite{ruzhansky2017nonharmonic}, it is defined by 
		$$
		\operatorname{Dom}\left(\mathfrak{L}^{\infty}\right):=\bigcap_{k=1}^{\infty} \operatorname{Dom}\left(\mathfrak{L}^{k}\right),
		$$
		where $\operatorname{Dom}\left(\mathfrak{L}^{k}\right)$ is the domain of $\mathfrak{L}^{k}$ and is  defined as
		$$
		\operatorname{Dom}\left(\mathfrak{L}^{k}\right):=\left\{f \in L^{2}(\Omega): \mathfrak{L}^{j} f \in L^{2}(\Omega), j=0,1,2, \ldots, k-1\right\}
		$$
		such that the boundary conditions (BC) are satisfied by each of  the operators $\mathfrak{L}^{j} $ for $j=0, 1 \cdots k-1.$
	\end{defn} 
	The Fréchet topology of $C_{\mathfrak{L}}^{\infty}(\overline{\Omega})$ is given by the family of norms
	$$
	\|\varphi\|_{C_{\mathfrak{L}}^{k}}:=\max _{j \leq k}\left\|\mathfrak{L}^{j} \varphi\right\|_{L^{2}(\Omega)},~~ k \in \mathbb{N}_{0}, ~\varphi \in C_{\mathfrak{L} }^{\infty}(\overline{\Omega}).
	$$
	Analogously, we can introduce the space $C_{\mathfrak{L}^{*}}^{\infty}(\Omega)$ corresponding to the adjoint $\mathfrak{L}_{\Omega}^{*}$ by
	$$
	C_{\mathfrak{L}^{*}}^{\infty}(\overline{\Omega}):=\operatorname{Dom}\left((\mathfrak{L}^{*})^{\infty}\right)=\bigcap_{k=1}^{\infty} \operatorname{Dom}\left(\left(\mathfrak{L}^{*}\right)^{k}\right),
	$$
	where $\operatorname{Dom}\left(\left(\mathfrak{L}^{*}\right)^{k}\right)$ is the domain of the operator $\left(\mathfrak{L}^{*}\right)^{k}$ given by
	$$
	\operatorname{Dom}\left(\left(\mathfrak{L}^{*}\right)^{k}\right):=\left\{f \in L^{2}(\Omega):\left(\mathfrak{L}^{*}\right)^{j} f \in L^{2}(\Omega), j=0,1,2, \ldots, k-1\right\}.
	$$
	which satisfy the adjoint boundary conditions corresponding to the operator $\mathfrak{L}^*.$
	The Fréchet topology of $C_{\mathfrak{L}^*}^{\infty}(\overline{\Omega})$ is given by the family of norms
	$$
	\|\varphi\|_{C_{\mathfrak{L}^*}^{k}}:=\max _{j \leq k}\left\|\left(\mathfrak{L}^{*}\right)^{j} \varphi\right\|_{L^{2}(\Omega)}, ~~k \in \mathbb{N}_{0}, ~\varphi \in C_{\mathfrak{L}^* }^{\infty}(\overline{\Omega}).
	$$
	\begin{rmk}For all $\xi \in \mathcal{I}$, since $u_\xi \in C_{\mathfrak{L}}^{\infty}(\overline{\Omega})$ and $v_\xi \in  C_{\mathfrak{L}^{*}}^{\infty}(\overline{\Omega}),$ we  observe that Assumption 1 implies that the spaces $C_{\mathfrak{L}}^{\infty}(\overline{\Omega})$ and $C_{\mathfrak{L}^*}^{\infty}(\overline{\Omega})$ are dense in $L^2(\Omega)$.	
	\end{rmk}
	The space
	$\mathcal{D}_\mathfrak{L}'(\Omega) :=\mathcal{L}(C_{\mathfrak{L}^{*}}^{\infty}(\overline{\Omega}), \mathbb{C} )$ (or $\mathcal{D}_{\mathfrak{L}^*}'(\Omega) :=\mathcal{L}(C_{\mathfrak{L}}^{\infty}(\overline{\Omega}), \mathbb{C} )$) of linear continuous functionals on $C_{\mathfrak{L}^{*}}^{\infty}(\overline{\Omega})$ (or $C_{\mathfrak{L}}^{\infty}(\overline{\Omega})$) is called the space of $\mathfrak{L}$-distributions (or $\mathfrak{L}^*$-distributions). Moreover, the space $ \mathcal{D}_{\mathfrak{L}}^{\prime}(\Omega) \otimes \mathcal{S}^{\prime}(\mathcal{I})$ is the collection of all functions $f:\Omega \times \mathcal{I} \rightarrow \mathbb{C}$ such that $f(x,.) \in \mathcal{S}^{\prime}(\mathcal{I})$ for every $x \in \Omega$ and $f(.,\xi) \in \mathcal{D}_{\mathfrak{L}}^{\prime}(\Omega)$ for every $\xi \in \mathcal{I}.$
	
	\subsection{$\mathfrak{L}$-Fourier transform}
	In this subsection we recall the definition and some important properties of $\mathfrak{L}$-Fourier transform generated by the boundary value problem $\mathfrak{L}$.
	
	Let $\mathcal{S}(\mathcal{I})$ denote the space of rapidly decaying functions $\phi: \mathcal{I} \rightarrow \mathbb{C}$ i.e., $\phi \in \mathcal{S}(\mathcal{I})$ if for every $\ell<\infty,$ there exists a constant $C_{\phi, \ell}$ such that
	$$
	|\phi(\xi)| \leq C_{\phi, \ell}\langle\xi\rangle^{-\ell}
	$$
	for all $\xi \in \mathcal{I}$. The topology on $\mathcal{S}(\mathcal{I})$ is given by the seminorms $p_{k}, k \in \mathbb{N}_{0}$ and is defined by 
	$$
	p_{k}(\phi):=\sup _{\xi \in \mathcal{I}}\langle\xi\rangle^{k}|\phi(\xi)|.
	$$
	The continuous linear functionals on $\mathcal{S}(\mathcal{I})$ are of the form $$\phi \mapsto\langle u, \phi \rangle=\sum_{\xi \in \mathcal{I}} u(\xi) \phi(\xi),$$ where  $u: \mathcal{I} \rightarrow \mathbb{C}$	has the property that it grows at most polynomially at infinity, i.e. there exist constants $\ell<\infty$ and $C_{u, \ell}$ such that 
	$$
	|u(\xi)| \leq C_{u, \ell}\langle\xi\rangle^{\ell}
	$$
	for all $\xi \in \mathcal{I}$. Such kind of distributions $u: \mathcal{I} \rightarrow \mathbb{C}$ form a distribution space which we  denote by $\mathcal{S}^{\prime}(\mathcal{I})$.
	\begin{defn}[Fourier transform]
		The $\mathfrak{L}$-Fourier transform
		$$
		\left(\mathcal{F}_{\mathfrak{L}} f\right)(\xi)=(f \mapsto \widehat{f}): C_{\mathfrak{L}}^{\infty}(\overline{ \Omega}) \rightarrow \mathcal{S}(\mathcal{I})
		$$
		is given by
		$$
		\hat{f}(\xi):=\left(\mathcal{F}_{\mathfrak{L}} f\right)(\xi)=\int_{\Omega} f(x) \overline{v_{\xi}(x)} ~d x.
		$$
	\end{defn}
	Analogously, the  $\mathfrak{L}^{*}$-Fourier transform 
	$$
	\left(\mathcal{F}_{\mathfrak{L} \cdot^{*}} f\right)(\xi)=\left(f \rightarrow \widehat{f}_{*}\right): C_{\mathfrak{L}^{*}}^{\infty}(\overline{ \Omega}) \rightarrow S(\mathcal{I})
	$$
	is defined by
	$$
	\hat{f}_*(\xi):=\left(\mathcal{F}_{\mathfrak{L}^{*}}  f\right)(\xi)=\int_{\Omega} f(x) \overline{u_{\xi}(x)} ~d x.
	$$
	\begin{prop}
		The $\mathfrak{L}$-Fourier transform $\mathcal{F}_{\mathfrak{L}}$ is a bijective homeomorphism from $C_{\mathfrak{L}}^{\infty}(\overline{ \Omega})$ into $ \mathcal{S}(\mathcal{I})$. The inverse of $\mathcal{F}_{\mathfrak{L}}$,
		$$\mathcal{F}_{\mathfrak{L}}^{-1}: \mathcal{S}(\mathcal{I}) \rightarrow C_{\mathfrak{L}}^{\infty}(\overline{ \Omega})$$ is given by
		$$
		\left(\mathcal{F}_{\mathfrak{L}}^{-1} h \right)(x)=\sum_{\xi \in \mathcal{I}} h(\xi) u_{\xi}(x),~ h \in \mathcal{S}(\mathcal{I})
		$$
		in order that the Fourier inversion formula is given by
		$$
		f(x)=\sum_{\xi \in \mathcal{I}} \hat{f}(\xi) u_\xi(x),~~ f \in C_{\mathfrak{L}}^{\infty}(\overline{ \Omega}).
		$$
	\end{prop}
	Similarly,  $\mathfrak{L}^*$-Fourier transform $\mathcal{F}_{\mathfrak{L}^{*}}: C_{\mathfrak{L}^{*}}^{\infty}(\overline{ \Omega}) \rightarrow \mathcal{S}(\mathcal{I})$ is a bijective homeomorphism and its inverse $\mathcal{F}_{\mathfrak{L}^{*}}^{-1}: \mathcal{S}(\mathcal{I}) \rightarrow C_{\mathfrak{L}^{*}}^{\infty}(\overline{ \Omega})$ is
	given by
	$$
	\left(\mathcal{F}_{\mathfrak{L}^{*}}^{-1} h\right)(x)=\sum_{\xi \in \mathcal{I}} h(\xi) v_{\xi}(x), h \in S(\mathcal{I})
	$$
	such  that the conjugate Fourier inversion formula is given by
	$$
	f(x)=\sum_{\xi \in \mathcal{I}} \widehat{f}_{*}(\xi) v_{\xi}(x), f \in C_{\mathfrak{L}^{*}}^{\infty}(\overline{ \Omega}).
	$$
	Since the systems  $u_\xi$ and  $v_\xi$ are Riesz bases, from the work of Bari \cite{bari},  we have the following relation between a $L^2$ function and its Fourier coffecient.

	\begin{rmk}\label{estimate of f and f hat}
		For every $f \in L^{2}(\Omega)$ there exist constants $m_1, m_2$  such that  
		$$
		m_1^{2}\|f\|_{L^{2}}^{2} \leq \sum_{\xi \in \mathcal{I}}|\hat{f}(\xi)|^{2} \leq m_2^{2}\|f\|_{L^{2}}^{2}
		$$
		and
		$$
		m_1^{2}\|f\|_{L^{2}}^{2} \leq \sum_{\xi \in I}\left|\hat{f}_{*}(\xi)\right|^{2} \leq m_2^{2}\|f\|_{L^{2}}^{2}.
		$$
	\end{rmk}
	\subsection{$\mathfrak{L}$-Convolution, Plancherel formula and Sobolev space}
	In this subsection we recall  Sobolev spaces generated by the model operator $\mathfrak{L}.$ We also recall the Plancherel identity obtained by introducing suitable sequence spaces $\ell^2(\mathfrak{L})$ and $\ell^2(\mathfrak{L}^*)$. We start this subsection by the definition of $\mathfrak{L}$-convolution using the $\mathfrak{L}$-Fourier transform.
	\begin{defn}
		For $f_1, f_2 \in C_{\mathfrak{L}}^{\infty}(\overline{ \Omega})$ the  $\mathfrak{L}$-convolution between $f_1$ and $f_2$ is defined by
		$$
		\left(f_1*_{\mathfrak{L}} f_2\right)(x):=\sum_{\xi \in \mathcal{I}} \widehat{f_1}(\xi) \widehat{f_2}(\xi)  u_\xi(x).
		$$
	\end{defn}
	Analogously, for $f_1, f_2 \in C_{\mathfrak{L}^{*}}^{\infty}(\overline{ \Omega})$ the  $\mathfrak{L}^{*}$-convolution using the $\mathfrak{L}^*$-Fourier transform is given by
	$$
	\left(f_1*_{\mathfrak{L}} f_2\right)(x):=\sum_{\xi \in \mathcal{I}} \widehat{f_{1_*}}(\xi) \widehat{f_{2_*}}(\xi)  v_\xi(x).
	$$ 
	\begin{prop}\cite{Ruz16}
		For any $f_1, f_2 \in C_{\mathfrak{L}}^{\infty}(\overline{ \Omega})$,
		$$
		\widehat{f_1*_{\mathfrak{L}} f_2}=\widehat{f_1} \widehat{f_2}.
		$$
	\end{prop}
	\begin{defn}
		Let $\ell^2(\mathfrak{L})$ be the linear space of complex-valued functions $g$ on  $\mathcal{I}$ such that $\mathcal{F}_{\mathfrak{L}}^{-1} g \in L^{2}(\Omega),$ i.e. if there exists $f \in L^{2}(\Omega)$ such that $\mathcal{F}_{\mathfrak{L}} f=g.$ Then the space of sequences $\ell^{2}({\mathfrak{L}})$ is a Hilbert space with the inner product
		$$
		\langle a, b\rangle_{\ell^{2}({\mathfrak{L}})}:=\sum_{\xi \in \mathcal{I}} a(\xi) \overline{\left(\mathcal{F}_{\mathfrak{L}^{*}} \circ \mathcal{F}_{\mathfrak{L}}^{-1} b\right)(\xi)}
		$$
		for arbitrary $a, b \in \ell^{2}({\mathfrak{L}}).$
	\end{defn}
	The norm on $\ell^{2}({\mathfrak{L}})$ is given by the formula  
	$$
	\|a\|_{\ell^{2}({\mathfrak{L}})}=\sum_{\xi \in \mathcal{I}} a(\xi) \overline{\left(\mathcal{F}_{\mathfrak{L}^{*}} \circ \mathcal{F}_{\mathfrak{L}}^{-1} a\right)(\xi)}, ~~\quad a\in \ell^{2}({\mathfrak{L}}).
	$$
	Analogously,  the Hilbert space  $\ell^2(\mathfrak{L}^*)$ be the linear space of complex-valued functions $g$ on  $\mathcal{I}$ such that $\mathcal{F}_{\mathfrak{L}^*}^{-1} g \in L^{2}(\Omega),$  with the inner product
	$$
	\langle a, b\rangle_{\ell^{2}({\mathfrak{L}^*})}:=\sum_{\xi \in \mathcal{I}} a(\xi) \overline{\left(\mathcal{F}_{\mathfrak{L}} \circ \mathcal{F}_{\mathfrak{L}^*}^{-1} b\right)(\xi)}
	$$
	for arbitrary $a, b \in \ell^2(\mathfrak{L}^*)$.
	The   sequence spaces $\ell^2(\mathfrak{L})$ 
	and  $\ell^2(\mathfrak{L}^*)$  are thus generated by the biorthogonal
	systems $\{u_\xi\}_{\xi \in \mathcal{I}}$ and $\{v_\xi\}_{\xi \in \mathcal{I}}$ respectively. The reason for this particular  choice of the definition is clear from the  following Plancherel's identity
	\begin{prop}[Plancherel's identity]
		If $f_1, f_2 \in L^{2}(\Omega)$ then $\widehat{f_1}, \widehat{f_2} \in \ell^{2}({\mathfrak{L}}), \widehat{f_1}_{*}, \widehat{f_2}_{*} \in \ell^{2}({\mathfrak{L}^*})$
		and the inner products  take the form
		$$
		\langle\widehat{f_1}, \widehat{f_2} \rangle_{\ell^{2}({\mathfrak{L}})}=\sum_{\xi \in \mathcal{I}} \widehat{f_1}(\xi) \overline{\widehat{f_2}_{*}(\xi)}
		$$
		and
		$$
		\langle\widehat{f_1}_{*}, \widehat{f_2}_{*} \rangle_{\ell^{2}({\mathfrak{L}^*})}=\sum_{\xi \in \mathcal{I}} \widehat{f_1}_{*}(\xi) \overline{\widehat{f_2}(\xi)}
		$$
		In particular, we have
		$$
		\overline{\langle\widehat{f_1}, \widehat{f_2} \rangle}_{\ell^{2}({\mathfrak{L}})}=\langle \widehat{f_2}_{*}, \widehat{f_1}_{*}\rangle_{\ell^{2}({\mathfrak{L}^*})}
		$$
		and the Parseval identity takes of the form
		$$
		\langle {f_1},  {f_2} \rangle_{L^2(\Omega)}=\langle\widehat{f_1}, \widehat{f_2} \rangle_{\ell^{2}({\mathfrak{L}})}=\sum_{\xi \in \mathcal{I}} \widehat{f_1}(\xi) \overline{\widehat{f_2}_{*}(\xi)}=\sum_{\xi \in \mathcal{I}} \widehat{f_1}_{*}(\xi) \overline{\widehat{f_2}(\xi)}.
		$$
		Moreover, for any $f \in L^{2}(\Omega),$ we have $\widehat{f} \in \ell^{2}({\mathfrak{L}}), \widehat{f}_{*} \in \ell^{2}({\mathfrak{L}^*})$ and
		$$
		\|f\|_{L^2(\Omega)}=\|\widehat{f}\|_{\ell^{2}({\mathfrak{L}})}=\|\widehat{f}_{*}\|_{\ell^{2}({\mathfrak{L}^*})}.
		$$
	\end{prop}
	
	We now recall Sobolev spaces generated by the model operator $\mathfrak{L}$.
	\begin{defn}[Sobolev spaces $\mathcal{H}_{\mathfrak{L}}^{s}(\Omega)$]
		For $f \in \mathcal{D}_{\mathfrak{L}}^{\prime}(\Omega) \cap \mathcal{D}_{\mathfrak{L}^{*}}^{\prime}(\Omega)$ and $s \in \mathbb{R},$ we
		say that $f \in \mathcal{H}_{\mathfrak{L}}^{s}(\Omega)$ if and only if $\langle\xi\rangle^{s} \widehat{f}(\xi) \in  \ell^{2}({\mathfrak{L}^*})$. We define the norm on $\mathcal{H}_{\mathfrak{L}}^{s}(\Omega)$ by
		$$
		\|f\|_{\mathcal{H}_{\mathfrak{L}}^{s}(\Omega)}:=\left(\sum_{\xi \in \mathcal{I}}\langle\xi\rangle^{2 s} \widehat{f}(\xi) \overline{\widehat{f}_{*}(\xi)}\right)^{1/2}.
		$$
	\end{defn}
	The Sobolev space $\mathcal{H}_{ \mathfrak{L}}^{s}(\Omega)$ is then the space of $\mathfrak{L}$-distributions $f$ for which $\|f\|_{\mathcal{H}_{\mathfrak{L}}^{s}(\Omega)}<\infty.$ Likewise,  $\mathcal{H}_{\mathfrak{L}^{*}}^{s}(\Omega)$  is the space of $\mathfrak{L}^*$-distributions $f$ for which 
	$$
	\|f\|_{\mathcal{H}_{\mathfrak{L}^{*}}^s(\Omega)}:=\left(\sum_{\xi \in \mathcal{I}}\langle\xi\rangle^{2 s} \widehat{f}_{*}(\xi) \overline{\widehat{f}(\xi)}\right)^{1 / 2}<\infty.
	$$
	Note that the spaces $\mathcal{H}_{ \mathfrak{L}}^{s}(\Omega)= \mathcal{H}_{ \mathfrak{L}^*}^{s}(\Omega)$. For every $s \in \mathbb{R},$ the Sobolev space $\mathcal{H}_{\mathfrak{L}}^{s}(\Omega)$ is a Hilbert space with the inner product given by 
	$$
	\langle f, g\rangle_{\mathcal{H}_{\mathfrak{L}}^{s}(\Omega)}:=\sum_{\xi \in \mathcal{I}}\langle\xi\rangle^{2 s} \widehat{f}(\xi) \overline{\widehat{g}_{*}(\xi)}
	$$
	Similarly, with respect to the inner product  
	$$
	\langle f, g\rangle_{\mathcal{H}_{\mathfrak{L}^*}^{s}(\Omega)}:=\sum_{\xi \in \mathcal{I}}\langle\xi\rangle^{2 s} \widehat{f}_{*}(\xi) \overline{\widehat{g}(\xi)},
	$$
	the Sobolev space $\mathcal{H}_{\mathfrak{L}^{*}}^{s}(\Omega)$ is a Hilbert space for every $s \in \mathbb{R}$. Note that the Sobolev spaces $ \mathcal{H}_{\mathfrak{L}}^{s}(\Omega)$ and $\mathcal{H}_{\mathfrak{L}^{*}}^{s}(\Omega)$ are isometrically isomorphic and in particular $\mathcal{H}_{\mathfrak{L}}^{0}(\Omega)=L^2(\Omega)$.
	
	Now we recall the $p$-Lebesgue version of the space of Fourier coffecients, i.e., $\ell^p(\mathfrak{L})$-space associated to the model operator $\mathfrak{L}.$
	\begin{defn}
		The $\ell^p(\mathfrak{L})$ space associated to the model operator $\mathfrak{L}$ is defined by 
		$$ \ell^p(\mathfrak{L}):= \left\{
		\begin{array}{ll}
			\left\{a: \mathcal{I} \to \mathbb{C}: \displaystyle\sum_{\xi \in \mathcal{I}}|a(\xi)|^{p}\left\|u_{\xi}\right\|_{L^{\infty}(\Omega)}^{2-p}<\infty \right\}&\text{if} ~1\leq p \leq 2, \\\\
			\left\{a: \mathcal{I} \to \mathbb{C}: \displaystyle\sum_{\xi \in \mathcal{I}}|a(\xi)|^{p}\left\|v_{\xi}\right\|_{L^{\infty}(\Omega)}^{2-p}<\infty \right\}&\text{if} ~2\leq p <\infty,  \\\\
			\left\{a: \mathcal{I} \to \mathbb{C}:  \displaystyle \sup _{\xi \in \mathcal{I}}\left(|a(\xi)| \left\|v_{\xi}\right\|_{L^{\infty}(\Omega)}^{-1}\right)<\infty \right\}&\text{if} ~~p=\infty.
		\end{array} 
		\right.$$
		Analogously,  $\ell^p(\mathfrak{L}^*)$-space is defined by 
		$$ \ell^p(\mathfrak{L}^*):= \left\{
		\begin{array}{ll}
			\left\{a: \mathcal{I} \to \mathbb{C}: \displaystyle\sum_{\xi \in \mathcal{I}}|a(\xi)|^{p}\left\|v_{\xi}\right\|_{L^{\infty}(\Omega)}^{2-p}<\infty \right\}&\text{if} ~1\leq p \leq 2, \\\\
			\left\{a: \mathcal{I} \to \mathbb{C}: \displaystyle\sum_{\xi \in \mathcal{I}}|a(\xi)|^{p}\left\|u_{\xi}\right\|_{L^{\infty}(\Omega)}^{2-p}<\infty \right\}&\text{if} ~2\leq p <\infty,  \\\\
			\left\{a: \mathcal{I} \to \mathbb{C}:  \displaystyle \sup _{\xi \in \mathcal{I}}\left(|a(\xi)| \left\|u_{\xi}\right\|_{L^{\infty}(\Omega)}^{-1}\right)<\infty \right\}&\text{if} ~~p=\infty.
		\end{array} 
		\right.$$
	\end{defn}
	\subsection{$\mathfrak{L}$-admissible, $\mathfrak{L}$-quantization and full symbols}
	In this subsection we discuss the $\mathfrak{L}$-quantization of the $\mathfrak{L}$-admissible operator induced by our model operator $\mathfrak{L}$. From now onwards we assume that $\{u_\xi: \xi\in \mathcal{I}\}$ is a WZ-system in the sense of Definition \ref{wz}. 
	
	We recall that for any linear continuous operator
	$
	A: C_{\mathfrak{L}}^{\infty}(\overline{\Omega}) \rightarrow \mathcal{D}_{\mathfrak{L}}^{\prime}(\Omega),
	$ there exists a kernel $K_{A} \in \mathcal{D}_{\mathfrak{L}}^{\prime}(\Omega \times \Omega)$ such that for all $f \in C_{\mathfrak{L}}^{\infty}(\overline{\Omega})$,
	$$
	A f(x)=\int_{\Omega} K_{A}(x, y) f(y) ~d y.
	$$
	in
	the sense of distributions. Here $K_{A}$ is called the Schwartz kernel of $A$. Using the Fourier series formula we have
	$$
	A f(x)=\sum_{\eta \in \mathcal{I}} \widehat{f}(\eta) \int_{\Omega} K_{A}(x, y) u_{\eta}(y) d y.
	$$
	Analogously, for any linear continuous operator
	$
	A: C_{\mathfrak{L}^{*}}^{\infty}(\overline{\Omega}) \rightarrow \mathcal{D}_{\mathfrak{L}^{*}}^{\prime}(\Omega)
	$
	there exists a kernel $\widetilde{K}_{A} \in \mathcal{D}_{\mathfrak{L}^{*}}^{\prime}(\Omega \times \Omega)$ such that for all $f \in C_{\mathfrak{L}^{*}}^{\infty}(\overline{\Omega}),$  in
	the sense of distributions we can write
	$$
	A f(x)=\int_{\Omega} \tilde{K}_{A}(x, y) f(y) d y.
	$$
	\begin{defn}
		A linear continuous operator $
		A: C_{\mathfrak{L}}^{\infty}(\overline{\Omega}) \rightarrow \mathcal{D}_{\mathfrak{L}}^{\prime}(\Omega)
		$
		belongs to the class of $\mathfrak{L}$-admissible operators if
		$$
		\sum_{\eta \in \mathcal{I}} u_{\eta}^{-1}(x) u_{\eta}(z) \int_{\Omega} K_{A}(x, y) u_{\eta}(y) d y, \quad x,z \in \Omega,
		$$ is in $ \mathcal{D}_{\mathfrak{L}}^{\prime}(\Omega \times \Omega)$. Here  we assume that $\{u_\xi: \xi\in \mathcal{I}\}$ is a WZ-system in the sense of Definition \ref{wz}.
	\end{defn}
	\begin{rmk}
		Noted that the expression
		$$
		u_{\eta}^{-1}(x) u_{\eta}(z) \int_{\Omega} K_{A}(x, y) u_{\eta}(y) d y
		$$
		exists for any operator $A$ from the class of $\mathfrak{L}$-admissible operators. Furthermore it is in $ \mathcal{D}_{\mathfrak{L}}^{\prime}(\Omega) \otimes \mathcal{S}^{\prime}(\mathcal{I}).$
	\end{rmk} 
	
	Now we now recall the definition of $\mathfrak{L}$-symbol of an  $\mathfrak{L}$-admissible operator.
	\begin{defn}[ $\mathfrak{L}$-symbols of operators]
		The $\mathfrak{L}$-symbol of a linear continuous $\mathfrak{L}$-admissible operator
		$$
		A: C_{\mathfrak{L}}^{\infty}(\overline{\Omega}) \rightarrow \mathcal{D}_{\mathfrak{L}}^{\prime}(\Omega)
		$$
		is defined by
		$$
		\sigma_{A}(x, \xi):=u_{\xi}^{-1}(x) \int_{\Omega} K_{A}(x, y) u_{\xi}(y) d y.
		$$
		Also, we can also write
		$$
		\sigma_A(x, \xi)=\int_{\Omega} k_A(x, y) \overline{v_{\xi}(y)} \mathrm{d} y=\left\langle k_A(x), \bar{v}_{\xi}\right\rangle,
		$$
		where by the $\mathfrak{L}$-Fourier inversion formula, the convolution kernel can be regained from the symbol:
		$$
		k_A(x, y)=\sum_{\xi \in \mathcal{I}} \sigma_A(x, \xi) u_{\xi}(y)
		$$
		all in the sense of $\mathfrak{L}$-distributions.
	\end{defn} 	 
	
	\begin{thm}[$\mathfrak{L}$-quantization] Let
		$$
		A: C_{\mathfrak{L}}^{\infty}(\overline{\Omega}) \rightarrow \mathcal{D}_{\mathfrak{L}}^{\prime}(\Omega)
		$$
		be a linear continuous $\mathfrak{L}$-admissible operator with $\mathfrak{L}$-symbol $\sigma_{A} \in \mathcal{D}_{\mathfrak{L}}^{\prime}(\Omega) \otimes \mathcal{S}^{\prime}(\mathcal{I}) .$ Then the $\mathfrak{L}$-quantization
		$$A f(x)=\sum_{\xi \in \mathcal{I}} \widehat{f}(\xi) \sigma_{A}(x, \xi) u_{\xi}(x)$$
		is true for all $f \in C_{\mathfrak{L}}^{\infty}(\overline{\Omega})$ and $x\in \Omega$. Then the $\mathfrak{L}$-symbol $\sigma_{A}$ of $A$ can be written as
		$$\sigma_{A}(x, \xi)=u_{\xi}^{-1}(x)\left(A u_{\xi}\right)(x)
		$$ for all $x\in \Omega$ and $\xi\in \mathcal{I}.$
	\end{thm}
	In view of the above theorem, from now onwards we will mainly interseted  in operators $A: C_{\mathfrak{L}}^{\infty}(\overline{\Omega}) \rightarrow \mathcal{D}_{\mathfrak{L}}^{\prime}(\Omega)$ from the class of $\mathfrak{L}$-admissible operators.
	\begin{rmk}
		Let $A: Span\{u_\xi\}\subseteq Dom(A)\subseteq L^2(\Omega)\to L^2(\Omega)$ be a linear continuous $\mathfrak{L}$-admissible operator. If there exist a measurable function $\sigma_A:\overline{ \Omega}\times \mathcal{I}\to \mathbb{C}$ such that
		\begin{align}\label{k}
			\sigma_{A}(x, \xi)u_{\xi}(x)=\left(A u_{\xi}\right)(x),
		\end{align} then the $\mathfrak{L} $-quantisation 
		$$A f(x)=\sum_{\xi \in \mathcal{I}} \widehat{f}(\xi) \sigma_{A}(x, \xi) u_{\xi}(x)$$ is true for all $f\in Span\{u_\xi\}$ and  the function $\sigma_{A}(x, \xi) $ does not need to be in  $$\mathcal{D}_{\mathfrak{L}}^{\prime}(\Omega) \otimes \mathcal{S}^{\prime}(\mathcal{I}).$$ It is enought that $$\sigma_{A}(\cdot, \xi) \in L^2(\overline{ \Omega}), ~\text{for each}~\xi \in \mathcal{I}. $$ For this reason we call the linear operator $A$ that satisfy  condition \ref{k} operator.
	\end{rmk}
	As a consequence of the above theorem and under the assumption that  the biorthogonal system $u_\xi$ is a WZ-system, we have the following formule for the symbol.
	
	\begin{Cor}
		The following  are  equivalent for the 
		$\mathfrak{L}$-symbols:
		\begin{enumerate}
			\item $\sigma_{A}(x, \xi)=\displaystyle \int_{\Omega} k_{A}(x, y) \overline{v_{\xi}(y)}~ d y$
			\item $\sigma_{A}(x, \xi)=u_{\xi}^{-1}(x)\left(A u_{\xi}\right)(x)$
			\item $\sigma_{A}(x, \xi)=\displaystyle u_{\xi}^{-1}(x) \int_{\Omega} K_{A}(x, y) u_{\xi}(y) d y$
			\item $\displaystyle K_{A}(x, y)=\sum_{\xi \in \mathcal{I}} u_{\xi}(x) \sigma_{A}(x, \xi) \overline{v_{\xi}(y)}$
		\end{enumerate}
	\end{Cor}
	Note that 	when $\left\{u_{\xi}: \xi \in \mathcal{I}\right\}$ is not a WZ-system, the $\mathfrak{L}$-symbol $\sigma_{A}$ of the operator $A$ is a function on $\overline{\Omega} \times \mathcal{I}$ for which the equality
	$$
	u_{\xi}(x) \sigma_{A}(x, \xi)=\int_{\Omega} K_{A}(x, y) u_{\xi}(y) d y
	$$
	holds for all $\xi$ in $\mathcal{I}$ and for $x \in \overline{\Omega}.$ 
	Similarly, we have an analogous notion of the $\mathfrak{L}^{*}$ -quantization.
	\begin{defn}
		A linear continuous operator $
		A: C_{\mathfrak{L}^*}^{\infty}(\overline{\Omega}) \rightarrow \mathcal{D}_{\mathfrak{L}^*}^{\prime}(\Omega)
		$
		belongs to the class of $\mathfrak{L}^*$-admissible operators if
		$$
		\sum_{\eta \in \mathcal{I}} v_{\eta}^{-1}(x) v_{\eta}(z) \int_{\Omega} K_{A}(x, y) v_{\eta}(y) d y
		$$ is in $ \mathcal{D}_{\mathfrak{L}^*}^{\prime}(\Omega \times \Omega)$. Here  we assume that $\{v_\xi: \xi\in \mathcal{I}\}$ is a WZ-system in the sense of Definition \ref{wz}.
	\end{defn}
	\begin{rmk}
		Noted that the expression
		$$
		v_{\eta}^{-1}(x) v_{\eta}(z) \int_{\Omega} K_{A}(x, y) v_{\eta}(y) d y
		$$
		exists for any operator $A$ from the class of $\mathfrak{L}^*$-admissible operators. Furthermore it is in $ \mathcal{D}_{\mathfrak{L}^*}^{\prime}(\Omega) \otimes \mathcal{S}^{\prime}(\mathcal{I}).$
	\end{rmk} 
	
	Now we now recall the definition of $\mathfrak{L}^*$-symbol of an  $\mathfrak{L}^*$-admissible operator.
	\begin{defn}[ $\mathfrak{L}^{*}$-symbols of operators]
		The $\mathfrak{L}^*$-symbol of a linear continuous $\mathfrak{L}^*$-admissible operator
		$$
		A: C_{\mathfrak{L}^*}^{\infty}(\overline{\Omega}) \rightarrow \mathcal{D}_{\mathfrak{L}^*}^{\prime}(\Omega)
		$$
		is defined by
		$$
		\tau_{A}(x, \xi):=v_{\xi}^{-1}(x) \int_{\Omega} K_{A}(x, y) v_{\xi}(y) d y.
		$$
		
	\end{defn} 	 
	
	\begin{thm}[$\mathfrak{L}^*$-quantization] Let
		$$
		A: C_{\mathfrak{L}^*}^{\infty}(\overline{\Omega}) \rightarrow \mathcal{D}_{\mathfrak{L}^*}^{\prime}(\Omega)
		$$
		be a linear continuous $\mathfrak{L}$-admissible operator with $\mathfrak{L}^*$-symbol $\tau_{A} \in \mathcal{D}_{\mathfrak{L}^*}^{\prime}(\Omega) \otimes \mathcal{S}^{\prime}(\mathcal{I}) .$ Then the $\mathfrak{L}^*$-quantization
		$$A f(x)=\sum_{\xi \in \mathcal{I}} \widehat{f}(\xi)\tau_{A}(x, \xi) v_{\xi}(x)$$
		is true for all $f \in C_{\mathfrak{L}^*}^{\infty}(\overline{\Omega})$ and $x\in \Omega$. Then the $\mathfrak{L}^{*}$-symbol $\tau_{A}$ of $A$ can be written as
		$$\tau_{A}(x, \xi)=v_{\xi}^{-1}(x)\left(A v_{\xi}\right)(x)
		$$ for all $x\in \Omega$ and $\xi\in \mathcal{I}.$
	\end{thm}
	In view of the above theorem, from now onwards we will mainly interseted  in operators $A: C_{\mathfrak{L}^*}^{\infty}(\overline{\Omega}) \rightarrow \mathcal{D}_{\mathfrak{L}^*}^{\prime}(\Omega)$ from the class of $\mathfrak{L}^*$-admissible operators.
	\begin{rmk}
		Let $A: Span\{v_\xi\}\subseteq Dom(A)\subseteq L^2(\Omega)\to L^2(\Omega)$ be a linear continuous $\mathfrak{L}^*$-admissible operator. If there exist a measurable function $\tau_A:\overline{ \Omega}\times \mathcal{I}\to \mathbb{C}$ such that
		\begin{align}\label{kk}
			\tau_{A}(x, \xi)v_{\xi}(x)=\left(A v_{\xi}\right)(x),
		\end{align} then the $\mathfrak{L}^* $-quantization 
		$$A f(x)=\sum_{\xi \in \mathcal{I}} \widehat{f}(\xi)\tau_{A}(x, \xi) v_{\xi}(x)$$
		is true for all $f\in Span\{v_\xi\}$ and  the function $\tau_{A}(x, \xi) $ does not need to be in  
		$$\mathcal{D}_{\mathfrak{L}^*}^{\prime}(\Omega) \otimes \mathcal{S}^{\prime}(\mathcal{I}).$$
		It is enought that $$\tau_{A}(\cdot, \xi) \in L^2(\overline{ \Omega}), ~\text{for each}~\xi \in \mathcal{I}. $$ Analogously to the $\mathfrak{L}$-quantization, we have the following
		corollary.
	\end{rmk}
	\begin{Cor}\label{formulae of * symbol}
		We also have the following equivalent formulae for the $\mathfrak{L}^*$-symbols.
		\begin{enumerate}
			\item $\tau_{A}(x, \xi)=\displaystyle \int_{\Omega} \tilde{K}_{A}(x, y) \overline{u_{\xi}(y)}~ d y$
			\item $\tau_{A}(x, \xi)=v_{\xi}^{-1}(x)\left(A v_{\xi}\right)(x)$
			\item $\tau_{A}(x, \xi)=\displaystyle v_{\xi}^{-1}(x) \int_{\Omega}\tilde{K}_{A}(x, y) v_{\xi}(y) d y$
			\item $\displaystyle \tilde{K}_{A}(x, y)=\sum_{\xi \in \mathcal{I}} v_{\xi}(x) \tau_{A}(x, \xi)\overline{u_{\xi}(y)}$.
		\end{enumerate}
	\end{Cor}
	In particular when the $\mathfrak{L}$-symbol of a continuous linear operator is independent of $x$ then the notaion of Fourier multiplier arises. Indeed we have the following definition related to our model operator $\mathfrak{L}.$
	\begin{defn}
		Let $A:  Dom(A)\subseteq L^2(\overline{\Omega})\to L^2(\overline{\Omega})$ be a continuous linear operator. We say  $A$ is a $\mathfrak{L}$-Fourier multiplier if it satisfies 
		$$\mathcal{F}_{\mathfrak{L}}\left( Af\right)(\xi)=\sigma(\xi)\mathcal{F}_{\mathfrak{L}}\left( f\right)(\xi), ~f\in Dom(A)$$
		for some $\sigma: \mathcal{I}\to \mathbb{C}.$ Similarly we have  $\mathfrak{L}^*$-Fourier multiplier.
		Let $B:  Dom(A)\subseteq L^2(\overline{\Omega})\to L^2(\overline{\Omega})$ be a continuous linear operator. We say  $B$ is a $\mathfrak{L}^*$-Fourier multiplier if it satisfies 
		$$\mathcal{F}_{\mathfrak{L}^*}\left( Bf\right)(\xi)=\tau(\xi)\mathcal{F}_{\mathfrak{L}^*}\left( f\right)(\xi), ~f\in Dom(B)$$
		for some $\tau: \mathcal{I}\to \mathbb{C}.$
	\end{defn}
	Realtion betwen symbols of a $\mathfrak{L}$-Fourier multiplier and its adjoint is given in the next proposition. 
	\begin{prop}\cite{Ruz16}
		The continuous linear operator $A$ is an $\mathfrak{L}$-Fourier multiplier by $\sigma(\xi)$ if and only if  $A^*$ is an $\mathfrak{L}^*$-Fourier multiplier by $\overline{\sigma(\xi)}.$ 
	\end{prop} 
	\subsection{Difference operator}
	In this subsection we recall  difference operators given in \cite{Ruz16,ruzhansky2017nonharmonic}.  We start this subsection by the following definitions.
	\begin{defn}[$\mathfrak{L}$-strongly admissible functions]
		Denote  $$C_b^{\infty}(\Omega \times \Omega)=C^{\infty}(\Omega \times \Omega)\cap C^{\infty}(\overline{\Omega} \times \overline{\Omega})$$ and 
		let $q_{j} \in C_b^{\infty}(\Omega \times \Omega), j=1, \ldots, l,$ be a given family of smooth functions. The collection of $q_{j}$'s are called  $\mathfrak{L}$-strongly admissible if the following properties holds:
		\begin{enumerate}
			\item For every $x\in \Omega,$ the multiplication by $q_{j}(x, \cdot )$ is a continuous linear mapping on $C_{\mathfrak{L}^{*}}^{\infty}(\overline{\Omega})$ for all $j=1, \ldots, \ell;$
			\item $ q_{j}(x, x)=0$ for  all $j=1, \ldots, \ell;$
			\item $\left.\operatorname{rank}\left(\nabla_{y} q_{1}(x, y), \ldots,; \nabla_{y} q_{l}(x, y)\right)\right|_{y=x}=\dim(\Omega)$ for all $x \in \Omega$;
			\item the diagonal in $\Omega \times \Omega$ is the only set when all of $q_{j}$'s vanish:
			$$
			\bigcap_{j=1}^{l}\left\{(x, y) \in \Omega \times \Omega: q_{j}(x, y)=0\right\}=\{(x, x): x \in \Omega\}.
			$$
		\end{enumerate}
	\end{defn}
	The collection of $q_{j}$'s which satisfies the above propertics generalizes the notion of a strongly admissible collection of functions for difference operators in the context of compact Lie groups introduced in \cite{wirth}. Here, we will use the multi-indices notation
	$$
	q^{\alpha}(x, y):=q_{1}^{\alpha_{1}}(x, y) \cdots q_{l}^{\alpha_{l}}(x, y).
	$$
	Analogously, the notion of $\mathfrak{L}^*$-strongly admissible functions is given by.
	\begin{defn}[$\mathfrak{L}^*$-strongly admissible functions]
		let $p_{j} \in C_b^{\infty}(\Omega \times \Omega), j=1, \ldots, l,$ be a given family of smooth functions. The collection of $p_{j}$'s are called  $\mathfrak{L}^*$-strongly admissible if the following properties holds:
		\begin{enumerate}
			\item For every $x\in \Omega,$ the multiplication by $p_{j}(x, \cdot )$ is a continuous linear mapping on $C_{\mathfrak{L}^{*}}^{\infty}(\overline{\Omega})$ for all $j=1, \ldots, \ell;$
			\item $ p_{j}(x, x)=0$ for  all $j=1, \ldots, \ell;$
			\item $\left.\operatorname{rank}\left(\nabla_{y} p_{1}(x, y), \ldots,; \nabla_{y} p_{l}(x, y)\right)\right|_{y=x}=\dim(\Omega)$ for all $x \in \Omega$;
			\item the diagonal in $\Omega \times \Omega$ is the only set when all of $p_{j}$'s vanish:
			$$
			\bigcap_{j=1}^{l}\left\{(x, y) \in \Omega \times \Omega: p_{j}(x, y)=0\right\}=\{(x, x): x \in \Omega\}.
			$$
		\end{enumerate}
	\end{defn}
	Here we have used
	$$
	p^{\alpha}(x, y):=p_{1}^{\alpha_{1}}(x, y) \cdots p_{l}^{\alpha_{l}}(x, y).
	$$
	
	Now we recall  the Taylor expansion formula with respect to  $\mathfrak{L}$-strongly admissible family of $q_{j}$'s follows from expansions of functions $g$ and $q^{\alpha}(e, \cdot)$ by the common Taylor series.
	\begin{prop}
		Any smooth function $g \in C^{\infty}(\Omega)$ can be approximated by Taylor polynomial type expansions, i.e. for any $e \in \Omega,$ we have
		
		\begin{align*}
			g(x)&=\sum_{|\alpha|<N} \frac{1}{\alpha !} D_{x}^{(\alpha)} g(x)\big|_{x=e} q^{\alpha}(e, x)+\sum_{|\alpha|=N} \frac{1}{\alpha !} q^{\alpha}(e, x) g_{N}(x) \\
			&\sim\sum_{\alpha\geq 0} \frac{1}{\alpha !} D_{x}^{(\alpha)} g(x)\big|_{x=e} q^{\alpha}(e, x)
		\end{align*}
		in a neighborhood of $e \in \Omega,$ where $g_{N} \in C^{\infty}(\Omega)$ and $D_{x}^{(\alpha)} g(x)\big|_{x=e} $ can be calculated from the recurrent relation: $D_{x}^{(0,\cdots, 0)}:=I$ and for $\alpha \in \mathbb{N}_{0}^{\ell}$ we have
		$$
		\partial_{x}^{\beta} g(x)\big|_{x=e}=\sum_{|\alpha| \leq|\beta|} \frac{1}{\alpha !}\left[\partial_{x}^{\beta}{q}^{\alpha}(e, x)\right]\big|_{x=e} {D}_{x}^{(\alpha)} g(x)\big|_{x=e},
		$$
		where $\beta=(\beta_1, \cdots, \beta_n)$ and $\partial_{x}^{\beta}=\frac{\partial^{\beta_1}}{\partial x_1^{\beta_1}}\cdots \frac{\partial^{\beta_n}}{\partial x_n^{\beta_n}}$.
	\end{prop}
	Similarly, for the adjoint problem, any smooth function $g \in C^{\infty}(\Omega)$ can be approximated by Taylor polynomial type expansions  given by 
	\begin{align*}i
		g(x)&=\sum_{|\alpha|<N} \frac{1}{\alpha !} \tilde{D}_{x}^{(\alpha)} g(x)\big|_{x=e} \tilde{q}^{\alpha}(e, x)+\sum_{|\alpha|=N} \frac{1}{\alpha !} \tilde{q}^{\alpha}(e, x) g_{N}(x) \\
		&\sim\sum_{\alpha\geq 0} \frac{1}{\alpha !}\tilde{D}_{x}^{(\alpha)} g(x)\big|_{x=e} \tilde{q}^{\alpha}(e, x)
	\end{align*}
	in a neighborhood of $e \in \Omega,$ where $g_{N} \in C^{\infty}(\Omega)$ and $\tilde{D}_{x}^{(\alpha)} g(x)\big|_{x=e} $ can be calculated from the recurrent relation: $\tilde{D}_{x}^{(0,\cdots, 0)}:=I$ and for $\alpha \in \mathbb{N}_{0}^{\ell},$ we have
	$$
	\partial_{x}^{\beta} g(x)\big|_{x=e}=\sum_{|\alpha| \leq|\beta|} \frac{1}{\alpha !}\left[\partial_{x}^{\beta}\tilde{q}^{\alpha}(e, x)\right]\big|_{x=e} \tilde{D}_{x}^{(\alpha)} g(x)\big|_{x=e},
	$$
	where $\beta=(\beta_1, \cdots, \beta_n)$ and $\partial_{x}^{\beta}=\frac{\partial^{\beta_1}}{\partial x_1^{\beta_1}}\cdots \frac{\partial^{\beta_n}}{\partial x_n^{\beta_n}}$.
	The operators $D^{(\alpha)}$ and $\tilde{D}^{(\alpha)}$ are differential operators of order $|\alpha| $ and the operator $\partial_{x}^{\beta}$ can be written as a linear combination of  $D^{(\alpha)}$ or $\tilde{D}^{(\alpha)}$ with smooth bounded coefficients. Also, note that difference operators applied to symbols will depend on a point $x$ since the problem may lack invariance or symmetry. Now we recall the definition of difference operators given in \cite{ruzhansky2017nonharmonic}. 
	\begin{defn}
		Let
		$$
		A: C_{\mathfrak{L}}^{\infty}(\overline{\Omega}) \rightarrow \mathcal{D}_{\mathfrak{L}}^{\prime}(\Omega)
		$$
		be a linear continuous $\mathfrak{L}$-admissible operator with $\mathfrak{L}$-symbol $\sigma_{A} \in \mathcal{D}_{\mathfrak{L}}^{\prime}(\Omega) \otimes \mathcal{S}^{\prime}(\mathcal{I}) $ with Schwartz kernel $ K_{A}(x, y)\in \mathcal{D}_{\mathfrak{L}}^{\prime}(\Omega\times \Omega)$. Then the difference operator 
		$$ \Delta_q^\alpha: \mathcal{D}_{\mathfrak{L}}^{\prime}(\Omega) \otimes \mathcal{S}^{\prime}(\mathcal{I})\to \mathcal{D}_{\mathfrak{L}}^{\prime}(\Omega) \otimes \mathcal{S}^{\prime}(\mathcal{I}) ,$$
		is an action on $\mathfrak{L}$-symbol $\sigma_{A},$ given by 
		$$\Delta_q^\alpha\sigma_{A}(x, \xi)=u_{\xi}^{-1}(x) \int_{\Omega}{q}^{\alpha}( x, y) K_{A}(x, y) u_{\xi}(y) d y.$$
	\end{defn}
	Analogously, for a $\mathfrak{L}^*$-admissible operator 
	$$
	A: C_{\mathfrak{L}^*}^{\infty}(\overline{\Omega}) \rightarrow \mathcal{D}_{\mathfrak{L}^*}^{\prime}(\Omega)
	$$
	with $\mathfrak{L}$-symbol $\tau_{A} \in \mathcal{D}_{\mathfrak{L}^*}^{\prime}(\Omega) \otimes \mathcal{S}^{\prime}(\mathcal{I}),$ with Schwartz kernel $ \tilde{K}_{A}(x, y)\in \mathcal{D}_{\mathfrak{L}^*}^{\prime}(\Omega\times \Omega)$ the difference operator 
	$$ \tilde{\Delta}_q^\alpha: \mathcal{D}_{\mathfrak{L}}^{\prime}(\Omega) \otimes \mathcal{S}^{\prime}(\mathcal{I})\to \mathcal{D}_{\mathfrak{L}}^{\prime}(\Omega) \otimes \mathcal{S}^{\prime}(\mathcal{I}) ,$$
	is an action on $\mathfrak{L}$-symbol $\tau_{A},$ given by 
	$$\tilde{\Delta}_q^\alpha\tau_{A}(x, \xi)=v_{\xi}^{-1}(x) \int_{\Omega}\tilde{q}^{\alpha}( x, y) \tilde{K}_{A}(x, y) v_{\xi}(y) d y.$$
	
	\subsection{Symbolic calculus}
	Using the difference operators and the derivatives $D^{(\alpha)}$ from the previous subsection we recall the classes of symbols  corresponding to our present setting. 
	\begin{defn}[Symbol classes $S_{\rho, \delta}^{m}(\Omega \times \mathcal{I})$]
		The $\mathfrak{L}$-symbol class   $S_{\rho, \delta, \Lambda}^{m}(\overline{ \Omega}\times \mathcal{I})$  consist of  such symbols $\sigma(x, \xi)$ which are in $ C_{\mathfrak{L}}^{\infty} \left(\overline{\Omega} \right) $ for all $\xi \in \mathcal{I} $ and  it satisfies
		\begin{align*}
			\left|{\Delta}_{(x)}^{\alpha}{D}_{x}^{(\beta)} \sigma(x, \xi)\right| \leq C_{\sigma, \alpha, \beta, m}\langle \xi \rangle^{m-\rho|\alpha|+\delta|\beta|}
		\end{align*}
		for all $\alpha, \beta \geq 0$ and $x\in \overline{\Omega}$.
	\end{defn} 
	
	Similarly, the $\mathfrak{L}^*$-symbol class   $\tilde{S}_{\rho, \delta, \Lambda}^{m}(\overline{ \Omega}\times \mathcal{I})$ consist of  such symbols $\tau(x, \xi)$ which are in $ C_{\mathfrak{L}^*}^{\infty} \left( \overline{\Omega} \right) $ for all $\xi \in \mathcal{I} $ and  it satisfies
	\begin{align*}
		\left|\tilde{\Delta}_{(x)}^{\alpha} \tilde{D}_{x}^{(\beta)}\tau(x, \xi)\right| \leq C_{\sigma, \alpha, \beta, m}\langle \xi \rangle^{m-\rho|\alpha|+\delta|\beta|}
	\end{align*}
	for all $\alpha, \beta \geq 0$ and $x\in \overline{\Omega}$.

	\begin{defn}[Pseudo-differential operators]
		Let $\sigma: \overline{ \Omega}\times \mathcal{I} \to \mathbb{C}$ be a measurable function such that $\sigma(\cdot, \xi) \in L^2(\Omega)$ for all $\xi \in \mathcal{I}$. Then the $\mathfrak{L}$-pseudo-differential operator corresponding to $\sigma $ is given by
		$$T_\sigma  f(x):=\sum_{\xi \in\mathcal{I}}\sigma(x, \xi) \hat{f}(\xi) u_{\xi}(x)$$
		for all $f\in Span\{u_\xi\}$. This $\sigma(x, \xi)$ is called the $\mathfrak{L}$-symbol of the pseudo-differential operator $T_\sigma.$
	\end{defn}
	Analogously, for a measurable function $\tau: \overline{ \Omega}\times \mathcal{I} \to \mathbb{C}$  such that $\tau(\cdot, \xi) \in L^2(\Omega)$ for all $\xi \in \mathcal{I}$. Then the $\mathfrak{L}^*$-pseudo-differential operator corresponding to $\tau $ is given by
	$$T_\tau  f(x):=\sum_{\xi \in\mathcal{I}}\tau(x, \xi) \hat{f}_*(\xi) v_{\xi}(x)$$
	for all $f\in Span\{v_\xi\}$. 
	\maketitle
	\section{Weighted $\mathfrak{L}$-Symbol Classes and $M$-elliptic pseudo-differential operators}\label{sec3}
 In this section, we introduce a weighted $\mathfrak{L}$-symbol class   $M_{\rho, 0, \Lambda}^{m}(\bar{\Omega} \times \mathcal{I}), m\in \mathbb{R},$ associated with a suitable weight function  $\Lambda$ on a countable set $\mathcal{I} $ and study several important properties of the weighted symbol class $M_{\rho, 0,\Lambda}^{m}(\bar{\Omega} \times \mathcal{I})$.  We also investigate elements of symbolic calculus for pseudo-differential operators associated with $M_{\rho, 0, \Lambda}^{m}, m\in \mathbb{R},$ by deriving formulae for the composition, adjoint, transpose.   We	construct the parametrix of $M$-elliptic  pseudo-differential operators by recalling the notion of $M$-ellipticity for symbols belonging to $\mathfrak{L}$-symbol class, $M_{\rho,0, \Lambda}^m(\bar{\Omega} \times \mathcal{I}).$

 We begin with the definition of the Weight function which will be used to define  $\mathfrak{L}$-symbol classes   $M_{\rho, 0, \Lambda}^{m}$.
	\begin{defn}
		[Weight function] Let $ \Lambda \in C_{ \mathfrak{L} }^{ \infty }\left(\mathcal{I} \right)$. We say that $\Lambda$ is a weight function if there exist
		suitable positive constants $\mu_{0} \leq \mu_{1} \leq \mu$ and $C_{0}, C_{1}$ such that for any multi-indices $ \alpha, \gamma \geq 0$
		such that $\gamma_j \in\{0,1\}$ for all $j$ and some positive $C_{\alpha, \gamma}$ depending only on $\alpha, \gamma$, such that 
		\begin{align} \label{1}
			C_{0}( 1+  |\xi|  )^{\mu_{0} } \leq \Lambda( \xi) \leq C_{1} (1+  |\xi| )^{\mu_{1} }, & ~\xi \in \mathcal{I}. \\
			\left| \xi^{\gamma} \Delta_{(x)}^{\alpha+\gamma } \Lambda(\xi)\right| \leq C_{\alpha, \gamma} \Lambda(\xi)^{1-(1 / \mu)|\alpha|}, &~ \xi \in\mathcal{I}.
		\end{align}
	\end{defn}

	\begin{defn}\label{020}
		[Class related to weight function] Let $m \in \mathbb{R}$ and let $\rho \in ( 0,1 / \mu]$. Then we define the $\mathfrak{L}$-symbol class   $S_{\rho, 0, \Lambda}^{m}$ to be the class of all such symbols $\sigma(x, \xi)$ which are in $ C_{\mathfrak{L}}^{\infty} \left( \bar{\Omega} \right) $ for all $\xi \in \mathcal{I} $ and  it satisfies
		\begin{align}
			\left|\Delta_{(x)}^{\alpha} D_{x}^{(\beta)} \sigma(x, \xi)\right| \leq C_{\sigma, \alpha, \beta, m}\Lambda(\xi)^{m-\rho|\alpha|}
		\end{align}
		for all  $(x, \xi) \in \bar{\Omega}\times \mathcal{I} $ and for all multi-indices $\alpha, \beta \geq 0$.
	\end{defn}
	Analogously,  the $\mathfrak{L}^*$-symbol class   $\tilde{S}_{\rho, 0, \Lambda}^{m}$ to be the class of all such symbols $\sigma(x, \xi)$ which are in $ C_{\mathfrak{L}^*}^{\infty} \left( \bar{\Omega} \right) $ for all $\xi \in \mathcal{I} $ and  it satisfies
	\begin{align}
		\left|\tilde{\Delta}_{(x)}^{\alpha} \tilde{D}_{x}^{(\beta)} \sigma(x, \xi)\right| \leq C_{\sigma, \alpha, \beta, m}\Lambda(\xi)^{m-\rho|\alpha|}
	\end{align}
	for all  $(x, \xi) \in \bar{\Omega}\times \mathcal{I} $, for all multi-indices $\alpha, \beta \geq 0$.
	
	\begin{rmk}
		When $\Lambda(\xi)=(1+|\lambda_\xi|^{2})^{\frac{1}{2m}},  \xi \in \mathcal{I}$, then $S_{\rho, 0, \Lambda}^{m}$ is the same as the Hörmander class $ S_{\rho, 0}^{m}(\bar{\Omega} \times \mathcal{I})$ for $m \in \mathbb{R}$ and $\rho \in (0,1] $.
	\end{rmk} 
	
	Let $\sigma  \in S_{\rho, 0, \Lambda}^{m}$. 
	Then, we define the $\mathfrak{L}$-pseudo-differential operator $T_{\sigma}$ associated to the symbol $\sigma$ by 
	$$
	T_\sigma f(x)=\sum_{\xi \in \mathcal{I}} \sigma (x, \xi) \widehat{f}(\xi) u_{\xi}(x), ~x\in \bar{\Omega}
	$$
	for every $f \in C_{\mathfrak{L}}^{\infty}(\bar{\Omega})$.
	
	Analogously, if $\tau  \in S_{\rho, 0, \Lambda}^{m}$,   the corresponding  $\mathfrak{L}^*$-pseudo-differential operator  is defined by
	$$
	T_\tau f(x)=\sum_{\xi \in \mathcal{I}} \tau (x, \xi) \widehat{f}_{*}(\xi) v_{\xi}(x), ~x\in \bar{\Omega}
	$$
	for every $f \in C_{\mathfrak{L}^*}^{\infty}(\bar{\Omega})$.

	\begin{defn}\label{0}
		For $m \in \mathbb{R}$ and  $\rho \in ( 0,1 / \mu]$,   we define the $\mathfrak{L}$-symbol class   $M_{\rho, 0, \Lambda}^{m}$ to be the class of all such symbols $\sigma(x, \xi)$ which are in $ C_{\mathfrak{L}}^{\infty} \left( \bar{\Omega} \right) $ for all $\xi \in \mathcal{I} $, such that for all $\gamma \geq 0$  with  $\gamma_j \in\{0,1\},$ we have
		\begin{align}
			\xi^\gamma \Delta_{(x)}^{\gamma} \sigma(x, \xi)\in S_{\rho, 0, \Lambda}^{m}.
		\end{align}
		
	\end{defn}
	
	\begin{rmk}	With the assumptions in Definitions  \ref{0}, let us notice that the class $M_{\rho, 0, \Lambda}^{m}$ also can be characterized by the following estimate: For every compact set $K\subset \Omega$	 \begin{align}	\left| \xi^\gamma \Delta_{(x)}^{\alpha+\gamma}D_{x}^{(\beta)} \sigma(x, \xi)\right| \leq C_{\sigma, \alpha, \beta, m}\Lambda(\xi)^{m-\rho|\alpha|}	\end{align}	for all  $(x, \xi) \in K \times \mathcal{I} $. 
	\end{rmk}
	
	Using the Definiton \ref{020} and \ref{0} we have the results.
	\begin{prop}\label{3}
		For any $m_1, m_2 \in \mathbb{R}, 0< \rho \leq \frac{1}{\mu},  a(x, \xi) \in M_{\rho, 0, \Lambda}^{m_1}, $ and $b(x, \xi) \in M_{\rho, 0,  \Lambda}^{m_2},$  the following properties hold:
		\begin{enumerate}
			\item If $m_1 \leq m_2$ then $M_{\rho, 0, \Lambda}^{m_1}\subset M_{\rho, 0, \Lambda}^{m_2}$
			\item $a(x, \xi) b(x, \xi) \in M_{\rho, 0, \Lambda}^{m_1+m_2}$
			\item  $ \Delta_{(x)}^{\alpha} D_{x}^{(\beta)} a(x, \xi)  \in M_{\rho, 0,  \Lambda}^{m_1-\rho|\alpha|}$ for all  $\alpha, \beta \geq 0$.
		\end{enumerate}
	\end{prop}
	
	\begin{defn}
		[$M$-elliptic] A symbol $\sigma \in C_{\mathfrak{L}}^{\infty} \left( \bar{\Omega} \times \mathcal{I} \right)$ in $M_{\rho, 0, \Lambda}^{m}$ is said to be $M$-elliptic if there exist positive constants $C$ and $R$ such that
		\begin{align}
			|\sigma(x, \xi)| \geq C \Lambda(\xi)^{m}, \quad |\xi|  \geq R.
		\end{align}
		
	\end{defn}

	\begin{defn} [$M$-hypoelliptic]
		Let $\ell , m \in \mathbb{R}$ with $\ell \leq m$. Then we define the $\mathfrak{L}$-symbol class   $HM_{\rho, 0, \Lambda}^{m, \ell}$ to be the class of all such symbols $\sigma(x, \xi)$ which are in $ C_{\mathfrak{L}}^{\infty} \left( \bar{\Omega} \right) $ for all $\xi \in \mathcal{I} $ and  it satisfies the following conditions:
		\begin{enumerate}
			\item  There exist positive constants $C_{1}, C_{2}$ and $R$ such that
			\begin{align}
				C_{1} \Lambda(\xi)^{\ell} \leq|\sigma(x, \xi)| \leq C_{2} \Lambda(\xi)^{m}, \quad x \in \bar{\Omega},|\xi| \geq R.
			\end{align}
			\item There exists positive constant \(R\)  and for all $\alpha, \beta, \gamma \geq 0 $, $\gamma_j \in\{0,1\}$ such that 
			\begin{align}
				\left| \xi^\gamma \Delta_{(x)}^{\alpha+\gamma} D_{x}^{(\beta)} \sigma(x, \xi)\right| \leq C_{\sigma, \alpha, \beta} |\sigma(x, \xi) |\Lambda(\xi)^{-\rho|\alpha|}\quad x \in \bar{\Omega}, |\xi| \geq R.
			\end{align}
			
		\end{enumerate}
	\end{defn}

	A symbol $\sigma$ in $HM_{\rho, 0, \Lambda}^{m, \ell}$ is said to be $M$-hypoelliptic of order $(m, \ell)$.

	\begin{rmk}
		$T_{\sigma}$ is $M$-elliptic of order $m$ if $\sigma \in HM_{\rho, 0, \Lambda}^{m, m}$.
	\end{rmk}

	\begin{rmk}
		$$
		HM_{\rho, 0, \Lambda}^{m, \ell}\subset M_{\rho, 0, \Lambda}^{m}\subset S_{\rho, 0, \Lambda}^{m}
		$$
		for every $\ell, m \in \mathbb{R}, \ell \leq  m$ and $\rho \in(0,1 / \mu]$.
		
	\end{rmk}
	
	\begin{lem}\label{2}
		For any $0 <\rho \leq \frac{1}{\mu}$, there exist  $N_{0}>0$ such that
		$$
		S_{\rho, 0, \Lambda}^{m-N_{0}} \subset M_{\rho, \Lambda}^{m} \subset S_{\rho, 0,  \Lambda}^{m}
		$$
		holds for every $m \in \mathbb{R}$. 
	\end{lem}
	\begin{proof}
		Let $\sigma(x, \xi) \in S_{\rho, 0, \Lambda}^{m-N_{0}}$. Then, for every compact
		set $K \subset {\Omega}$ and  for all multi-indices $\alpha, \beta ,  \gamma \geq 0 $   with  $\gamma_j \in\{0,1\},$ we have the following estimate: 
		$$
		\left| \xi^\gamma \Delta_{(x)}^{\alpha+\gamma} D_{x}^{(\beta)} \sigma(x, \xi)\right| \leq C_{\alpha, \beta, \gamma, K}(1+|\xi| )^{|\gamma|} \Lambda(\xi)^{m-N_{0}-\rho|\alpha+\gamma|}
		$$
		for  all  $(x, \xi) \in K \times \mathcal{I} $ and some positive constant $C_{\alpha, \beta, \gamma, K}$.  On the other hand, from
		the polynomial growth of $\Lambda(\xi)$ given in (\ref{1}), we have
		$$
		(1+|\xi|)^{|\gamma|} \leq C_{\gamma} \Lambda(\xi)^{\frac{1}{\mu_{0}}|\gamma|}, \quad \xi \in  \mathcal{I}.
		$$
		Since $\rho \leq \frac{1}{\mu} \leq \frac{1}{\mu_{0}}$ 
		\begin{align*}\left| \xi^\gamma \Delta_{(x)}^{\alpha+\gamma} D_{x}^{(\beta)} \sigma(x, \xi)\right| 
			& \leq C_{\alpha, \beta, \gamma, K}^{\prime} \Lambda(\xi)^{m-N_{0}+\left(\frac{1}{\mu_{0}}-\rho\right)|\gamma|-\rho|\alpha|} \\ 
			& \leq C_{\alpha, \beta, \gamma, K}^{\prime} \Lambda(\xi)^{m-\rho|\alpha|}
		\end{align*} 
		for  all  $(x, \xi) \in K \times \mathcal{I} $, therefore $\sigma(x, \xi) \in M_{\rho, 0, \Lambda}^{m}$. 
	\end{proof}
	\begin{thm} [Asymptotic sums of symbols]\label{4}  
		Suppose that  $\sigma_{j} \in M_{\rho, 0, \Lambda}^{m_{j}}$ for all $j \in \mathbb{N}_{0}$, where $\left\{ m_{j}\right\}_{j=0}^{\infty} \subset \mathbb{R}$ be a sequence such that $m_{j}>m_{j+1},$ and $m_{j} \rightarrow-\infty$ as $j \rightarrow \infty$.  Then there exists a $\mathfrak{L}$-symbol $\sigma \in M_{\rho, 0, \Lambda}^{m_{0}}$ such that for all $N \in \mathbb{N}_{0}$ $$\sigma \stackrel{m_{N, \rho, \delta}}{\sim} \sum_{j=0}^{N-1} \sigma_{j}.$$
	\end{thm}
	\begin{proof}
		Let us choose a function $\chi \in C^{\infty}(\mathbb{R}^k)$ by

		$$  \chi(\xi)= \left\{
		\begin{array}{ll}
			1, &\text{if} ~|\xi| \geq 1 \\
			0, & \text{ if} ~ |\xi| \leq \frac{1}{2} \\
		\end{array} 
		\right.,$$
		otherwise \(0 \leq \chi(\xi) \leq 1\). Let us consider a sequence
		$\left(\varepsilon_{j}\right)_{j=0}^{\infty}$ of positive real numbers such that $\varepsilon_{j}>\varepsilon_{j+1} \rightarrow 0$ as $j \rightarrow \infty$. Define  the function $\chi_{j} \in C^{\infty}(\mathbb{R}^k)$ by $$\chi_{j}(\xi):=\chi\left(\varepsilon_{j} \xi\right) .$$ Since $\chi_{j}(\xi)=1$ for sufficiently large $\xi,$ we
		get $\chi_{j} \sigma_{j} \in S_{\rho, 0, \Lambda}^{m_{j}}$ for each $j $. When \(j\) is large enough, for any fixed $\xi \in \mathcal{I},$ the function $\chi_{j}(\xi) \sigma_{j}(x, \xi)$
		vanishes. This justifies the definition
		$$
		\sigma(x, \xi):=\sum_{j=0}^{\infty} \chi_{j}(\xi) \sigma_{j}(x, \xi)
		$$
		and clearly $\sigma \in S_{\rho,  0, \Lambda}^{m_{0}}$. Moreover
		$$
		\sigma(x, \xi)-\sum_{j=0}^{N-1} \sigma_{j}(x, \xi)=\sum_{j=0}^{N-1}\left[\chi_{j}(\xi)-1\right] \sigma_{j}(x, \xi)+\sum_{j=N}^{\infty} \chi_{j}(\xi) \sigma_{j}(x, \xi).
		$$
		Since $\varepsilon_{j}>\varepsilon_{j+1},$ and $\varepsilon_{j} \rightarrow 0$ as $j \rightarrow \infty,$  then the  sum $\displaystyle\sum_{j=0}^{N-1}\left[\chi_{j}(\xi)-1\right] \sigma_{j}(x, \xi) $ 
		vanishes  whenever $\xi$ is large enough. This shows that
		$$
		\sigma(x, \xi)-\sum_{j=0}^{N-1} \chi_{j}(\xi) \sigma_{j}(x, \xi) \in S_{\rho, 0, \Lambda}^{m_{N}}.
		$$

		Since $m_{j} \rightarrow-\infty,$ as $j \rightarrow \infty,$ using the left
		inclusions in Lemma \ref{2} we have that $\displaystyle \sigma-\sum_{j=0}^{N-1} \sigma_{j} \in S_{\rho, 0, \Lambda}^{m_N} \subset M_{\rho, 0, \Lambda}^{m_0}$ for a sufficiently	large \(N .\) Then we infer \(\sigma(x, \xi) \in M_{\rho, 0,  \Lambda}^{m_0} .\) Furthermore, for all \(N \geq 2\) and \(N^{\prime} \geq N\)
		$$
		\sigma-\sum_{j=0}^{N-1} \sigma_{j} =\sum_{j=N}^{N'-1} \sigma_{j}+r_{N'}
		$$
		with $r_{N'} \in S_{\rho, 0, \Lambda}^{m_{N'}}$.  By choosing a sufficiently large $N'$ so that $m_{N'}<m_{N}-N_{0},$
		we have $r_{N'} \in M_{\rho, 0,  \Lambda}^{m_{N}}$ and then $\displaystyle \sigma-\sum_{j=0}^{N-1} \sigma_{j}  \in M_{\rho,0,  \Lambda}^{m_{N}}.$
	\end{proof}

	
	\begin{thm}[Adjoint operators]
		Let $T: C_{\mathfrak{L}}^\infty(\bar{\Omega})  \to C_{\mathfrak{L}}^\infty(\bar{\Omega}) $ be a continuous linear operator such that its $\mathfrak{L}$-symbol $\sigma_T \in M_{\rho, 0, \Lambda}^{m}$.	Assume that the conjugate symbol class $\tilde{S}_{\rho, 0, \Lambda}^{m}$ is defined with strongly admissible	functions $\widetilde{q}_{j}(x, y):=\bar{q}_{j}(x, y)$ which are $\mathfrak{L}^*$-strongly admissible. Then the adjoint  $T^*$ of $T$ is a $\mathfrak{L}^{*}$-pseudo-differential operator  with  $\mathfrak{L}^{*}$-symbol $\sigma_{T^*} \in M_{\rho, 0, \Lambda}^{m}$   having the asymptotic expansion
		$$
		\sigma_{T^*}(x, \xi) \sim \sum_{\alpha} \frac{1}{\alpha !} \widetilde{\Delta}_{(x)}^{\alpha} D_{x}^{(\alpha)} \overline{\sigma_T(x, \xi)}
		$$
	\end{thm}
	\begin{proof}
		From the definition of adjoint $$\langle T  u_\xi, v_\eta \rangle =\langle  u_\xi, T^* v_\eta \rangle$$
		$$
		\implies \int_{\Omega} T u_{\xi}(x) \overline{v_{\eta}(x)} d x=\int_{\Omega} u_{\xi}(x) \overline{T^{*} v_{\eta}(x)} d x
		$$
		for $\xi, \eta \in \mathcal{I}$.  Since $T_\sigma $ is an integral operator with Schwartz kernel $K_{A}$, therefore 
		\begin{align*}
			\int_{\Omega}\left[\int_{\Omega} K_{T }(x, y) u_{\xi}(y) d y\right] \overline{v_{\eta}(x)} d x
			&=\int_{\Omega} u_{\xi}(x)\left[\int_{\Omega} K_{T^{*}}(x, y) v_{\eta}(y) d y\right] d x\\
			&=\int_{\Omega} u_{\xi}(y)\left[\int_{\Omega} K_{T^{*}}(y, x) v_{\eta}(x) d x\right] d y\\
		\end{align*}
		for all $\xi, \eta \in \mathcal{I}$.  Consequently, we get the familiar property
		$$
		K_{T^{*}}(x, y)=\overline{K_{T}(y, x)}.
		$$
		Using the fact that $\tau$ is the $\mathfrak{L}^*$-symbol of \(T_\sigma ^{*}\), Corollary \ref{formulae of * symbol} and Taylor series expansion, we get
		\begin{align*}
			u_{\xi}(x) \sigma_{T^*}(x, \xi)
			&=\int_{\Omega} K_{T^{*}}(x, y) u_{\xi}(y) d y\\
			&=\int_{\Omega} \overline{K_{T}(y, x)}u_{\xi}(y) d y\\
			&=\int_{\Omega} \left[ \sum_{\eta \in \mathcal{I}} \overline{u_{\eta}(y) \sigma_T(y, \eta)} v_{\eta}(x) \right] v_{\xi}(y) d y\\
			&\sim \int_{\Omega} \left[ \sum_{\eta \in \mathcal{I}} \overline{u_{\eta}(y)} \left(\sum_{\alpha} \frac{1}{\alpha !} \overline{D_{x}^{(\alpha)} \sigma_T (x, \eta) q^{\alpha}(x, y)} \right)v_{\eta}(x)\right] v_{\xi}(y) d y\\
			&\sim\sum_{\eta \in \mathcal{I}} \left(\sum_{\alpha} \frac{1}{\alpha !} \overline{D_{x}^{(\alpha)} \sigma_T (x, \eta)} \right)v_{\eta}(x) \int_{\Omega} \overline{ q^{\alpha}(x, y) u_{\eta}(y)}   v_{\xi}(y) d y\\
			&\sim  \sum_{\alpha} \frac{1}{\alpha !}   \sum_{\eta \in \mathcal{I}} v_{\eta}(x) \overline{D_{x}^{(\alpha)} \sigma_T (x, \eta)}  \int_{\Omega} \overline{ q^{\alpha}(x, y) u_{\eta}(y)}   v_{\xi}(y) d y.
		\end{align*}
		Using the $ \mathfrak{L}^*$-version of the difference formula and by taking
		$$
		\widetilde{q}(x, y):=\overline{q(x, y)}
		$$	we get,
		
		\begin{align*}
			\sigma_{T^*}(x, \xi)
			& 	\sim  \sum_{\alpha} \frac{1}{\alpha !}   \sum_{\eta \in \mathcal{I}} v_{\eta}(x) \overline{D_{x}^{(\alpha)} \sigma_T (x, \eta)}  u_{\xi}(x)^{-1} \int_{\Omega} {\tilde{q}}^{\alpha}(x, y) \overline{u_{\eta}(y)}   v_{\xi}(y) d y\\
			&\sim \sum_{\alpha} \frac{1}{\alpha !} \widetilde{\Delta}_{(x)}^{\alpha} D_{x}^{(\alpha)} \overline{\sigma_T (x, \xi)}.
		\end{align*}
	\end{proof}
	
	\begin{thm} \label{comp}
		Let $m_1, m_2 \in \mathbb{R}$. 	Let $A, B: C_{\mathfrak{L}}^\infty(\bar{\Omega})  \to C_{\mathfrak{L}}^\infty(\bar{\Omega}) $  be continuous linear operator such that their  $\mathfrak{L}$-symbols  $\sigma_A  \in M_{\rho, 0, \Lambda}^{m_1}$ and $\sigma_B  \in M_{\rho, 0, \Lambda}^{m_2}$. Then $A B$  has  $\mathfrak{L}$-symbol  $\sigma_{AB}\in   M_{\rho, 0, \Lambda}^{m_1+m_2}$
		having asymptotic expansion
		\begin{align*}
			{\sigma_{AB}(x, \xi)} & {\sim \sum_{\alpha \geq 0} \frac{1}{\alpha !}\left(\Delta_{(x)}^{\alpha} \sigma_A (x, \xi)\right) D_x^{(\alpha)} \sigma_{B} (x, \xi)},
		\end{align*}
		where the asymptotic expansion means that for every $N \in \mathbb{N},$ we have $$\sigma_{AB} (x, \xi)-\sum_{|\alpha|<N} \frac{1}{\alpha !}\left(\Delta_{(x)}^{\alpha} \sigma_A (x, \xi)\right) D_x^{(\alpha)} \sigma_B (x, \xi) \in M_{\rho, 0, \Lambda} ^{m_{1}+m_{2}-\rho N}. $$
		
	\end{thm} 
	
	\begin{proof}
		Since $A$ is an integral operator with Schwartz kernel $k_{A }$, therefore  
		\begin{align*} 
			A B f(x) &=\left(k_{A}(x,.) *_{\mathfrak{L}} B f\right)(x) \\
			&=\int_{\Omega}\left[\int_{\Omega} F(x, y, z) k_{A}(x, z) d z\right](B  f)(y) d y \\ 
			&=\int_{\Omega}\left(\left[\int_{\Omega} F(x, y, z) k_{A}(x, z) d z\right]\right.\\ 
			&\left.\times \int_{\Omega}\left[\int_{\Omega} F(y, s, t) k_{B}(y, t) d t\right] f(s) d s\right) d y
		\end{align*}
		where the $\mathfrak{L}$-convoluton is defined by $$\left(f *_{\mathfrak{L}} g\right)(x):=\int_{\Omega} \int_{\Omega} F(x, y, z) f(y) g(z) d y d z,$$$$F(x, y, z)=\sum_{\xi \in \mathcal{I}} u_{\xi}(x) \overline{v_{\xi}(y)} \overline{v_{\xi}(z)}.$$
		Therefore 
		\begin{align*} \sigma_{AB}(x, \xi) &=u_{\xi}^{-1}(x)\left(A \left( T_\tau u_{\xi}\right)\right)(x) \\ &=u_{\xi}^{-1}(x) \int_{\Omega}\left(\left[\int_{\Omega} F(x, y, z) k_{A}(x, z) d z\right]\int_{\Omega}\left[\int_{\Omega} F(y, s, t) k_{B}(y, t) d t\right] u_{\xi}(s) d s\right) d y. \end{align*}
		Using Taylor polynomial expansion on $k_{B}(\cdot, t) \in C_{\mathfrak{L}}^{\infty}(\bar{\Omega})$ we have
		\begin{align*} 
			\sigma_{AB}(x, \xi) & \sim u_{\xi}^{-1}(x) \int_{\Omega}\left(\left[\int_{\Omega} F(x, y, z) k_{A}(x, z) d z\right]\right.\\ 
			&\left.\times \int_{\Omega}\left[\int_{\Omega} F(y, s, t) \sum_{\alpha \geq 0} \frac{1}{\alpha !} q^{\alpha}(x, y) D_{x}^{(\alpha)} k_{B}(x, t) d t\right] u_{\xi}(s) d s\right) d y \\ 
			&=\sum_{\alpha \geq 0} \frac{1}{\alpha !} u_{\xi}^{-1}(x) \int_{\Omega}\left(\left[\int_{\Omega} F(x, y, z) q^{\alpha}(x, y) k_{A}(x, z) d z\right]\right.\\ 
			&\left.\times \int_{\Omega}\left[\int_{\Omega} F(y, s, t) D_{x}^{(\alpha)} k_{B}(x, t) d t\right] u_{\xi}(s) d s\right) d y.
		\end{align*}
		Using the biorthogonality of $\left\{u_{\xi}\right\}_{\xi \in \mathcal{I}}$ and $\left\{v_{\xi}\right\}_{\xi \in \mathcal{I}},$ we have 
		\begin{align*}
			&u_{\xi}^{-1}(y) \int_{\Omega}\left[\int_{\Omega} F(y, s, t) D_{x}^{(\alpha)} k_{B}(x, t) d t\right] u_{\xi}(s) d s\\
			&=u_{\xi}^{-1}(y) \int_{\Omega}\left[\int_{\Omega} \left(\sum_{\eta \in \mathcal{I}} u_{\eta}(y) \overline{v_{\eta}(s)} \overline{v_{\eta}(t)} \right) D_{x}^{(\alpha)} k_{B}(x, t) d t\right] u_{\xi}(s) d s\\
			&=\sum_{\eta \in \mathcal{I}} u_{\xi}^{-1}(y) u_{\eta}(y) \int_{\Omega}\left[\int_{\Omega} D_{x}^{(\alpha)} k_{B}(x, t) \overline{v_{\eta}(t)} dt\right] u_{\xi}(s) \overline{v_{\eta}(s)} d s\\
			&=\sum_{\eta \in \mathcal{I}} u_{\xi}^{-1}(y) u_{\eta}(y)\left[\int_{\Omega} u_{\xi}(s) \overline{v_{\eta}(s)} d s\right] \times\left[\int_{\Omega} D_{x}^{(\alpha)} k_{B}(x, t) \overline{v_{\eta}(t)} d t\right]\\
			&=u_{\xi}^{-1}(y) u_{\xi}(y) D_{x}^{(\alpha)} (\widehat{k_{B}}(x, \xi))\\
			&=D_{x}^{(\alpha)} \overline{k_{B}}(x, \xi)\\
			&=D_{x}^{(\alpha)} \sigma_B(x, \xi).
		\end{align*}
		Therefore 
		\begin{align*} 
			\sigma_{AB}(x, \xi) & \sim \sum_{\alpha \geq 0} \frac{1}{\alpha !} u_{\xi}^{-1}(x) \int_{\Omega}\left(\left[\int_{\Omega} F(x, y, z) q^{\alpha}(x, y) k_{B}(x, z) d z\right]  u_{\xi}(y)  D_{x}^{(\alpha)} \sigma_{B}(x, \xi)\right) d y\\
			&=\sum_{\alpha \geq 0} \frac{1}{\alpha !} u_{\xi}^{-1}(x) D_{x}^{(\alpha)} \sigma_{B}(x, \xi)\int_{\Omega}\left(\left[\int_{\Omega} F(x, y, z) q^{\alpha}(x, y) k_{A}(x, z) d z\right]  u_{\xi}(y)  \right) d y\\
			&=\sum_{\alpha \geq 0} \frac{1}{\alpha !} D_{x}^{(\alpha)} \sigma_{B}(x, \xi)  u_{\xi}^{-1}(x) \int_{\Omega} q^{\alpha}(x, y)  \left[\int_{\Omega} F(x, y, z) k_{A}(x, z) d z\right]  u_{\xi}(y)  d y\\
			&=\sum_{\alpha \geq 0} \frac{1}{\alpha !} D_{x}^{(\alpha)} \sigma_{B}(x, \xi) u_{\xi}^{-1}(x) \int_{\Omega} q^{\alpha}(x, y)  k_{A}(x, y)  u_{\xi}(y)  d y\\
			&={ \sum_{\alpha \geq 0} \frac{1}{\alpha !}\left(\Delta_{(x)}^{\alpha} \sigma_A (x, \xi)\right) D_x^{(\alpha)} \sigma_{B} (x, \xi)}.
		\end{align*}
		
		Since  $\sigma_A  \in M_{\rho, 0, \Lambda}^{m_1}$  and $\sigma_{B}\in M_{\rho, 0, \Lambda}^{m_2},$ so  by the virtue of  Proposition \ref{3}, for every positive integer $N$
		$$ \sum_{|\alpha| = N} \frac{1}{\alpha !}\left(\Delta_{(x)}^{\alpha} \sigma_A (x, \xi)\right) D_x^{(\alpha)} \sigma_{B} (x, \xi) \in M_{\rho, 0, \Lambda}^{m_1+m_2-\rho N}.$$
		Applying Theorem \ref{4} to the sequence $\displaystyle \sum_{|\alpha| =N} \frac{1}{\alpha !}\left(\Delta_{(x)}^{\alpha} \sigma_A (x, \xi)\right) D_x^{(\alpha)} \sigma_{B}(x, \xi),$ we can conclude that $\sigma_{AB}\in M_{\rho, 0,  \Lambda}^{m_1+m_2}$.

	\end{proof}
	
	The following theorem is a special case of \cite[Theorem 14.1]{Ruz16}, which gives us the $L^2$- bounededness of the pseudo-differential operator with $\mathfrak{L}$-symbol in the class $M_{\rho, 0, \Lambda}^{0}$.
	\begin{thm}
		Let $\sigma \in M_{\rho, 0, \Lambda}^{0}$. Then  the associated $\mathfrak{L}$-pseudo differential operator $T_{\sigma}: L^2(\bar{\Omega}) \to L^2(\bar{\Omega})$ is a bounded linear operator.
	\end{thm}
	
	\begin{defn}\label{7}		
		A symbol $\sigma$ in $M_{\rho, 0, \Lambda}^{m}$ is said to be $M$-elliptic if there exist constants $C>0$ and $N_{0} \in \mathbb{N}$ such that $$|\sigma(x, \xi)|\geq C \Lambda(\xi)^m$$
		for all $(x, \xi) \in \bar{\Omega} \times \mathcal{I}$ for which $|\xi|\geq N_{0} $; this is equivalent to assuming that there exists $\sigma_{B} \in M_{\rho, 0, \Lambda}^{-m}$ such that $I-AB, I-BA$ are in $\mathrm{Op}_{\mathfrak{L}} M^{-\infty}$.
	\end{defn}

	The following lemma is analogous to the Lemma 2 in \cite{GM}, which can be proved in a similar way. So we skip the proof here.
	\begin{lem}\label{inverse of symbol class}
		Let $\sigma(x,\xi) \in M_{\rho, 0,\Lambda}^{m}(\overline{\Omega} \times \mathcal{I})$ and $\psi(x, \xi) \in C^{\infty}\left(\overline{\Omega} \times \mathbb{R}^{k}\right)$, be such that there exist two sufficiently large positive constants $R^{\prime},$ $R^{\prime \prime}$ with the property that $R^{\prime \prime}>R^{\prime},$ and $\psi(x, \xi)=0$, for all $x \in \overline{\Omega}$ and $|\xi| \leq R^{\prime}$, and $\psi(x, \xi)=1$, for all $x \in \overline{\Omega}$ and $|\xi| \geq R^{\prime \prime}$. Then $q(x,\xi)=\frac{\psi(x, \xi)}{\sigma(x, \xi)} \in  M_{\rho, \Lambda}^{-m}(\overline{\Omega} \times \mathcal{I})$.
	\end{lem}
	Using the above lemma, we obtain the parametrix of the elliptic operators.
	\begin{thm}\label{8}	Let $A: C_{\mathfrak{L}}^\infty(\bar{\Omega})  \to C_{\mathfrak{L}}^\infty(\bar{\Omega}) $ continuous linear operator such that its $\mathfrak{L}$-symbol $\sigma_A$ is  $M$-elliptic. Then  	there exist a symbol $\sigma_{B}\in M_{\rho, 0, \Lambda}^{-m}$ such that
		\begin{eqnarray}\label{left side parametrix}
			BA =I+R
		\end{eqnarray}
		and
		\begin{eqnarray}\label{right side parametrix}
			AB =I+S,
		\end{eqnarray}	
		where the pseudo differential operators $R, S$ are in $\mathrm{Op}_{\mathfrak{L}} M^{-\infty}$.
	\end{thm}
	\begin{proof} Given that $\sigma_{A}$ is an $M$-elliptic symbol of order $m$. Then there exist constants $C>0$ and $N_{0} \in \mathbb{N}$ such that $$|\sigma(x, \xi)|\geq C \Lambda(\xi)^m$$
		for all $(x, \xi) \in \bar{\Omega} \times \mathcal{I}$ for which $|\xi|\geq N_{0} $. The main idea in this proof is to find a sequence of symbols $\sigma_{j} \in M_{\rho,0,\Lambda}^{-m-\rho j}( \overline{\Omega} \times \mathcal{I}), j=0,1,2, \ldots.$ Let us assume that this can be done. Then, by Theorem \ref{4}, there exists a symbol $\sigma_{B} \in M_{\rho,0,\Lambda}^{-m}(\overline{\Omega} \times \mathcal{I})$ such that $\sigma_{B} \sim \sum_{j=0}^{\infty} \sigma_{j}$, and by Theorem \ref{comp}, the symbol $\lambda$ of the product $BA$ is in $M_{\rho,0,\Lambda}^{0}(\overline{\Omega} \times \mathcal{I})$ such that
		\begin{eqnarray}\label{product formula on lambda}
			\lambda-\sum_{|\gamma|<N} \frac{(-i)^{|\gamma|}}{\gamma !}\left(D_x^{(\gamma)} \sigma\right)(\Delta_{(x)}^{\gamma} \sigma_{B}) \in M_{\rho,0,\Lambda}^{-\rho N}(\overline{\Omega} \times \mathcal{I}),
		\end{eqnarray}
		for every positive integer $N$. Also $\sigma_{B} \sim \sum_{j=0}^{\infty} \sigma_{j}$ implies that
		\begin{eqnarray}\label{ass. exp. on tau}
			\sigma_{B}-\sum_{j=0}^{N-1} \sigma_{j} \in M_{\rho,0,\Lambda}^{-m-\rho N}(\overline{\Omega} \times \mathcal{I}),
		\end{eqnarray}
		for every positive integer $N$. Hence, by \eqref{product formula on lambda} and \eqref{ass. exp. on tau},
		\begin{eqnarray}\label{lambda minus forward difference tau class}
			\lambda-\sum_{|\gamma|<N} \frac{(-i)^{|\gamma|}}{\gamma !}\left(D_x^{(\gamma)}\sigma\right) \sum_{j=0}^{N-1}(\Delta_{(x)}^{\gamma} \sigma_{j}) \in M_{\rho,0,\Lambda}^{-\rho N}(\overline{\Omega} \times \mathcal{I}),
		\end{eqnarray}
		for every positive integer $N$. But we can write
		\begin{equation}\label{separation of sum}
			\begin{aligned}[b]
				& \sum_{|\gamma|<N} \frac{(-i)^{|\gamma|}}{\gamma !} \sum_{j=0}^{N-1}(\Delta_{(x)}^{\gamma} \sigma_{j})\left(D_x^{(\gamma)} \sigma\right) \\
				=& \,\sigma_{0} \sigma+\sum_{l=1}^{N-1}\left\{\sigma_{l} \sigma+\sum_{\substack{|\gamma|+j=l \\
						j<l}} \frac{(-i)^{|\gamma|}}{\gamma !}(\Delta_{(x)}^{\gamma} \sigma_{j})\left(D_x^{(\gamma)} \sigma\right)\right\} \\
				+& \sum_{\substack{|\gamma|+j \geq N \\
						|\gamma|<N,\, j<N}} \frac{(-i)^{|\gamma|}}{\gamma !}(\Delta_{(x)}^{\gamma} \sigma_{j})\left(D_x^{(\gamma)} \sigma\right) .
			\end{aligned}
		\end{equation}
		To find a sequence of symbols $\sigma_{j} \in M_{\rho,0,\Lambda}^{-m-\rho j}(\overline{\Omega} \times \mathcal{I}), j=0,1,2, \ldots$, we choose $\psi$ to be any function in $C^{\infty}\left(\mathbb{R}^{k}\right)$ such that $\psi(\xi)=1$, for $|\xi| \geq 2 R$ and $\psi(\xi)=0$, for $|\xi| \leq R$. Define
		\begin{eqnarray}\label{tau_0 defi}
			\sigma_{0}(x, \xi)= \begin{cases}\frac{\psi(\xi)}{\sigma(x, \xi)}, & |\xi|>R, \\ 0, & |\xi| \leq R,\end{cases}, \quad (x,\xi) \in \overline{\Omega} \times \mathcal{I}\text{.}
		\end{eqnarray}
		From Lemma \ref{inverse of symbol class}, it is clear that $\sigma_0 \in M_{\rho,0, \Lambda}^{-m}(\overline{\Omega} \times \mathcal{I})$. Now define, $\sigma_{l}$, for $l \geq 1$, inductively by
		\begin{eqnarray}\label{tau_l defi}
			\sigma_{l}=-\left\{\sum_{\substack{|\gamma|+j=l \\ j<l}} \frac{(-i)^{|\gamma|}}{\gamma !}\left(D_x^{(\gamma)} \sigma\right)(\Delta_{(x)}^{\gamma} \sigma_{j})\right\} \sigma_{0}.
		\end{eqnarray}
		Then it can be shown that $\sigma_{j} \in M_{\rho,0,\Lambda}^{-m-\rho j}( \overline{\Omega} \times \mathcal{I}), j=0,1,2, \ldots$. Now, by \eqref{tau_0 defi}, $\sigma_{0} \sigma=1$, for $|\xi| \geq 2 R$. The second term on the right hand side of \eqref{separation of sum} vanishes for $|\xi| \geq 2 R$ by \eqref{tau_0 defi} and \eqref{tau_l defi}. Also the third term in \eqref{separation of sum} satisfies
		$$
		(\Delta_{(x)}^{\gamma} \sigma_{j})\left(D_x^{(\gamma)} \sigma\right) \in M_{\rho,0,\Lambda}^{-\rho N}(\overline{\Omega} \times \mathcal{I}),
		$$
		whenever $|\gamma|+j \geq N$. Hence, by \eqref{separation of sum},
		\begin{eqnarray}\label{class for third term}
			\sum_{|\gamma|<N} \frac{(-i)^{|\gamma|}}{\gamma !} \sum_{j=0}^{N-1}(\Delta_{(x)}^{\gamma} \sigma_{j})\left(D_x^{(\gamma)} \sigma\right)-1 \in M_{\rho,0,\Lambda}^{-\rho N}(\overline{\Omega} \times \mathcal{I}),
		\end{eqnarray}
		for every positive integer $N$. Thus, by \eqref{lambda minus forward difference tau class} and \eqref{class for third term},
		$$
		\lambda-1 \in M_{\rho,0,\Lambda}^{-\rho N}(\overline{\Omega} \times \mathcal{I}),
		$$
		for every positive integer $N$. Hence, if we pick $R$ to be the pseudo-differential operator with symbol $\lambda-1$, then the proof of \eqref{left side parametrix} is complete.
		By a similar argument, we can find another symbol, $\sigma_{C}\in M_{\rho,0,\Lambda}^{-m}(\overline{\Omega} \times \mathcal{I})$, such that
		\begin{eqnarray}\label{second left side parametrix}
			AC=I+R^{\prime},
		\end{eqnarray}
		where $R^{\prime}$ is a pseudo-differential operator in $\mathrm{Op}_{\mathfrak{L}} M^{-\infty}$. By \eqref{left side parametrix} and \eqref{second left side parametrix},
		$$
		C+R C=B+B R^{\prime} .
		$$
		Since $R C$ and $B R^{\prime}$ are pseudo-differential operators in $\mathrm{Op}_{\mathfrak{L}} M^{-\infty}$, it follows that
		\begin{eqnarray}\label{relation in kappa and tau}
			C=B+R^{\prime \prime},
		\end{eqnarray}
		where
		$$
		R^{\prime \prime}=B R^{\prime}-R C,
		$$
		is another pseudo-differential operator in $\mathrm{Op}_{\mathfrak{L}} M^{-\infty}$. Hence, by \eqref{second left side parametrix} and \eqref{relation in kappa and tau},
		$$
		A B=I+S,
		$$
		where
		$$
		S=R^{\prime}-A R^{\prime \prime} .
		$$
		Since $S$ is a pseudo-differential operator in $\mathrm{Op}_{\mathfrak{L}} M^{-\infty}$, it follows that \eqref{right side parametrix} is proved.
	\end{proof}

	From the duality, the inverse $\mathfrak{L}$-Fourier transform $\mathcal{F}_{\mathfrak{L}}^{-1}:  \mathcal{S}(\mathcal{I}) \to C_{\mathfrak{L}}^{\infty}(\bar{\Omega}) $, the $\mathfrak{L}$-Fourier transform extends uniquely to the mapping $\mathcal{F}_{\mathfrak{L}}: \mathcal{D}_{\mathfrak{L}}'(\Omega)\to \mathcal{S}'(\mathcal{I}) $ by the formula 
	$$\langle \mathcal{F}_{\mathfrak{L}} w, u \rangle =\langle w, \overline{ \mathcal{F}_{\mathfrak{L}^*}^{-1}\bar{u}} \rangle, ~ w\in  \mathcal{D}_{\mathfrak{L}}'(\Omega), u\in \mathcal{S}(\mathcal{I}).$$
	
	For $s\in \mathbb{R}.$ Then $\Lambda(D)^s$ is defined to be the Fourier multiplier given by $$\Lambda(D)^s w= \mathcal{F}_{\mathfrak{L}^*}^{-1} \Lambda^s   \mathcal{F}_{\mathfrak{L}}w, ~w\in  \mathcal{D}_{\mathfrak{L}}'(\Omega).$$ 
	For $s\in \mathbb{R}$, the weighted Sobolev space $\mathcal{H}_{\mathfrak{L}, \Lambda}^{s, 2}(\Omega)$ is defined by
	$$
	\mathcal{H}_{\mathfrak{L}, \Lambda}^{s, 2}(\Omega)=\left\{w \in \mathcal{D}^{\prime}\left(\Omega \right): \Lambda(D)^{s} w \in L^{2}\left(\Omega \right)\right\}.
	$$
	We denote   $ \mathcal{H}_{\mathfrak{L}, \Lambda}^{s, 2}(\Omega)$ as $ \mathcal{H}_{\mathfrak{L}, \Lambda}^{s,2}$.  Moreover $ \mathcal{H}_{\mathfrak{L}, \Lambda}^{s,2}$  is a Banach space
	in which the norm $\|.\|_{\mathcal{H}_{\mathfrak{L}, \Lambda}^{s, 2}}$ is given by
	$$
	\|w\|_{\mathcal{H}_{\mathfrak{L}, \Lambda}^{s, 2}}=\left\|\Lambda(D)^{s} w\right\|_{L^{2}\left(\Omega \right)}, \quad w \in \mathcal{H}_{\mathfrak{L}, \Lambda}^{s, 2}.
	$$
	\begin{rmk}\label{006}
		For $m_1 \leq m_2, \mathcal{H}_{\mathfrak{L}, \Lambda}^{m_2, 2} \subseteq  \mathcal{H}_{\mathfrak{L}, \Lambda}^{m_1, 2}$ and
		$$
		\|f\|_{ \mathcal{H}_{\mathfrak{L}, \Lambda}^{m_1, 2}} \leq \|f\|_{\mathcal{H}_{\mathfrak{L}, \Lambda}^{m_2, 2}}, ~~ f\in \mathcal{H}_{\mathfrak{L}, \Lambda}^{m_1, 2}.$$
	\end{rmk}
	Using the above Sobolev embedding, we have the following boundedness result for $M$-elliptic pseudo-differential operators.
	\begin{prop}\label{5}
		Let the $\mathfrak{L}$-symbol $\sigma \in M_{\rho, 0, \Lambda}^{m}(\bar{\Omega} \times \mathcal{I})$. Then the operator $Op_{\mathfrak{L}} (\sigma )$   extends to a bounded operator from $ \mathcal{H}_{\mathfrak{L}, \Lambda}^{s, 2} $  to $ \mathcal{H}_{\mathfrak{L}, \Lambda}^{s-m, 2}$ for any $s \in \mathbb{R}.$
	\end{prop}
	An immidiate application of above Proposition and  Theorem \ref{8}, we have the following result.
	
	\begin{prop}\label{11}
		Let $T_{\sigma}$ be a pseudo-differential operator with ${\mathfrak{L}}$-symbol $\sigma \in M_{\rho, 0, \Lambda}^{m}, 
		m>0$.  Assume that $\sigma$ is $M$-elliptic. Then there exist positive constants $C$ and $D>0$ such that
		$$
		C\|u\|_{\mathcal{H}_{\mathfrak{L}, \Lambda}^{m, 2}} \leq \|T_\sigma u\|_{\mathcal{H}_{\mathfrak{L}, \Lambda}^{0, 2}}+\|u\|_{\mathcal{H}_{\mathfrak{L}, \Lambda}^{0, 2}} \leq D\|u\|_{\mathcal{H}_{\mathfrak{L}, \Lambda}^{m, 2}}
		$$
		for every $u\in \mathcal{H}_{\mathfrak{L}, \Lambda}^{m, 2}.$
	\end{prop}

	\begin{proof} Since $T_{\sigma}$ is a pseudo-differential operator with $\mathfrak{L}$-symbol $\sigma \in M_{\rho, 0, \Lambda}^{m}  $ be $M$-elliptic. Then from  Theorem \ref{8}, there exists a pseudo-differential operator $T_{\tau}$ with symbol $\tau \in M_{\rho, 0, \Lambda}^{-m}  $ such that
		$$T_\tau T_\sigma=I+R,$$  where the pseudo-differential operators $R\in \mathrm{Op}_{\mathfrak{L}} M^{-\infty}$. Thus
		$$
		T_\tau (T_\sigma u)=(I+R) u,
		$$
		and therefore
		$$
		u=T_\tau (T_\sigma u)-R u
		$$
		for every $u\in \mathcal{H}_{\mathfrak{L}, \Lambda}^{m, 2}.$ Now,  using  Corollary \ref{5}, we have
		\begin{align*}
			\|u\|_{\mathcal{H}_{\mathfrak{L}, \Lambda}^{m, 2}} &= \|T_\tau (T_\sigma u)-R u\|_{\mathcal{H}_{\mathfrak{L}, \Lambda}^{m, 2}} \\
			&\leq \|T_\tau (T_\sigma u)\|_{\mathcal{H}_{\mathfrak{L}, \Lambda}^{m, 2}} +\|Ru\|_{\mathcal{H}_{\mathfrak{L}, \Lambda}^{m, 2}} \\
			&\leq C_1\|T_\sigma u\|_{\mathcal{H}_{\mathfrak{L}, \Lambda}^{0, 2}}+C_2\|u\|_{\mathcal{H}_{\mathfrak{L}, \Lambda}^{0, 2}}\\
			&\leq \max\{ C_1, C_2\}\left(\|T_\sigma u\|_{\mathcal{H}_{\mathfrak{L}, \Lambda}^{0, 2}}+\|u\|_{\mathcal{H}_{\mathfrak{L}, \Lambda}^{0, 2}}\right)
		\end{align*}
		for every $u\in \mathcal{H}_{\mathfrak{L}, \Lambda}^{m, 2}.$	For other side inequality, using   Corollary \ref{5} and  Remark \ref{006}, we have 
		\begin{align*}
			\|T_\sigma u\|_{\mathcal{H}_{\mathfrak{L}, \Lambda}^{0, 2}}+\|u\|_{\mathcal{H}_{\mathfrak{L}, \Lambda}^{0, 2}} &\leq C\|u\|_{\mathcal{H}_{\mathfrak{L}, \Lambda}^{m,2}}+\|u\|_{\mathcal{H}_{\mathfrak{L}, \Lambda}^{m,2}}\\
			&=(C+1) \|u\|_{\mathcal{H}_{\mathfrak{L}, \Lambda}^{m,2}}
		\end{align*}
		for every $u\in \mathcal{H}_{\mathfrak{L}, \Lambda}^{m, 2}.$
	\end{proof}

	\section{Minimal and maximal extension operators}\label{sec4}
	In this section we see that pseudo-differential operators with $\mathfrak{L}$-symbol in $\sigma \in M_{\rho, 0, \Lambda}^{m}, m>0,$  are closable.  Let $\sigma \in M_{\rho, 0, \Lambda}^{m}, m>0$. Then $T_\sigma: C_{\mathfrak{L}}^{\infty}(\bar{\Omega})  \rightarrow C_{\mathfrak{L}}^{\infty}(\bar{\Omega}).$  So, we can consider $T_\sigma$ as a linear operator from $L^{2}(\Omega)$ into $L^{2}(\Omega)$ with dense domain  $C_{\mathfrak{L}}^{\infty}(\bar{\Omega}) $. However, it is not closed. But one can see, it has a closed extension. Regarding this, let us first recall the definition of a closable operator.
	\begin{defn}
		A operator $T$ from a Banach space $X$ to another  Banach space $Y$ is said to be closable if and only if for any sequence $\{  x_k\} $ in $Dom(T)$  such that $x_k \to 0 $ in $X$ and $T(x_k) \to y $ as $k\to \infty $, then $y=0.$
	\end{defn}
	
	\begin{prop}
		Let $\sigma \in M_{\rho, 0, \Lambda}^{m}, m>0$.	The linear operator $T_{\sigma}: L^{2}(\Omega) \rightarrow L^{2}(\Omega)$ is closable with dense domain  $C_{\mathfrak{L}}^{\infty}(\bar{\Omega})$ in $ L^{2}(\Omega).$
	\end{prop}
	\begin{proof}
		Let $\left(\phi_{k}\right)_{k \in \mathbb{N}}$ be a sequence in $C_{\mathfrak{L}}^{\infty}(\bar{\Omega})  $ such that $\phi_{k} \rightarrow 0$ and $T_{\sigma} \phi_{k} \rightarrow f$ for some $f$ in
		$L^{2}(\Omega)$ as $k \rightarrow \infty$.  
		For any $ \psi \in C_{\mathfrak{L}^{*}}^{\infty}(\bar{\Omega}),$ we have 
		\begin{align*}
			\langle T_{\sigma} \phi_{k}, \psi\rangle =\langle \phi_{k}, T_{\sigma}^*  \psi\rangle,
		\end{align*}
		where $T_{\sigma}^* $ is the adjoint of $T_{\sigma}$.	Taking limit as $k \rightarrow \infty$ we have 
		$\langle f, \psi\rangle =0$ for all $\psi \in C_{\mathfrak{L}^{*}}^{\infty}(\bar{\Omega}) 
		$. Since  $ C_{\mathfrak{L}^{*}}^{\infty}(\bar{\Omega}) $ is dense in   $ L^{2}(\Omega)$, it follows that $f=0.$ Hence $T_{\sigma}: L^{2}(\Omega) \rightarrow L^{2}(\Omega)$ is closable.
	\end{proof}
	
	For $\sigma \in M_{\rho, 0, \Lambda}^{m}, m>0$, consider the pseudo-differential operator $T_{\sigma}: L^{2}(\Omega) \rightarrow L^{2}(\Omega)$ with domain $C_{\mathfrak{L}}^{\infty}(\bar{\Omega})$. Then by the previous proposition  $T_{\sigma}$  has a closed extension. Let $T_{\sigma, 0}$ be the minimal operator for $T_{\sigma}$. Let us recall,  the domain ${Dom}\left(T_{\sigma, 0}\right)$ of $T_{\sigma, 0}$ consists of all functions $g \in L^{2}(\Omega)$ for which there exists a sequence $\left(\phi_{k}\right)_{k \in \mathbb{N}}$ in $C_{\mathfrak{L}}^{\infty}(\bar{\Omega})$ such that $\phi_{k} \rightarrow g$ in $L^{2}(\Omega)$ and $T_{\sigma} \phi_{k} \rightarrow f$
	for some $f \in L^{2}(\Omega)$ as $k \rightarrow \infty$.  Moreover  it can be shown that \(f\) does not depend on the choice of \(\left(\phi_{k}\right)_{k \in \mathbb{N}}\) and $T_{\sigma,0} g=f$.	In order to define the maximal operator of $T_{\sigma}$,  we recall the following  definition.
	\begin{defn}\label{8'}
		Let $T_{\sigma, 1}$  be a  linear operator on $L^{2}(\Omega)$. Let $f$ and $g$ be  two function in $L^{2}(\Omega).$  Then we say that $g \in {Dom}\left(T_{\sigma, 1}\right)$ and $T_{\sigma, 1} g=f $
		if and only if
		\begin{align}\label{7'}
			\langle g, T_{\sigma}^{*} \psi\rangle =\langle f, \psi \rangle,
		\end{align}
		for all $\psi \in C_{\mathfrak{L}^{*}}^{\infty}(\bar{\Omega})$, where $T_{\sigma}^{*}$ is the adjoint of $T_{\sigma}$ and $\langle f, g\rangle =\int_{\Omega} f(x)\overline{g(x) }~dx.$
	\end{defn}

	\begin{prop}
		$T_{\sigma, 1}$ is a closed linear operator from $L^{2}(\Omega)$ into $L^{2}(\Omega)$ with domain ${Dom}\left(T_{\sigma, 1}\right)$ containing $C_{\mathfrak{L}}^{\infty}(\bar{\Omega}).$  Moreover $T_{\sigma,1} u=T_{\sigma}u$ in ${Dom}\left(T_{\sigma, 1}\right)$.
	\end{prop}
	\begin{proof}
		From the above definition, clearly $C_{\mathfrak{L}}^{\infty}(\bar{\Omega}) \subset {Dom}\left(T_{\sigma, 1}\right)$. 	Let \(\left\{\phi_{j}\right\}\) be a sequence  in ${Dom}\left(T_{\sigma, 1}\right)$ such that $\phi_{j} \rightarrow u$ and $T_{\sigma, 1} \phi_{j} \rightarrow f$ in \(L^{2}(\Omega)\) for
		some \(u\) and \(f\) in \(L^{2}(\Omega)\) as \(j \rightarrow \infty\). Then, by \eqref{7'}
		$$\langle\phi_{j}, T_{\sigma}^{*} \psi\rangle=\langle T_{\sigma, 1} \phi_{j}, \psi\rangle, ~~\quad \psi \in C_{\mathfrak{L}^{*}}^{\infty}(\bar{\Omega}).$$
		Letting \(j \rightarrow \infty\),  we have $
		\langle u, T_{\sigma}^{*} \psi\rangle=\langle  f, \psi\rangle$  for all $\psi \in C_{\mathfrak{L}^{*}}^{\infty}(\bar{\Omega}).$ By Definition \ref{8'}, $u\in {Dom}\left(T_{\sigma, 1}\right)$   and  $T_{\sigma, 1}u=f.$   Hence $T_{\sigma, 1}$ is a closed operator from $L^{2}(\Omega)$ into $L^{2}(\Omega)$ with domain ${Dom}\left(T_{\sigma, 1}\right)$ containing $C_{\mathfrak{L}}^{\infty}(\bar{\Omega}).$ It can be prove easily that for $u\in {Dom}\left(T_{\sigma, 1}\right)$, $T_{\sigma,1} u=T_{\sigma}u$. It is also easy to prove the linearity.
	\end{proof}
	From the above result we see that $T_{\sigma, 1}$ is a closed operator on $L^{2}(\Omega)$ with domain ${Dom}\left(T_{\sigma, 1}\right)$. Furthermore,   the operator $T_{\sigma, 1}$  is actually  an extension of the operator  $T_{\sigma, 0},$ can be seen  in the following result.
	\begin{prop}\label{13}
		$T_{\sigma, 1}$ is an extension of $T_{\sigma, 0} .$ 
	\end{prop}
	\begin{proof}
		Let \(u \in {Dom}\left(T_{\sigma, 0}\right)\) and \(T_{\sigma, 0} u=f .\) Then, by the definition of  ${Dom}\left(T_{\sigma, 0}\right)$, there exists a sequence \(\left\{\phi_{j}\right\}\) in
		\(C_{\mathfrak{L}}^{\infty}(\bar{\Omega})\) such that \(\phi_{j} \rightarrow u\) and \(T_{\sigma} \phi_{j} \rightarrow f\) in \(L^{2}\left(\Omega\right)\) as \(j \rightarrow \infty \) for some $f\in L^2(\Omega)$. Now $$
		\langle \phi_{j}, T_{\sigma}^{*}\psi \rangle=\langle T_{\sigma} \phi_{j}, \psi\rangle$$
		for all $\psi \in C_{\mathfrak{L}^{*}}^{\infty}(\bar{\Omega})$. Letting \(j \rightarrow \infty\),  we have $
		\langle u, T_{\sigma}^{*} \psi\rangle=\langle  f, \psi\rangle$ for all $ \psi \in C_{\mathfrak{L}^{*}}^{\infty}(\bar{\Omega})$. Therefore  by the Definition \ref{8'}$, u\in {Dom}\left(T_{\sigma, 1}\right)$ and $T_{\sigma,1} u=f$. Hence $T_{\sigma, 1}$ is an extension of $T_{\sigma, 0}.$ 
	\end{proof}
	\begin{prop}
		$C_{\mathfrak{L}^{*}}^{\infty}(\bar{\Omega}) \subseteq  Dom(T_{\sigma,1} ^ { t } )$, where $T_{\sigma,1} ^ { t } $ is the transpose of $T_{\sigma, 1}$.
	\end{prop}
	\begin{prop}
		$T_{\sigma, 1}$ is the largest closed extension
		of $T_{\sigma}$ in the sense that if $B$ is any closed extension of $T_{\sigma}$ such that $C_{\Omega^{*}}^{\infty}(\bar{\Omega}) \subseteq Dom\left(B^{t}\right), $ then
		$T_{\sigma, 1}$ is an extension of $B$. 
	\end{prop}
	\begin{rmk}
		In view of the above, we call $T_{\sigma, 0} $ and $T_{\sigma, 1}$  are the minimal and maximal pseudo differential operator of $T_\sigma.$
	\end{rmk}
The following result will be useful for determining the domain of the minimal and maximal operators of  $T_{\sigma}$.
	\begin{prop}\label{10} 	$C_{\mathfrak{L}}^{\infty}(\bar{\Omega})$ is dense  in $\mathcal{H}_{\mathfrak{L}, \Lambda}^{m, 2}.$	\end{prop}
	\begin{proof}
		Let $g\in \mathcal{H}_{\mathfrak{L}, \Lambda}^{m, 2}.$ Then, from the definition $\Lambda(D)^m g\in L^2(\Omega).$ Density of  $C_{\mathfrak{L}}^{\infty}(\bar{\Omega})$  in $ L^2(\Omega)$ implies there exists a sequence $\{g_k\} \in C_{\mathfrak{L}}^{\infty}(\bar{\Omega})$  such that $ g_k\to\Lambda(D)^m g $ in $L^2(\Omega).$ Let us define $h_k=\Lambda(D)^{-m}g_k.$ Then $h_k\in C_{\mathfrak{L}}^{\infty}(\bar{\Omega})$ for $k=1, 2, \cdots,$ and 
		\begin{align*}
			\|h_k-g\|_{\mathcal{H}_{\mathfrak{L}, \Lambda}^{m, 2}}&=\|\Lambda(D)^{m}(h_k-g)\|_{L^2(\Omega)}\\
			&=\|\Lambda(D)^{m}h_k-\Lambda(D)^{m}g\|_{L^2(\Omega)}\\
			&=\|g_k-\Lambda(D)^{m}g\|_{L^2(\Omega)}\to 0
		\end{align*}
		as $k\to \infty.$ Therefore $C_{\mathfrak{L}}^{\infty}(\bar{\Omega})$ is dense  in $\mathcal{H}_{\mathfrak{L}, \Lambda}^{m, 2}.$
	\end{proof}
	Using Proposition \ref{11}  and Proposition \ref{10}, we have the domain of the minimal operator $T_{\sigma, 0}$ for an $M$-elliptic pseudo-differential operator $T_\sigma$.
	\begin{prop}\label{15}
		Let $T_{\sigma}$ be a pseudo-differential operator with $\mathfrak{L}$-symbol $\sigma \in  M_{\rho, 0, \Lambda}^{m}$ be $M$-elliptic.  Then $Dom( T_{\sigma, 0}) =\mathcal{H}_{\mathfrak{L}, \Lambda}^{m, 2}$.
	\end{prop}
	\begin{proof}
		Let $g \in \mathcal{H}_{\mathfrak{L}, \Lambda}^{m, 2}.$ Then from Proposition \ref{10}, there exists a sequence $\left(g_{k}\right)_{k \in \mathbb{N}}$ in $C_{\mathfrak{L}}^{\infty}(\bar{\Omega})$ such that $$g_{k} \rightarrow g \quad  \text{in $\mathcal{H}_{\mathfrak{L}, \Lambda }^{m, 2}$ }$$
		as $k \rightarrow \infty$. By Proposition \ref{11}, $\left(g_{k}\right)_{k \in \mathbb{N}}$ and $\left(T_{\sigma} g_{k}\right)_{k \in \mathbb{N}}$ are Cauchy sequences in $L^{2}(\Omega)$. Hence $g_{k} \rightarrow g$	and $T_{\sigma} g_{k} \rightarrow f$ for some $f \in L^{2}(\Omega)$ as $k \rightarrow \infty$. By the definition of $T_{\sigma, 0}$, this implies that $g \in {Dom}\left(T_{\sigma, 0}\right)$ and $T_{\sigma, 0} g=f$. Therefore $$ \mathcal{H}_{\mathfrak{L}, \Lambda}^{m, 2}\subseteq {Dom}\left(T_{\sigma, 0}\right).$$
		
		On the other side,  suppose that  $g \in {Dom}\left(T_{\sigma, 0}\right).$ Then from the definition of $ {Dom}\left(T_{\sigma, 0}\right)$, there exists a sequence $\left(g_{k}\right)_{k \in \mathbb{N}}$ in $C_{\mathfrak{L}}^{\infty}(\bar{\Omega})$ such that
		\begin{align}\label{12}
			g_{k} \rightarrow g \quad \text{in $L^{2}(\Omega) $ and }~T_{\sigma} g_{k} \rightarrow f
		\end{align}
		for some $f \in L^{2}(\Omega) .$ Therefore, by the left inequality of Proposition \ref{11}, $\left(g_{k}\right)_{k \in \mathbb{N}}$ is
		a Cauchy sequence in $\mathcal{H}_{\mathfrak{L}, \Lambda}^{m, 2}$. Since $\mathcal{H}_{\mathfrak{L}, \Lambda}^{m, 2}$ is complete, there exists a function $h \in \mathcal{H}_{\mathfrak{L}, \Lambda}^{m, 2}$ such that
		$g_{k} \rightarrow h$ in $\mathcal{H}_{\mathfrak{L}, \Lambda}^{m, 2} .$  From Proposition \ref{006}, this implies that $g_{k} \rightarrow h$ in $L^{2}(\Omega)$. Therefore $ h=g \in \mathcal{H}_{\mathfrak{L}, \Lambda}^{m, 2} $ and  completes the proof.
	\end{proof}
	
	The following theorem is the main result  of this section which states that,  under   the  $M$-ellipticity assumption of the  symbols $\sigma$,   the minimal extension $T_{\sigma, 0}$ and maximal extension $T_{\sigma, 1}$  of a pseudo-differential operator $T_\sigma$ coincide on $ L^{2}(\Omega)$.
	\begin{thm}
		Let $T_{\sigma}$ be a pseudo-differential operator with $\mathfrak{L}$-symbol $\sigma \in M_{\rho, 0, \Lambda}^{m}, m>0$ be $M$-elliptic.  Then $T_{\sigma, 0}=T_{\sigma, 1}$.
	\end{thm}
	
	\begin{proof}
		Since  $T_{\sigma, 1}$ is a closed extension of $T_{\sigma, 0}$,  by Proposition \ref{15}, it is sufficient to prove that ${Dom}\left(T_{\sigma, 1}\right) \subseteq \mathcal{H}_{\mathfrak{L}, \Lambda}^{m, 2} $. Let $g\in {Dom}\left(T_{\sigma, 1}\right) .$ Since $\sigma$ is an $M$-elliptic, there exists $\tau \in M_{\rho, 0, \Lambda}^{-m}$
		such that
		$$
		g=T_{\tau} T_{\sigma} g-Rg,
		$$
		where the pseudo-differential operator $R \in Op_{\mathfrak{L}} M^{-\infty}$. Since $T_{\sigma} g=T_{\sigma,1} g \in L^{2}(\Omega)$, it follows from 
		Corollary  \ref{5} that $T_{\tau} T_{\sigma} g \in \mathcal{H}_{\mathfrak{L}, \Lambda}^{m, 2}$. Since $g\in L^2(\Omega)$ and $R$ has symbol in $ M_{\rho, 0, \Lambda}^{-m}$, again using Corollary  \ref{5}, we can conclude that $Rg \in  \mathcal{H}_{\mathfrak{L}, \Lambda}^{m, 2}$. Therefore $g \in  \mathcal{H}_{\mathfrak{L}, \Lambda}^{m, 2}.$
	\end{proof}
	Now, we will discuss the resolvent set of  $M$-elliptic pseudo-differential operator  whose symbol is  independent of $ x \in \bar{\Omega}.$ For this, we will recall some basic definitions and results.
	Let $A: X \rightarrow X$ be a  bounded  linear operator  from a complex Banach
	space $X$ to itself with dense domain $D(A)$. We define the spectrum $\Sigma(A)$ and essential spectrum $\Sigma_{e}(A)$ of $A$ by
	$$
	\Sigma(A)=\mathbb{C} \backslash \rho(A),
	$$
	where $\rho(A)$, the resolvent set of $A$ is defined by $$\rho(A)=\{\lambda \in \mathbb{C}: A-\lambda I \text { is one to one and onto}\}.$$
	We start with the following theorem, which is a straightforward application of closed graph theorem.
	\begin{thm}\label{007}
		Let $A$ be a closed linear operator from a complex Banach space $X$ onto $X$ with dense domain $Dom(X)$. Then a complex number $\lambda$ is in the resolvent set $\rho(A)$ of
		$A$ if and only if the range $R(A-\lambda I)$ of $A-\lambda I$ is dense in $X$ and there exists a positive
		constant $C$ such that
		$$
		\|x\| \leq C\|(A-\lambda I) x\|, ~x \in {Dom}(A).
		$$
	\end{thm}
	The following theorem gives us the resolvent set of the  $M$-elliptic pseudo-differential operators.
	
	\begin{thm}\label{008}
		Let $\sigma \in M_{\rho, 0, \Lambda}^{m}, m>0$ be $M$-elliptic.  Suppose that $\sigma$ is independent of $x$ in $\bar{\Omega}$. If $\lambda$ is a complex number such that
		$$
		\sigma(\xi) \neq \lambda ,~~ \xi \in \mathcal{I},
		$$
		then $\lambda \in \rho\left(T_{\sigma, 0}\right)$.
	\end{thm}
	\begin{proof}
		Let us  consider the function $\tau$ on $\mathcal{I}$  by $\tau(\xi)=\frac{1}{\sigma(\xi)-\lambda}, ~\xi \in \mathcal{I}$. Then $\tau \in M_{\rho, 0, \Lambda}^{-m}$ and from Corollary \ref{5}, $T_{\tau}: \mathcal{H}_{\mathfrak{L}, \Lambda}^{s, 2} \rightarrow  \mathcal{H}_{\mathfrak{L}, \Lambda}^{s+m, 2}$ is a bounded linear operator for any
		$-\infty<s<\infty$. By  Remark \ref{5}  and Proposition \ref{11}, there exists a positive constant $C$ such that
		$$
		\|\phi\|_{\mathcal{H}_{\mathfrak{L}, \Lambda}^{0, 2}}=\left\|T_{\tau} T_{\sigma-\lambda} \phi\right\|_{\mathcal{H}_{\mathfrak{L}, \Lambda}^{0, 2}} \leq C\left \|T_{\sigma-\lambda} \phi\right\|_{\mathcal{H}_{\mathfrak{L}, \Lambda}^{-m, 2}} \leq C\left\|\left(T_{\sigma}-\lambda \right)  \phi \right\|_{\mathcal{H}_{\mathfrak{L}, \Lambda}^{0, 2}}
		$$
		for all $\phi$ in $C_{\mathfrak{L}}^{\infty}(\bar{\Omega})$. By a density argument, we can conclude that
		$$
		\|u\|_{\mathcal{H}_{\mathfrak{L}, \Lambda}^{0, 2}} \leq C\left\|\left(T_{{\sigma}, 0}-\lambda I\right) u\right\|_{\mathcal{H}_{\mathfrak{L}, \Lambda}^{0, 2}},~ \forall  u \in {Dom}\left(T_{\sigma, 0}\right).
		$$

		The proof will complete if we can show that   the range $R\left(T_{\sigma, 0}-\lambda I\right)$ is dense in $L^{2}\left({\Omega}\right)$. For that it is enought to show  $C_{\mathfrak{L}}^{\infty}(\bar{\Omega}) \subseteq R\left(T_{\sigma, 0}-\lambda I\right).$  Let $\psi \in C_{\mathfrak{L}}^{\infty}(\bar{\Omega})$. Then the function $\frac{\hat{\psi}}{\sigma-\lambda}$ is in $C_{\mathfrak{L}}^{\infty}(\bar{\Omega})$.  Let us define a 
		function $\phi$ on $\Omega$  by
		$$
		\phi(x)=\sum_{\xi\in \mathcal{I}} \frac{\hat{\psi}(\xi)}{\sigma(\xi)-\lambda} u_{\xi}(x), \quad x \in \Omega,
		$$
		is in $C_{\mathfrak{L}}^{\infty}(\bar{\Omega})$. The $\mathfrak{L}$-Fourier transform of $\phi $ is given by  $\hat{\phi}(\xi)=\frac{\hat{\psi}(\xi)}{\sigma(\xi)-\lambda}.$
		Now \begin{align*}
			\left(T_{\sigma}-\lambda I\right) \phi(x)&=\sum_{\xi\in \mathcal{I}} (\sigma(\xi)-\lambda)\hat{\phi}(\xi) u_{\xi}(x)\\
			&=\sum_{\xi\in \mathcal{I}}\hat{\psi}(\xi) u_{\xi}(x)=\psi(x).
		\end{align*}
		Therefore $C_{\mathfrak{L}}^{\infty}(\bar{\Omega}) \subseteq R\left(T_{\sigma, 0}-\lambda I\right)$ and hence $R\left(T_{\sigma, 0}-\lambda I\right)$ is dense in $L^2(\Omega).$ From Theorem \ref{007},  we have $\lambda \in \rho\left(T_{\sigma, 0}\right)$.
	\end{proof}

	\section{Gohberg's lemma, Compact operators and Riesz operators}\label{sec5}
	This section  is devoted to  study the compact and Riesz $M$-elliptic pseudo-differential operators on $\bar{\Omega} \times \mathcal{I}.$
	First, we provide a necessary and sufficient condition for compactness   of pseudo-differential operators $T_{\sigma}$ on $L^{2}(\Omega)$ with $\mathfrak{L}$-symbols and $\mathfrak{L}^{*}$-symbols in the Hörmander class $M_{\rho, 0,\Lambda}^{0}(\bar{\Omega} \times \mathcal{I})$ and $\tilde{M}_{\rho, 0,\Lambda}^{0}(\bar{\Omega} \times \mathcal{I})$, respectively, where  $0 <  \rho \leq 1.$ Moreover, we prove that the same condition ensures that $T_{\sigma}$ extends to a Riesz operator in $L^{p}(\Omega), 1<p<\infty.$
	Let $ \mathfrak{K}\left(L^{p}(\Omega)\right), 1\leq p \leq \infty,$ be the ideal of compact operators of $ \mathcal{L}(L^{p}(\Omega)).$ $\mathfrak{R}\left(L^p(\Omega)\right), 1\leq p \leq \infty,$ is the largest ideal in $ \mathcal{L}(L^{p}(\Omega))$ which only consists of Riesz operators. For a  detailed study on  algebraic considerations of Riesz opeators on $L^{p}(\Omega)$,   we refer to  \cite{MR&JPVR}. Also, we assume that $\Omega$ has finite measure 1 throughout this section.
	We start this section by the following result, which   is known as Gohberg's lemma in the literature.
	\begin{lem}\textbf{(Gohberg's lemma)}\label{gohberg lemma}
		Let $1<p<\infty.$ Let $T_\sigma$ be a pseudo-differential operator with $\mathfrak{L}$-symbol $\sigma \in M_{\rho, 0,\Lambda}^0(\bar{\Omega} \times \mathcal{I}), 0<\rho \leq 1$. Then
		$$
		\left\|T_\sigma-K\right\|_{\mathcal{L}\left(L^2(\Omega)\right)} \geq d_\sigma
		$$
		for all compact operators $K \in \mathfrak{K}\left(L^2(\Omega)\right)$, where
		$$
		d_\sigma:=\limsup _{\langle\xi\rangle \rightarrow \infty}\left\{\sup _{x \in \Omega}|\sigma(x, \xi)|\right\} .
		$$
	\end{lem}
	The proof of the above  Gohberg's lemma will go similar lines of   \cite[Lemma 3.1]{MR&JPVR} with slightly appropriate  modification. Using Gohberg's lemma, we have the following compactness  result  of pseudo-differential operators $T_{\sigma}$ on $L^{2}(\Omega)$. 
	\begin{thm}\label{compact operator}
		Let $T_{\sigma}$ be a pseudo-differential operator with $\mathfrak{L}$-symbol, $\sigma,$ in the Hörmander class $M_{\rho, 0,\Lambda}^{0}(\overline{\Omega} \times \mathcal{I}),0 <  \rho \leq 1.$ Then $T_{\sigma}$ extends to a compact operator in $L^{2}(\Omega)$ if and only if
		
		$$
		d_{\sigma}:=\limsup _{\langle\xi\rangle \rightarrow \infty}\left\{\sup _{x \in \Omega}|\sigma(x, \xi)|\right\}=0
		$$
	\end{thm}
	\begin{proof}
		Assume that $d_{\sigma}=0$. Then for $f \in C_{\mathfrak{L}}^{\infty}(\bar{\Omega}),$ we have
		$$
		\begin{aligned}
			\left(T_{\sigma} f\right)(x) & =\sum_{\xi \in \mathcal{I}} \sigma(x, \xi) \widehat{f}(\xi) u_{\xi}(x) \\
			& =\sum_{\xi \in \mathcal{I}}\left(\sum_{\eta \in \mathcal{I}} \hat{\sigma}(\eta, \xi) u_{\eta}(x)\right) \widehat{f}(\xi) u_{\xi}(x) \\
			& =\sum_{\eta \in \mathcal{I}} u_{\eta}(x)\left(\sum_{\xi \in \mathcal{I}} \widehat{\sigma}(\eta, \xi) \widehat{f}(\xi) u_{\xi}(x)\right) \\
			& =\sum_{\eta \in \mathcal{I}} u_{\eta}(x)\left(T_{\widehat{\sigma}_{\eta}} f\right)(x)\\
			&=\sum_{\eta \in \mathcal{I}}\left(A_{\eta} T_{\widehat{\sigma}_{\eta}} f\right)(x), \quad x \in \bar{\Omega},
		\end{aligned}
		$$
		where $\widehat{\sigma}_{\eta}(\xi):=\widehat{\sigma}(\eta, \xi),$ $\left(A_{\eta} f\right)(x):=u_{\eta}(x) f(x)$ is a multiplication operator, and the change in the order of summation is justified by Fubini-Tonelli's theorem. Since  
		$$
		\left\|A_{\eta} f\right\|_{L^{2}(\Omega)} \leq\left\|u_{\eta}\right\|_{L^{\infty}(\Omega)}\|f\|_{L^{2}(\Omega)} \leq C \langle\eta\rangle^{\mu_{0}}\|f\|_{L^{2}(\Omega)}, ~~\text{for each}~ \eta \in \mathcal{I}~~(\text{Using the initial assumption}),
		$$
		we can conclude that  $A_{\eta} \in \mathcal{L}\left(L^{2}(\Omega)\right).$ 
		Now, for each $\eta \in \mathcal{I},$ we get
		\begin{equation}\label{limit of sigmna_n}
			\begin{aligned}[b]
				\lim _{\langle\xi\rangle \rightarrow \infty}|\widehat{\sigma}(\eta, \xi)| & =\lim _{\langle\xi\rangle \rightarrow \infty}\left|\int_{\Omega} \sigma(x, \xi) \overline{u_{\eta}(x)} d x\right| \\
				& \leq \lim _{\langle\xi\rangle \rightarrow \infty}\|\sigma(\cdot, \xi)\|_{L^{2}(\Omega)}  \\
				& \leq \lim _{\langle\xi\rangle \rightarrow \infty}\left\{\sup _{x \in \Omega}|\sigma(x, \xi)|\right\} \\
				& \leq \limsup _{\langle\xi\rangle \rightarrow \infty}\left\{\sup _{x \in \Omega}|\sigma(x, \xi)|\right\}=0.
			\end{aligned}
		\end{equation}
		Also,  for each $\eta \in \mathcal{I}$, the operator $T_{\widehat{\sigma}_{\eta}}$ is a Fourier multiplier. From \cite{MR&JPVR}, it is well-known that a pseudo-differential operator with symbol $\sigma(\xi)$ is a compact operator in $L^{2}(\Omega)$ if and only if
		$$
		\lim _{\langle\xi\rangle \rightarrow \infty}|\sigma(\xi)|=0, 
		$$ 
		from the estimate \eqref{limit of sigmna_n}, 	 for each  $\eta \in \mathcal{I}$,   $T_{\widehat{\sigma}_{\eta}}$ is a compact operator.   
		Moreover, since    $A_{\eta} \in \mathcal{L}\left(L^{2}(\Omega)\right)$,  
		$A_{\eta} T_{\widehat{\sigma}_{\eta}}$ is compact for each $ \eta \in \mathcal{I}.$  This implies that for all $N \in \mathbb{N}$, the operator
		$$
		\sum_{\langle\eta\rangle \leq N} A_{\eta} T_{\widehat{\sigma}_{\eta}}
		$$
		is also compact,  since the set of compact operators $\mathfrak{K}\left(L^{2}(\Omega)\right)$ form a closed two sided ideal in $\mathcal{L}\left(L^{2}(\Omega)\right)$ in the operator norm topology. Hence, if the series
		$$
		\sum_{\eta \in \mathcal{I}} A_{\eta} T_{\widehat{\sigma}_{\eta}}
		$$
		converges in the operator norm topology, then
		$$
		T_{\sigma}=\lim _{N \rightarrow \infty} \sum_{\langle\eta\rangle \leq N} A_{\eta} T_{\widehat{\sigma}_{\eta}}, 
		$$
		the limit of a sequence of compact operators, 	is also a compact operator.  Now it remains to show that   the series
		$$
		\sum_{\eta \in \mathcal{I}} A_{\eta} T_{\widehat{\sigma}_{\eta}}
		$$
		converges in the operator norm topology. 
		Since $\sigma \in M_{\rho, 0,\Lambda}^{0}(\bar{\Omega} \times \mathcal{I}),$ then
		\begin{align}\label{garding convergence}
			\sum_{\eta \in \mathcal{I}}\left\|A_{\eta} T_{\widehat{\sigma}_{\eta}}\right\|_{\mathcal{L}\left(L^{2}(\Omega)\right)} \leq \sum_{\eta \in \mathcal{I}} C \langle\eta\rangle^{\mu_{0}}\left\|T_{\widehat{\sigma}_{\eta}}\right\|_{\mathcal{L}\left(L^{2}(\Omega)\right)} \leq C^{\prime} \frac{m_{2}}{m_{1}} \sum_{\eta \in \mathcal{I}} \Lambda(\eta) \sup _{\xi \in \mathcal{I}}|\widehat{\sigma}(\eta, \xi)|,
		\end{align}
		where $m_{1}, m_{2}$ are as in Remark \ref{estimate of f and f hat}. Note that 
		$$
		\sum_{\eta \in \mathcal{I}}\Lambda(\eta) \sup _{\xi \in \mathcal{I}}|\widehat{\sigma}(\eta, \xi)|=\sum_{\eta \in \mathcal{I}}(\Lambda(\eta))^{-s_{0}} \sup_{\xi \in \mathcal{I}} (\Lambda(\eta))^{1+s_{0}}|\widehat{\sigma}(\eta, \xi)| \leq \sum_{\eta \in \mathcal{I}}(\Lambda(\eta))^{-s_{0}} \sup _{\xi \in \mathcal{I}}\|\sigma(\cdot, \xi)\|_{\mathcal{H}_{\mathfrak{L},\Lambda}^{1+s_{0},2}}<\infty,
		$$
		thus 	 the above sum \eqref{garding convergence} converges. 
		So, $T_{\sigma}$ is a compact operator. 
		
		To show the converse part, assume that $d_{\sigma} \neq 0$. We need only to show that $T_{\sigma}$ is not compact on $L^{2}(\Omega)$. Suppose that $T_{\sigma}$ is compact. If we set $T_{\sigma}=K$ in Lemma \ref{gohberg lemma}, then it contradicts our assumption that $d_{\sigma} \neq 0$. This completes the proof of the theorem.
	\end{proof} 
	Analogously,  one can prove the following theorem without making any substantial modification  to the proof of Theorem \ref{compact operator}.
	\begin{thm}
		Let $T_{\tau}$ be a pseudo-differential operator with $\mathfrak{L}^{*}$-symbol $\tau \in \tilde{M}_{\rho,0,\Lambda}^{0}(\bar{\Omega} \times \mathcal{I})$. Then $T_{\tau}$ extends to a compact operator in $L^{2}(\Omega)$ if and only if
		$$
		d_{\tau}:=\limsup _{\langle\xi\rangle \rightarrow \infty}\left\{\sup _{x \in \Omega}|\tau(x, \xi)|\right\}=0
		$$
	\end{thm}
	Now, we will recall some basic definition and results related to Riesz operators. Riesz operators are the generalization of compact operators which have many equivalent definitions present in literature. We refer to \cite{caradus,pietsch,ruston,west} for some fundamental results on Riesz operator. Here, we will use the definition given in \cite{pietsch}.
	\begin{defn}
		Let $X$ be a Banach space. We say that $T \in \mathcal{L}(X)$ is a Riesz operator if for every $\varepsilon>0$, there exist an exponent $s$ and points $v_{1}, \ldots, v_{l} \in X$ which depend on $\varepsilon$ such that
		$$
		T^{s}\left(B_{X}\right) \subseteq \bigcup_{j=1}^{l} v_{j}+\varepsilon^{s} B_{X},
		$$
		where $B_{X}$ is the unit ball of $X$.
	\end{defn}
	The following lemma lists several prerequisites and sufficient requirements for becoming a Riesz operator for any arbitrary Banach space, which can be found in \cite{caradus,pietsch}.
	\begin{lem}\label{riesz basic result}
		Let $X$ be a Banach space and suppose $T \in \mathcal{L}(X)$. The following statements are equivalent:
		\begin{enumerate}[(i)]
			\item $T$ is a Riesz operator.
			\item $T^{m}$ is Riesz for some exponent $m$.
			
			\item $\lambda T$ is iterative compact for all $\lambda \in \mathbb{C}$.
			
			\item Non-zero points in $ {Spec} (T)$ are isolated eigenvalues with finite-dimensional associated eigenspace and zero as the only cluster point. If $X$ is infinite dimensional then $0 \in$ $ {Spec}(T)$.
			
			\item $T-\lambda I$ is a Fredholm operator for all $\lambda \in \mathbb{C}$.
		\end{enumerate}
	\end{lem}
	In the next result we provide a  necessary and sufficient condition to ensures that $T_{\sigma}$ extends to a Riesz operator in $L^{p}(\Omega).$
	\begin{thm}
		Let $T_{\sigma}$ be a pseudo-differential operator with symbol $\sigma \in M_{\rho,0,\Lambda}^{0}(\overline{ \Omega}\times \mathcal{I})$ and $1<p<\infty$. Then $T_{\sigma}$ is a Riesz operator on $L^{p}(\Omega)$ if and only if
		
		$$
		d_{\sigma}^{\prime}:=\lim_{\langle \xi \rangle \rightarrow \infty }\left\{\sup _{x \in  \Omega}|\sigma(x, \xi)|\right\}=0 .
		$$
	\end{thm}
	\begin{proof}
		First, let us suppose  that $d_{\sigma}^{\prime}=0.$   Then, similarly as in Theorem \ref{compact operator}, for all $ f \in  C_{\mathfrak{L}}^{\infty} \left(\overline{\Omega} \right) $, we can write
		$$
		\left(T_\sigma f\right)(x)=\sum_{\eta \in \mathcal{I}} u_\eta(x)\left(\sum_{\xi \in \mathcal{I}} \widehat{\sigma}(\eta, \xi) \widehat{f}(\xi) u_{\xi}(x)\right) =\sum_{\eta \in \mathcal{I}} u_\eta(x) (T_{\widehat{\sigma}_\eta}f)(x),
		$$
		where 
		$$
		\left(T_{\widehat{\sigma}_\eta} f\right)(x)=\sum_{\xi \in \mathcal{I}} \hat{\sigma}\left(\eta, {\xi}\right) \hat{f}(\xi) u_{\xi}(x).
		$$
		Since 
		$$
		\lim _{\langle \xi \rangle \rightarrow \infty}\left|\widehat{\sigma}\left(\eta,  {\xi}\right)\right| \leq \lim_{\langle \xi \rangle \rightarrow \infty}\left\{\sup _{x \in \Omega}\left|\sigma\left(x,  {\xi}\right)\right|\right\}=0 ,
		$$
		from   \cite[Proposition 3.8]{MR&JPVR},  each operator $T_{\widehat{\sigma}_\eta}$  that belongs to the operator ideal $\mathfrak{R}$, 	is a Riesz operator.
		As $\mathfrak{R}\left(L^p(\Omega)\right)$ is a closed operator ideal,   using the same approximation argument  as in Theorem \ref{compact operator} and  \cite[Theorem 3.2]{JPVR}, one can observe that $T_\sigma$ belongs to $\mathfrak{R}\left(L^p(\Omega)\right)$. 
		
		To prove converse, first we show that if a pseudo-differential operator with symbol in $M_{\rho,0,\Lambda}^{0}(\overline{ \Omega}\times \mathcal{I})$ is a Riesz operator then $d_{\sigma^{n}}^{\prime}=0$ for some $n \in \mathbb{N}$, where $\sigma^{n}$ denotes the symbol of $\left(T_{\sigma}\right)^{n}$. Once we have that,  we further  show that $d_{\sigma^{n}}^{\prime}=0$ if and only if $d_{\sigma}^{\prime}=0$.
		
		Now suppose that $T_{\sigma}$ is a Riesz operator and $d_{\sigma^{n}}^{\prime}>0$ for all $n$. By using (iii) in Lemma \ref{riesz basic result}, we get that, $\left(T_{\sigma}\right)^{k}$ is a compact operator  for some $k \in \mathbb{N}$. Then, by Gohberg's lemma, we have
		$$
		0<d_{\sigma^{k}}^{\prime} \leq d_{\sigma^{k}} \leq\left\|\left(T_{\sigma}\right)^{k}-K\right\|_{\mathcal{L}\left(L^{p}( \Omega)\right)}
		$$
		for all compact operators $K$. If we set $K=\left(T_{\sigma}\right)^{k}$, the above chain of inequalities contradicts our assumption that $d_{\sigma^{k}}^{\prime} \neq 0$.  This implies that $d_{\sigma^{k}}^{\prime} = 0$.
		Now its remains to show that $d_{\sigma^{n}}^{\prime}=0$ if and only if $d_{\sigma}^{\prime}=0$.   By Theorem \ref{comp}, taking the first two terms in the asymptotic expansion of $T_{\sigma^{n}},$ we have
		$$
		\sigma^{n}(x, \xi )=(\sigma(x, \xi))\left(\sigma^{n-1}(x, \xi )\right)+r_{n}(x, \xi)
		$$
		with $r_{n}(x, \xi ) \in M_{\rho,0,\Lambda}^{-\rho}(\overline{ \Omega}\times \mathcal{I})$. Repeating this process, we obtain
		$$
		\begin{aligned}
			\sigma^{n}(x, \xi ) & =(\sigma(x, \xi ))^{2}\left(\sigma^{n-2}(x, \xi )\right)+(\sigma(x, \xi ))\left(r_{n-1}(x, \xi )\right)+r_{n}(x, \xi ) \\
			& ~\vdots \\
			& =(\sigma(x, \xi ))^{n}+\sum_{j=0}^{n-1}(\sigma(x, \xi ))^{j}\left(r_{n-j}(x, \xi )\right),
		\end{aligned}
		$$
		with each $r_{j}(x, \xi ) \in M_{\rho,0,\Lambda}^{-\rho}(\overline{ \Omega}\times \mathcal{I})$. From this, we have
		$$
		(\sigma(x, \xi ))^{n}=\sigma^{n}(x, \xi )-\sum_{j=0}^{n-1}(\sigma(x, \xi ))^{j}\left(r_{n-j}(x, \xi )\right),
		$$
		which implies that
		$$
		|\sigma(x, \xi )|^{n} \leq\left|\sigma^{n}(x, \xi )\right|+\sum_{j=0}^{n-1}|\sigma(x, \xi )|^{j}\left|r_{n-j}(x, \xi )\right|.
		$$
		Hence
		$$
		\sup _{x \in\Omega}|\sigma(x, \xi )|^{n} \leq \sup _{x \in\Omega}\left|\sigma^{n}(x, \xi )\right|+\sum_{j=0}^{n-1} C^{j} \sup _{x \in\Omega}\left|r_{n-j}(x, \xi )\right|,
		$$
		where $\sigma(x,\xi)$ is uniformly bounded by a positive number $C$ because of the fact that $\sigma \in M_{\rho,0,\Lambda}^{0}(\overline{ \Omega}\times \mathcal{I})$. Therefore, 
		$$
		\lim _{\langle \xi \rangle \rightarrow \infty} \sup _{x \in \Omega}|\sigma(x, \xi )|^{n} \leq \lim _{\langle \xi \rangle \rightarrow \infty} \sup _{x \in \Omega}\left|\sigma^{n}(x, \xi )\right|+\sum_{j=0}^{n-1} \lim _{\langle \xi \rangle \rightarrow \infty} C^{j} \sup _{x \in\Omega}\left|r_{n-j}(x, \xi )\right| .
		$$
		On the other hand, $d_{\sigma}^{\prime}=0 \implies d_{\sigma^{n}}^{\prime}=0$ follows directly by definition of $d_{\sigma}^{\prime}.$ This completes  the proof of the theorem.
	\end{proof}
	\section{The G\r{a}rding inequality}\label{sec6}
 This section is devoted investigate G\"arding's inequality for  $\mathfrak{L}$-elliptic   pseudo-differential operators associated with  the symbol   from    $M_{\rho, 0,\Lambda}^{0}(\bar{\Omega} \times \mathcal{I})$ in the setting of  non-harmonic. As an application of G\r{a}rding's inequality,  we provide sufficient conditions for the existence and uniqueness of strong solutions in $L^{2}(\overline{\Omega})$ for the pseudo-differential operator $T_{\sigma}$.

	Let $m>0$ and $0 <\rho \leqslant 1$. Let $I$ be an interval (finite or infinite) in $\mathbb{R},$ and $\Lambda^{\prime}=\left\{\gamma(t): t \in I\right\}$ be an analytic curve in the complex plane $\mathbb{C}$. For simplicity, we assume that, if $I$ is a finite interval then $\Lambda^{\prime}$ is a closed curve, and if $I$ is an infinite interval then $\Lambda^{\prime}$ is homotopy equivalent to the line $\Lambda^{\prime}_{i \mathbb{R}}:=\{i y:-\infty<y<\infty\}$. Let $\sigma=\sigma(x, \xi) \in M_{\rho,0, \Lambda}^m(\overline{\Omega} \times \mathcal{I})$. We also assume   that $R_\lambda(x, \xi)^{-1}:=\sigma(x, \xi)-\lambda \neq 0$ for every $(x, \xi) \in \overline{\Omega} \times \mathcal{I}$, and $\lambda \in \Lambda^{\prime}$. We say that $\sigma$ is parameter $\mathfrak{L}$-elliptic with respect to $\Lambda^{\prime}$, if
	$$
	\sup _{\lambda \in \Lambda^{\prime}} \sup _{(x, \xi) \in \overline{\Omega} \times \mathcal{I}}\left|\left(|\lambda|^{\frac{1}{m}}+\Lambda(\xi)\right)^m R_\lambda(x, \xi)\right|<\infty.
	$$
	
	In order to  prove the G\"arding's inequality for the global pseudo-differential calculus, the following preliminaries are essential.
	\begin{lem}\label{sobolev estimates}
		Let us assume that $s \geqslant t \geqslant 0$ or that $s, t<0$. Then, for every $\varepsilon>0$, there exists $C_{\varepsilon}>0$ such that
		$$
		\|u\|_{\mathcal{H}_{\mathfrak{L}, \Lambda}^{t, 2}} \leqslant \varepsilon \|u\|_{\mathcal{H}_{\mathfrak{L}, \Lambda}^{s, 2}}+C_{\varepsilon}\|u\|_{\mathcal{H}_{\mathfrak{L}, \Lambda}^{0, 2}},
		$$
		holds true for every $u \in C_{\mathfrak{L}}^{\infty}(\overline{\Omega})$.
	\end{lem}
	\begin{proof}
		Let $\varepsilon>0$. Then, there exists $C_{\varepsilon}>0$ such that
		$$
		(\Lambda(\xi))^{2 t}-\varepsilon(\Lambda(\xi))^{2 s} \leqslant C_{\varepsilon},
		$$
		uniformly in $\xi \in \mathcal{I}$. So, we get
		$$
		\begin{aligned}
			\|u\|_{\mathcal{H}_{\mathfrak{L}, \Lambda}^{t, 2}}& =\sum_{\xi \in \mathcal{I}}(\Lambda(\xi))^{2 t} \widehat{u}(\xi) \overline{\widehat{u}_{*}(\xi)} \leqslant \sum_{\xi \in \mathcal{I}}\left(\varepsilon(\Lambda(\xi))^{2 s}+C_{\varepsilon}\right)\widehat{u}(\xi) \overline{\widehat{u}_{*}(\xi)} \\
			& =\epsilon \|u\|_{\mathcal{H}_{\mathfrak{L}, \Lambda}^{s, 2}}+C_{\varepsilon}\|u\|_{\mathcal{H}_{\mathfrak{L}, \Lambda}^{0, 2}}.
		\end{aligned}
		$$
	\end{proof}
	\begin{prop}\label{exponential symbol}
		Let $m>0$ and $0 <\rho \leqslant 1$. Let $\sigma \in M_{\rho,0,\Lambda}^{m}(\overline{\Omega} \times \mathcal{I})$ be an $\mathfrak{L}$-elliptic symbol with the assumption that $\sigma>0$. Then $\sigma$ is parameter $\mathfrak{L}$-elliptic with respect to $\mathbb{R}_{-}:=\{z=x+i 0: x<0\} \subset \mathbb{C}$. Furthermore, for any number $s \geq 0$,  if we define a symbol  
		$$
		\widehat{B}_s(x, \xi) \equiv \sigma(x, \xi)^{s}:=\exp (s \log (\sigma(x, \xi))),\quad (x, \xi) \in \overline{\Omega} \times \mathcal{I},
		$$
		then  $\widehat{B}_{s}(x, \xi) \in M_{\rho,0,\Lambda}^{m.s}(\overline{\Omega} \times \mathcal{I})$.
	\end{prop}
	\begin{proof}
		Since $\sigma \in M_{\rho,0,\Lambda}^{m}(\overline{\Omega} \times \mathcal{I})$ is $\mathfrak{L}$-elliptic, we obtained that
		$$
		\sup _{(x, \xi) \in \overline{\Omega} \times \mathcal{I}}\left|(\Lambda(\xi))^{-m} \sigma(x, \xi)\right|<\infty, \quad \sup _{(x, \xi) \in \overline{\Omega} \times \mathcal{I}}\left|(\Lambda(\xi))^{m} \sigma(x, \xi)^{-1}\right|<\infty.
		$$
		Now, for every $\lambda \in \mathbb{R}_{-},$ we have
		$$
		\begin{aligned}
			&\left|\left(|\lambda|^{\frac{1}{m}}+\Lambda(\xi)\right)^{m}(\sigma(x, \xi)-\lambda)^{-1}\right| \\
			=&\left|\left(|\lambda|^{\frac{1}{m}}+\Lambda(\xi)\right)^{m}\left((\Lambda(\xi))^{m}-\lambda\right)^{-1}\left((\Lambda(\xi))^{m}-\lambda\right)(\sigma(x, \xi)-\lambda)^{-1}\right| \\
			\lesssim&\left|\left(|\lambda|^{\frac{1}{m}}+\Lambda(\xi)\right)^{m}\left((\Lambda(\xi))^{m}-\lambda\right)^{-1}\right |\left|\left((\Lambda(\xi))^{m}-\lambda\right)(\sigma(x, \xi)-\lambda)^{-1}\right| \\
			\lesssim&\left|\left(|\lambda|^{\frac{1}{m}}+\Lambda(\xi)\right)^{m}\left((\Lambda(\xi))^{m}-\lambda\right)^{-1}\right| .
		\end{aligned}
		$$
		By using the compactness of $[0,1 / 2],$ we deduce that
		$$
		\begin{aligned}
			\sup _{0 \leqslant |\lambda| \leqslant 1 / 2}\left|\left(|\lambda|^{\frac{1}{m}}+\Lambda(\xi)\right)^{m}\left((\Lambda(\xi))^{m}-\lambda\right)^{-1}\right| & \lesssim \sup _{0 \leqslant |\lambda| \leqslant 1 / 2}\left|(\Lambda(\xi))^{m}\left((\Lambda(\xi))^{m}-\lambda\right)^{-1}\right| \\
			& =\sup _{0 \leqslant |\lambda| \leqslant 1 / 2}\left|\left(1-\lambda(\Lambda(\xi))^{-m}\right)^{-1}\right| \\
			& \lesssim 1.
		\end{aligned}
		$$
		On the other hand, we have
		$$
		\begin{aligned}
			& \sup _{|\lambda| \geqslant 1 / 2}\left|\left(|\lambda|^{\frac{1}{m}}+\Lambda(\xi)\right)^{m}\left((\Lambda(\xi))^{m}-\lambda\right)^{-1}\right| \\
			&=\sup _{1/2 \leqslant |\lambda| \leqslant (\Lambda(\xi))^{m}}\left|\left(|\lambda|^{\frac{1}{m}}+\Lambda(\xi)\right)^{m}\left((\Lambda(\xi))^{m}-\lambda\right)^{-1}\right| + \sup _{|\lambda| \geqslant (\Lambda(\xi))^{m}}\left|\left(|\lambda|^{\frac{1}{m}}+\Lambda(\xi)\right)^{m}\left((\Lambda(\xi))^{m}-\lambda\right)^{-1}\right| \\
			& =\sup _{1/2 \leqslant |\lambda| \leqslant (\Lambda(\xi))^{m}}\left|\left(|\lambda|^{\frac{1}{m}}(\Lambda(\xi))^{-1}+1\right)^{m}\left(1-(\Lambda(\xi))^{-m} \lambda\right)^{-1}\right|\\
			&\quad + \sup _{|\lambda| \geqslant (\Lambda(\xi))^{m}}\left|\left(|\lambda|^{\frac{1}{m}}+\Lambda(\xi)\right)^{m}\left((\Lambda(\xi))^{m}-\lambda\right)^{-1}\right| \\
			& =\sup _{1/2 \leqslant |\lambda| \leqslant (\Lambda(\xi))^{m}}\left|\left((\Lambda(\xi))^{-1}+|\lambda|^{-\frac{1}{m}}\right)^{m}| \lambda|\left(1-(\Lambda(\xi))^{-m} \lambda\right)^{-1}\right|\\
			&\quad + \sup _{|\lambda| \geqslant (\Lambda(\xi))^{m}}\left|\left(|\lambda|^{\frac{1}{m}}+\Lambda(\xi)\right)^{m}|\lambda|^{-1}\left((\Lambda(\xi))^{m}|\lambda|^{-1}+1\right)^{-1}\right|\\
			& \lesssim 1.
		\end{aligned}
		$$
		Hence, we obtained that $\sigma$ is parameter $\mathfrak{L}$-elliptic with respect to $\mathbb{R}_{-}$. Since $s \geqslant 0,$ then there exists $k \in \mathbb{N}$ such that $s-k<0$. By using the spectral calculus of matrices, we deduce that $\sigma(x, \xi)^{s-k} \in$ $M_{\rho,0,\Lambda}^{m.(s-k)}(\overline{\Omega} \times \mathcal{I})$. So, from the global pseudo-differential calculus we conclude that
		$$
		\widehat{B}_{s}(x, \xi)= \sigma(x, \xi)^{s}=\sigma(x, \xi)^{s-k} \sigma(x, \xi)^{k} \in M_{\rho,0,\Lambda}^{m.s}(\overline{\Omega} \times \mathcal{I}).
		$$
		This completes the proof.
	\end{proof}
	An immediate consequences of the above proposition is the following one.
	\begin{Cor}\label{exponentional 1/2-symbol}
		Let $m>0$, and $0 <\rho \leqslant 1$. Let $\sigma \in M_{\rho,0,\Lambda}^{m}(\overline{\Omega} \times \mathcal{I})$ be a $\mathfrak{L}$-elliptic symbol and let us assume that $\sigma>0$ Then $\widehat{B}(x, \xi) \equiv \sigma(x, \xi)^{\frac{1}{2}}:=$ $\exp \left(\frac{1}{2} \log (\sigma(x, \xi))\right) \in M_{\rho,0,\Lambda}^{\frac{m}{2}}(\overline{\Omega} \times \mathcal{I})$.
	\end{Cor}
	The following lower bound result is the main result of this section, which   is known as  G\"arding's inequality in the literature. 
	\begin{thm}\label{garding inequality}
		Let $m>0,$ and $0 <\rho \leqslant 1$. Let $\sigma(x, D): C_{\mathfrak{L}}^{\infty}(\overline{\Omega}) \rightarrow \mathcal{D}_{\mathfrak{L}}^{\prime}(\Omega)$ be an operator with symbol $\sigma \in M_{\rho,0, \Lambda}^{m}(\overline{\Omega} \times \mathcal{I}).$ For $ \sigma \in M_{\rho,0,\Lambda}^{m}(\overline{\Omega} \times \mathcal{I}),$ let us assume that
		$$
		A(x, \xi):=\frac{1}{2}(\sigma(x, \xi)+\overline{\sigma(x, \xi)}),\quad (x, \xi) \in \overline{\Omega} \times \mathcal{I},
		$$
		satisfies
		\begin{equation}\label{estimate of symbol in garding}
			\left|(\Lambda(\xi))^{m} A(x, \xi)^{-1}\right| \leqslant C_{0},
		\end{equation}
		for some $C_{0}>0.$ Then, there exist $C_{1}, C_{2}>0$, such that
		$$
		\operatorname{Re}(\sigma(x, D) u, u) \geqslant C_{1}\|u\|_{\mathcal{H}_{\mathfrak{L}, \Lambda}^{\frac{m}{2}, 2}}-C_{2}\|u\|_{\mathcal{H}_{\mathfrak{L}, \Lambda}^{0, 2}}^{2}
		$$
		holds true for every $u \in C_{\mathfrak{L}}^{\infty}(\overline{\Omega})$.
	\end{thm}
	\begin{proof}
		From the given condition
		$$
		\left|(\Lambda(\xi))^{m} A(x, \xi)^{-1}\right| \leqslant C_{0},
		$$
		we get
		$$
		(\Lambda(\xi))^{-m} A(x, \xi)\geqslant \frac{1}{C_{0}}.
		$$
		This implies that
		$$
		A(x, \xi) \geqslant \frac{1}{C_{0}}(\Lambda(\xi))^{m},
		$$
		and for $C_{1} \in\left(0, \frac{1}{C_{0}}\right),$ we have that
		$$
		A(x, \xi)-C_{1}(\Lambda(\xi))^{m} \geqslant\left(\frac{1}{C_{0}}-C_{1}\right)(\Lambda(\xi))^{m}>0 .
		$$
		Now, by using Corollary \ref{exponentional 1/2-symbol}, we have
		$$
		q(x, \xi):=\left(A(x, \xi)-C_{1}(\Lambda(\xi))^{m}\right)^{\frac{1}{2}} \in M_{\rho,0,\Lambda}^{\frac{m}{2}}(\overline{\Omega} \times \mathcal{I}).
		$$
		Using the asymptotic expansion for the product in Theorem \ref{comp}, we obtain
		$$
		q(x, D) q(x, D)^{*}=A(x, D)-C_{1} \mathrm{Op}\left((\Lambda(\xi))^{m}\right)+r(x, D), \quad r(x, \xi) \in M_{\rho,0,\Lambda}^{m-\rho}(\overline{\Omega} \times \mathcal{I}).
		$$
		Now, let us assume that $u \in C_{\mathfrak{L}}^{\infty}(\overline{\Omega})$. Denoting $\mathcal{M}_{m}:=\mathrm{Op}\left((\Lambda(\xi))^{m}\right)$, $m \in \mathbb{R}$, we have
		$$
		\begin{aligned}
			\operatorname{Re}(\sigma(x, D) u, u) & =\frac{1}{2}\left(\left(\sigma(x, D)+\sigma(x, D)^{*}\right) u, u\right)=(A(x, D) u, u) \\
			& =C_{1}\left(\mathcal{M}_{m} u, u\right)+\left(q(x, D) q(x, D)^{*} u, u\right)-(r(x, D) u, u) \\
			& =C_{1}\left(\mathcal{M}_{m} u, u\right)+\left(q(x, D)^{*} u, q(x, D)^{*} u\right)-(r(x, D) u, u) \\
			& \geqslant C_{1}\|u\|_{\mathcal{H}_{\mathfrak{L}, \Lambda}^{\frac{m}{2}, 2}}-(r(x, D) u, u) \\
			& =C_{1}\|u\|_{\mathcal{H}_{\mathfrak{L}, \Lambda}^{\frac{m}{2}, 2}}-\left(\mathcal{M}_{-\frac{m-\rho}{2}} r(x, D) u, \mathcal{M}_{\frac{m-\rho}{2}} u\right).
		\end{aligned}
		$$
		By using Corollary \ref{5}, one can observe that
		$$
		\begin{aligned}
			\left(\mathcal{M}_{-\frac{m-\rho}{2}} r(x, D) u, \mathcal{M}_{\frac{m-\rho}{2}} u\right) & \leqslant\left\|\mathcal{M}_{-\frac{m-\rho}{2}} r(x, D) u\right\|_{\mathcal{H}_{\mathfrak{L}, \Lambda}^{0, 2}}\|u\|_{\mathcal{H}_{\mathfrak{L}, \Lambda}^{\frac{m-\rho}{2}, 2}} \\
			& =\|r(x, D) u\|_{\mathcal{H}_{\mathfrak{L}, \Lambda}^{-\frac{m-\rho}{2}, 2}}\|u\|_{\mathcal{H}_{\mathfrak{L}, \Lambda}^{\frac{m-\rho}{2}, 2}} \\
			& \leqslant C_{1}\|u\|_{\mathcal{H}_{\mathfrak{L}, \Lambda}^{\frac{m-\rho}{2}, 2}}\|u\|_{\mathcal{H}_{\mathfrak{L}, \Lambda}^{\frac{m-\rho}{2}, 2}}.
		\end{aligned}
		$$
		Hence, by applying Lemma \ref{sobolev estimates}, we can find positive constants $C_{1}^{\prime}, C_{2}^{\prime}>0$, such that
		$$
		\operatorname{Re}(\sigma(x, D) u, u) \geqslant C_{1}^{\prime}\|u\|_{\mathcal{H}_{\mathfrak{L}, \Lambda}^{\frac{m}{2}, 2}}-C_{2}^{\prime}\|u\|_{\mathcal{H}_{\mathfrak{L}, \Lambda}^{0, 2}}^{2}
		$$
		holds true for every $u \in C_{\mathfrak{L}}^{\infty}(\overline{\Omega})$. This completes the proof of the theorem.
	\end{proof}
	
	Now, we present an application of G\r{a}rding's inequality for the   class $M_{\rho, \Lambda}^{m}\left(\Omega \times \mathcal{I}\right).$ Let  $T_{\sigma, 0} $ and $T_{\sigma, 1}$  are the minimal and maximal pseudo differential operator of $T_\sigma$ on $L^2(\Omega)$. First, we give the definition of strong solutions.
	\begin{defn}
		Let $\sigma \in M_{\rho,\Lambda}^{2m}(\Omega \times \mathcal{I}), m>0$, and let $f \in L^{p}\left(\Omega\right), 1<p<\infty$. A function $u \in L^{p}\left(\Omega\right)$ is   a strong solution of the equation $T_{\sigma} u=f$ if $u \in \mathcal{D}\left(T_{\sigma, 0}\right)$ and $T_{\sigma, 0} u=f$.
	\end{defn}
	Proceeding similarly as in Theorem 18.2 of \cite{MWbook}, we have the following result related to  the strong solution of the equation $T_{\sigma} u=f$ on $\Omega$.
	\begin{lem}\label{Lemma in garding application}
		Let $\sigma \in M_{\rho,\Lambda}^{2m}(\Omega \times \mathcal{I}), m>0$, be an $\mathfrak{L}$-elliptic symbol such that 
		\begin{align}\label{garding app}
			\operatorname{Re}(T_{\sigma} u, u) \geq C\|u\|_{\mathcal{H}_{\mathfrak{L}, \Lambda}^{m, 2}}, \quad u\in \mathcal{H}_{\mathfrak{L}, \Lambda}^{m, 2},
		\end{align}
		for some  positive  constant $C$.  Then for every function $f$ in $L^{2}\left(\Omega\right)$, the pseudo-differential equation $T_{\sigma} u=f$   has a unique strong solution $u$ in $L^{2}\left(\Omega\right)$.
	\end{lem} 
	In the next result, we provide sufficient conditions for the existence and uniqueness of strong solutions in $L^{2}(\Omega)$ for the pseudo-differential operator $T_{\sigma}$ with symbol satisfying \eqref{estimate of symbol in garding}.
	\begin{thm} Let $\sigma \in  M_{\rho,\Lambda}^{2m}(\Omega \times \mathcal{I}), m>0$, be such that it satisfies the condition \eqref{estimate of symbol in garding}. Then  for all $f$ in $L^{2}\left(\Omega\right)$  there exists a real number $\lambda_{0}$ such that  for   all $\lambda \geq \lambda_{0}$, the pseudo-differential equation $\left(T_{\sigma}+\lambda I\right) u=f$ on $\Omega$ has a unique strong solution $u$ in $L^{2}\left(\Omega\right)$, where $I$ is the identity operator on $L^{2}\left(\Omega\right)$.
	\end{thm}
	
	\begin{proof}
		From  Gårding's inequality, there exist constants $A>0$ and $\lambda_{0}$ such that  
		$$
		\operatorname{Re}\left(T_{\sigma} u, u\right) \geq A\|u\|_{\mathcal{H}_{\mathfrak{L}, \Lambda}^{m, 2}}-\lambda_{0}\|u\|_{\mathcal{H}_{\mathfrak{L}, \Lambda}^{0, 2}}^{2}, \quad u \in \mathcal{H}_{\mathfrak{L}, \Lambda}^{m, 2}.
		$$
		Now  for any $\lambda \geq \lambda_{0}$, we get 
		\begin{align*}
			\operatorname{Re}\left(\left(T_{\sigma}+\lambda I\right) u, u\right) &=\operatorname{Re}\left (T_{\sigma}   u, u\right)+\left(\lambda-\lambda_{0}\right)\|u\|_{\mathcal{H}_{\mathfrak{L}, \Lambda}^{0, 2}}^{2} +\lambda_{0}\|u\|_{\mathcal{H}_{\mathfrak{L}, \Lambda}^{0, 2}}^{2}  \\
			&\geq A \|u\|_{\mathcal{H}_{\mathfrak{L}, \Lambda}^{m, 2}}+\left(\lambda-\lambda_{0}\right)\|u\|_{\mathcal{H}_{\mathfrak{L}, \Lambda}^{0, 2}}^{2} \geq A\|u\|_{\mathcal{H}_{\mathfrak{L}, \Lambda}^{m, 2}}.
		\end{align*}
		This shows that  $	T_{\sigma}$  satisfies the condition (\ref{garding app}).  Thus by  Lemma \ref{Lemma  in garding application},   the pseudo-differential equation $\left(T_{\sigma}+\lambda I\right) u=f$  has a unique strong solution $u$ in $L^{2}\left(\Omega\right)$ and this completes the proof.
	\end{proof}  
	
	\section{Acknowledgement}
	The first was supported by Core Research Grant(RP03890G), Science and Engineering Research Board (SERB), DST, India.   The second author was supported by the FWO Odysseus 1 grant G.0H94.18N: Analysis and Partial Differential Equations, the Methusalem program of the Ghent University Special Research Fund (BOF) (grant number 01M01021), and by FWO Senior Research Grant G011522N.  The fourth author was supported by Indian Institute of Technology Delhi Institute Post-doctoral Fellowship.


\begin{thebibliography}{aaa}

 \bibitem{Agmon} S.  Agmon, The coerciveness problem for integro-differential forms, J. Analyse Math. 6, 183--223, 1958. \vspace{6pt}
		
		\bibitem{AM} M. Alimohammady and M. K. Kalleji, Fredholmness property of $M$-elliptic pseudo-differential operator under change variable in its symbol, J. Pseudo-Differ. Oper. Appl., 4(3), 371–392, 2013.\vspace{6pt}
		
		\bibitem{bari}  N. K. Bari, Biorthogonal systems and bases in Hilbert space, Moskov. Gos. Univ. Ucenye Zapiski Matematika 148(4), 69–107, 1951.\vspace{6pt}
		
		\bibitem{beals} R. Beals, A general class of pseudo-differential operators, Duke Math. J., 57, 1-42, 1975.\vspace{6pt}
		
		\bibitem{caradus} S. R. Caradus and W. E. Pfaffenberger, Calkin algebras and algebras of operators on Banach spaces, 1st edn. Lecture Notes in Pure and Applied Mathematics, Marcel Dekker, New York, 1974.\vspace{6pt}
		
		\bibitem{DC&MR} D. Cardona, V. Kumar,  M. Ruzhansky, and N.  Tokmagambetov, Global functional calculus, lower/upper bounds and evolution equations on manifolds with boundary, ( to appear) Adv. Oper. Theory, 2023.  arXiv:2101.02519, 2021.\vspace{6pt}
		
		\bibitem{vish} A.  Dasgupta and V.  Kumar, Ellipticity and Fredholmness of pseudo-differential operators on $\ell^2(\mathbb{Z}^n)$, Proc. Amer. Math. Soc., 150(7), 2849–2860,  2022. \vspace{6pt}
		
		\bibitem{AD&LM} A. Dasgupta and  L. Mohan, $M$-Ellipticity of Fredholm pseudo-differential operators on and G\r{a}rding’s inequality, J. Geom. Anal., 33(3), 97, 2023. \vspace{6pt}
		
		\bibitem{dasgupta} A.  Dasgupta and M.  Ruzhansky, The Gohberg lemma, compactness, and essential spectrum of operators on compact Lie groups, J. Anal. Math., 128, 179-190, 2016.\vspace{6pt}
		
		\bibitem{DMT17} J. Delgado, M. Ruzhansky, and N. Tokmagambetov, Schatten classes, nuclearity and nonharmonic analysis on compact manifolds with boundary, J. Math. Pures Appl., 107(6), 758-783, 2017.\vspace{6pt}
		
		\bibitem{apara} A. Dasgupta and M. Ruzhansky, Eigenfunction expansions of ultradifferentiable functions and ultradistributions, 
		Trans. Amer. Math. Soc., 36(12), 8481–8498, 2016. \vspace{6pt}
		
		\bibitem{dasgupta1} A. Dasgupta and M. W. Wong, Spectral invariance of SG pseudo-differential operators on $L^p(\mathbb{R}^n)$, In Pseudo-differential operators: complex analysis and partial differential equations, volume 205 of Oper. Theory Adv. Appl., pages 51–57. Birkh\"{a}user Verlag, Basel, 2010. \vspace{6pt}
		
		\bibitem{delgado1} J. Delgado and M. Ruzhansky, Fourier multipliers, symbols and nuclearity on compact manifolds, J. Math. Anal. Appl.,
		135, 757–800, 2018.\vspace{6pt}
		
		\bibitem{delgado2} J. Delgado and M. Ruzhansky, Kernel and symbol criteria for Schatten classes and r-nuclearity on compact manifolds, C. R. Math. Acad. Sci. Paris, 352(10),  779–784, 2014.\vspace{6pt} 
		
		\bibitem{feff} C. Fefferman, $L^p$-bounds for pseudo-differential operators,  Israel J. Math., 14, 413-417, 1973.\vspace{6pt}
		
		\bibitem{fis1}  V. Fischer and M. Ruzhansky, Quantization on nilpotent Lie groups, Progress in Mathematics, Vol. 314,   Birkh\"auser, Basel, 2016.\vspace{6pt}
		
		\bibitem{LG} L. G\r{a}rding, Dirichlet’s problem for linear elliptic partial differential equations, Math. Scand., 1, 55-72, 1953.\vspace{6pt}
		
		
		\bibitem{GM} G. Garello and A. Morando, A class of $L^p$ bounded pseudo-differential operators, In Progress in analysis, Vol. I, II (Berlin, 2001), {World Sci. Publ., River Edge, NJ}, pages 689-696, 2003.\vspace{6pt}
		
		\bibitem{mor05}   G. Garello and A. Morando, $L^p$-bounded pseudo-differential opearors and regularity for multi-quasi elliptic equations,  Integr. Equ. Oper. Theory, 51, 501-517, 2005.\vspace{6pt}
		
		\bibitem{mor16}   G. Garello and  A. Morando, $m$-microlocal elliptic pseudo-differential operators acting on $L_{loc}^p (\Omega)$, Math. Nachr., 289, 1820-1837, 2016.\vspace{6pt}
		
		\bibitem{Goh} 	I. C. Gohberg, On the theory of multi-dimensional singular integral equations, Soviet Math. Dokl., 1, 960–963, 1960.\vspace{6pt}
		
		\bibitem{VVG} V. V. Gru\v{s}in, Pseudo-differential operators in $\mathbb{R}^n$ with bounded symbols, Funkcional. Anal. i Prilo\v{z}en, 4(3), 37–50, 1970.\vspace{6pt}
		
		\bibitem{Hor}  L. H\"ormander, The analysis of linear partial differential operators, Vol. III. Pseudo-differential operators, Springer-Verlag, Berlin Heidelberg NewYork Tokyo, 1985.\vspace{6pt}
		
		\bibitem{inoue2016non} H. Inoue   and M. Takakura, Non-self-adjoint Hamiltonians defined by generalized Riesz bases, J. Math. Phys., 57(8), 083505, 2016.\vspace{6pt}
		
		\bibitem{kal} M. K. Kalleji, Essential spectrum of $M$-hypoelliptic pseudo-differential on the torus, J. Pseudo-Differ. Oper. Appl., 6, 439-459, 2016.\vspace{6pt}
		
		\bibitem{KT15} B. Kanguzhin, N. Tokmagambetov, and K. Tulenov, Pseudo-differential operators generated by a non-local boundary value problem, Complex Var. Elliptic Equ., 60(1), 107-117, 2015. \vspace{6pt}
		
		\bibitem{Niren} J. J. Kohn and L. Nirenberg, An algebra of pseudo-differential operators, Comm. Pure Appl. Math., 18, 269-305, 1965.\vspace{6pt}
		
		\bibitem{SSM&VK} V. Kumar and S. S.  Mondal,  Symbolic calculus and $M$-ellipticity
		of pseudo-differential operators on $\mathbb{Z}^n$, arXiv:2111.10224, 2021.\vspace{6pt}
		
		\bibitem{SM&MW} S. Molahajloo and M. W. Wong, Ellipticity, Fredholmness and spectral invariance of pseudo-differential operators on $\mathbb{S}^1$,  J. Pseudo-Differ. Oper. Appl., 1, 183-205, 2010.\vspace{6pt}
		
		\bibitem{mostafazadeh2010pseudo} A. Mostafazadeh,  Pseudo-Hermitian representation of quantum mechanics, Int. J. Geom. Methods Mod. Phys., 7(07), 1191-1306, 2010.\vspace{6pt}
		
		\bibitem{paley1934fourier} R. Paley and N. Wiener, Fourier transforms in the complex domain, Colloq,  Publ. Amer. Math. Soc., 19,  1934.\vspace{6pt}
		
		\bibitem{pietsch} A. Pietsch,  Eigenvalues and S-Numbers, Cambridge studies in advanced mathematics, vol. 13. Cambridge university press, Cambridge, 1987.\vspace{6pt}
		
		\bibitem{ruston} A. F. Ruston, Operators with a fredholm theory, J. Lond. Math. Soc., 3, 318–326, 1964.\vspace{6pt}
		
		\bibitem{Ruz16} M. Ruzhansky and N. Tokmagambetov, Non-harmonic analysis of boundary value problems, Int. Math. Res. Not. IMRN, (12), 3548-3615, 2016. \vspace{6pt}
		
		\bibitem{ruzhansky2017nonharmonic} M. Ruzhansky and N. Tokmagambetov,
		Non-harmonic analysis of boundary value problems without WZ condition, Math. Model. Nat. Phenom., 12(1),  115-140, 2017. \vspace{6pt}
		
		\bibitem{ruzhansky2017very} M. Ruzhansky and N. Tokmagambetov, Very weak solutions of wave equation for Landau Hamiltonian with irregular electromagnetic field, Lett. Math. Phys., 107, 591-618, 2017.\vspace{6pt}
		
		\bibitem{ruzhansky2017wave} M. Ruzhansky and N. Tokmagambetov, Wave equation for operators with discrete spectrum and irregular propagation speed, Arch. Ration. Mech. Anal., 226(3), 1161-1207, 2017.\vspace{6pt}
		
		\bibitem{MR&VT book} M. Ruzhansky and V. Turunen, Pseudo-differential Operators and Symmetries: Background Analysis and Advanced Topics, In volume 2 of Pseudo-differential operators. Theory and Applications, {Birkh\"{a}user-Verlag, Basel}, 2010.\vspace{6pt}
		
		\bibitem{wirth}   M. Ruzhansky, V. Turunen, and J. Wirth, H\"ormander class of pseudo-differential operatorson compact Lie groups and global hypoellipticity, J. Fourier Anal. Appl., 20(3), 476–499, 2014.\vspace{6pt}
		
		\bibitem{MR&JPVR} M. Ruzhansky and J. P. Velasquez-Rodriguez, Non-harmonic Gohberg’s lemma, Gershgorin theory and heat equation on manifolds with boundary, Math. Nachr., 294(9), 1783–1820, 2021.\vspace{6pt}
		
		\bibitem{MR&JW}	M. Ruzhansky and J. Wirth, Global functional calculus for operators on compact Lie groups, J. Funct. Anal., 267, 772–798, 2014.\vspace{6pt}

  
		
		\bibitem{saki1}  A. M. Sedletskii, Nonharmonic analysis, J. Math. Sci. (N.Y.), 116(5), 3551–3619, 2003. \vspace{6pt}
		
		\bibitem{seeley1} R. T. Seeley, Integro-differential operators on vector bundles, Trans. Amer. Math. Soc., 117, 167–204, 1965.\vspace{6pt}
		
		\bibitem{seeley2} R. T. Seeley, Eigenfunction expansions of analytic functions, Proc. Amer. Math. Soc., 21, 734–738, 1969.\vspace{6pt}

  \bibitem{Smith}K. T. Smith,
Inequalities for formally positive integro-differential forms,
Bull. Amer. Math. Soc. 67, 368–370, 1961.
 
 \bibitem{Sche} M. Schechter,
Integral inequalities for partial differential operators and functions satisfying
general boundary conditions,
Comm. Pure Appl. Math. 12, 37--66, 1959.


 \bibitem{Sche2} M. Schechter, Coerciveness in $L^p$, Trans. Amer. Math. Soc. 107, 10--29, 1963.
		
		\bibitem{taylor81}   M. E. Taylor, Pseudo-differential operators,  Princton university press, Princton, 1981.\vspace{6pt}
		
		\bibitem{JPVR} J. P. Velasquez-Rodriguez, On some spectral properties of pseudo-differential operators on $\mathbb{T}$, J. Fourier Anal. Appl., 25, 2703-2732, 2019.\vspace{6pt}
		
		\bibitem{west} T. T. West, Riesz operators in banach spaces, Proc. Lond. Math. Soc., 16(1), 131–140, 1996.\vspace{6pt}
		
		
		\bibitem{MWspectral} M. W. Wong, Spectral theory of pseudo-differential operators, Adv. in Appl. Math., 15(4), 437-451, 1994.\vspace{6pt}
		
		\bibitem{MWMelliptic} M. W. Wong. $M$-elliptic pseudo-differential operators on $L^p(\mathbb{R}^n)$, Math. Nachr., 279(3), 319–326, 2006.\vspace{6pt}
		
		\bibitem{MWbook} M. W. Wong,  {An introduction to pseudo-differential operators}, volume 6 of Series on Analysis, Applications and Computation,  {World Scientific Publishing Co. Pte. Ltd., Hackensack, NJ}, 2014.\vspace{6pt}
		
		
		
	\end{thebibliography}
\end{document}